\documentclass[10pt,letterpaper]{article}
\usepackage[textwidth=5.5in,textheight=9in,top=1in,centering,footskip=30pt]{geometry}
\usepackage[T1]{fontenc}
\usepackage{times}

\usepackage{amsthm}

\usepackage{amsmath,amssymb,mathtools,bm,microtype}
\usepackage{graphicx,xcolor,booktabs,longtable,array,calc,etoolbox}
\usepackage{placeins,float,pifont,needspace,pdflscape}
\usepackage{caption,subcaption}
\usepackage[authoryear,round]{natbib}
\bibpunct[, ]{(}{)}{,}{a}{}{,}

\def\BIBand{and}
\usepackage{xurl}
\usepackage[hidelinks]{hyperref}
\usepackage{bookmark}
\hypersetup{pdftitle={Financial Connectedness: Heavy Tails and Response Attribution},
pdfauthor={Yuzhao Zhu; Jiahui Shao; Dan Li; Dongqi Liu; Xiaoxiao Hu}}
\newtheorem{theorem}{Theorem}
\newtheorem{lemma}{Lemma}
\newtheorem{proposition}{Proposition}
\newtheorem{corollary}{Corollary}
\theoremstyle{definition}

\theoremstyle{remark}
\newtheorem{remark}{Remark}

\newcommand{\spacingset}[1]{}

\makeatletter
\renewcommand{\normalsize}{\@setfontsize\normalsize{10}{11}%
  \abovedisplayskip 7pt plus 2pt minus 2pt
  \belowdisplayskip \abovedisplayskip
  \abovedisplayshortskip 0pt plus 2pt
  \belowdisplayshortskip 4pt plus 2pt minus 2pt}
\normalsize
\renewcommand{\section}{\@startsection{section}{1}{\z@}%
  {-2.2ex plus -0.5ex minus -0.2ex}{1.1ex plus 0.2ex}{\large\bfseries}}
\renewcommand{\subsection}{\@startsection{subsection}{2}{\z@}%
  {-1.8ex plus -0.5ex minus -0.2ex}{0.8ex plus 0.2ex}{\normalsize\bfseries}}
\renewcommand{\subsubsection}{\@startsection{subsubsection}{3}{\z@}%
  {-1.5ex plus -0.4ex minus -0.2ex}{0.6ex plus 0.2ex}{\normalsize\bfseries}}
\def\ps@preprintfirst{\def\@oddhead{}\let\@evenhead\@oddhead
  \def\@oddfoot{\footnotesize Preprint.\hfil}\let\@evenfoot\@oddfoot}
\renewcommand{\maketitle}{\par\begingroup
  \thispagestyle{preprintfirst}\vspace*{-3.8pt}\hrule height 4pt\vskip 14pt
  {\centering\fontsize{17.28}{20}\selectfont\bfseries\@title\par}
  \vskip 14pt\hrule height 1pt\vskip 24pt
  {\centering\@author\par}\vskip 18pt
  \endgroup}
\makeatother
\pretocmd{\subsection}{\Needspace{6\baselineskip}}{}{}
\renewenvironment{abstract}{\par\begingroup
  \begin{center}\large\bfseries Abstract\end{center}\vspace{-3pt}
  \begin{list}{}{\setlength{\leftmargin}{0.5in}\setlength{\rightmargin}{0.5in}%
    \setlength{\topsep}{0pt}\setlength{\parsep}{0pt}}\item\relax\normalsize
}{\end{list}\endgroup\vspace{4pt}}
\newenvironment{keywords}{\par\small\noindent\textit{Keywords:} }{\par\normalsize\medskip}
\newcommand{\preprintparagraph}[1]{\par\Needspace{4\baselineskip}\addvspace{4pt}\noindent\textbf{#1}\quad}
\makeatletter
\renewcommand{\paragraph}{\@ifstar{\preprintparagraph}{\preprintparagraph}}
\makeatother
\title{Financial Connectedness: Heavy Tails and Response Attribution}
\author{\begin{tabular}{@{}c@{\hspace{0.25in}}c@{\hspace{0.25in}}c@{}}
\large Yuzhao Zhu\textsuperscript{1} & \large Jiahui Shao\textsuperscript{1} & \large Dan Li\textsuperscript{1}\\[3pt]
\multicolumn{3}{c}{\large Dongqi Liu\textsuperscript{1}\hspace{0.3in}Xiaoxiao Hu\textsuperscript{1,*}}
\end{tabular}\\[10pt]
\textsuperscript{1}School of Mathematics and Statistics, Lanzhou University\\
Lanzhou, Gansu 730000, China\\[4pt]
\texttt{\{zhuyzh2025,shjiahui2025,ldan2024\}@lzu.edu.cn}\\[2pt]
\texttt{\{liudq2026,huxx\}@lzu.edu.cn}\\[4pt]
\small\textsuperscript{*}Corresponding author}
\date{}
\begin{document}
\maketitle
\begin{abstract}
A financial market can appear mainly dependent on its own shocks under a variance account and mainly dependent on external shocks under a cumulative-response account. We establish how these opposite diagnoses arise from the same propagation system. For a given cumulative magnitude, squaring gives greater weight to responses concentrated across sources and time. Source priorities can therefore reverse even in stable Gaussian systems with known parameters. Heavy tails raise a separate question: whether the reported shares have a fixed population basis. One innovation with an infinite second moment removes the full network's finite covariance input. We characterize when normalized shares nevertheless stabilize and construct rank-amplitude connectedness, a cumulative-deviation target without finite first or second moments. It can be consistently estimated under local dependence and identification conditions. With fitted dynamics and rank scatter fixed, changing attribution replaces 3.14 names in the daily Top 10 of forty stocks, selecting recipients with greater external persistence relative to their own adjustment. In five futures contracts, mean external share rises from 28.0\% to 55.9\%, and 70.9\% of recipient--dates switch from own to external dominance.
\end{abstract}

\section{Introduction}\label{sec-intro}
Financial surveillance uses connectedness networks to locate sources of stress and identify recipients exposed to shocks elsewhere. The Diebold--Yilmaz framework assigns each source a share of a recipient's forecast-error variance \citep{diebold2009measuring,diebold2014network}. Designed for crisis warning and crisis tracking \citep{diebold2012better}, connectedness also informs the IMF's SysMo toolkit and ESMA's analysis of transmission among EU investment funds \citep{blancher2013sysmo,bouveret2023connectedness}. Yet what underlies the judgments about sources of stress drawn from these networks? Do the computed shares correspond to a fixed population network? Even when such a network exists, do its source priorities still depend on how shock responses are attributed?

Consider a risk committee using a connectedness dashboard to choose markets for sustained review. In one data history, an exceptional currency move makes foreign exchange appear to be the main source of stress; in another, an equity-market shock shifts attention to stocks. Suppose the propagation system is unchanged. If the influence of such extremes survives as the record grows, the committee can keep reallocating its attention without approaching a stable source assessment. The question is whether more data will reveal a persistent pattern or keep changing which realized event dominates the dashboard.

Generalized forecast-error variance decomposition (GFEVD) requires a finite covariance matrix for its population variance benchmark. One innovation with an infinite second moment removes that input for the full network, even if every other coordinate has finite variance. Sample shares remain computable. Their limits depend on how extreme observations enter the numerator and denominator: normalization may cancel their fluctuations or leave a random limit. Even when covariance exists, the shares can have nonstandard estimation precision. Section~\ref{sec-main-method} separates these outcomes.

Now consider a treasury desk renewing a fixed amount of overnight funding each day. One market shock raises its funding spread by 100 basis points for one day; another raises it by 20 basis points for ten days. Summing squared daily spread responses ranks the first source higher, but the second generates twice the additional interest expense over the ten-day window. A source list designed to capture large squared responses can therefore direct attention away from the largest cumulative funding cost. Even with both response paths known exactly, the appropriate review priority depends on what the treasury desk is trying to measure.

The funding example separates two accounting objectives. \emph{Response-energy attribution} sums squared responses; \emph{cumulative-deviation attribution} sums their absolute magnitudes across horizons. For a given cumulative magnitude, squaring gives greater weight to concentration. Cumulative accounting counts each unit of departure, so a later reversal does not erase earlier deviations. Our Gaussian construction in Section~\ref{sec-gaussian} gives opposite source priorities with all moments finite and with known dynamics and shock inputs held fixed. We call the dependence of shares and priorities on the accounting rule, with response paths fixed, \emph{attribution dependence}.

\begin{figure}[!htbp]
\centering
\captionsetup{font={small,stretch=1.05}}
\includegraphics[width=\linewidth,height=.49\textheight,keepaspectratio]{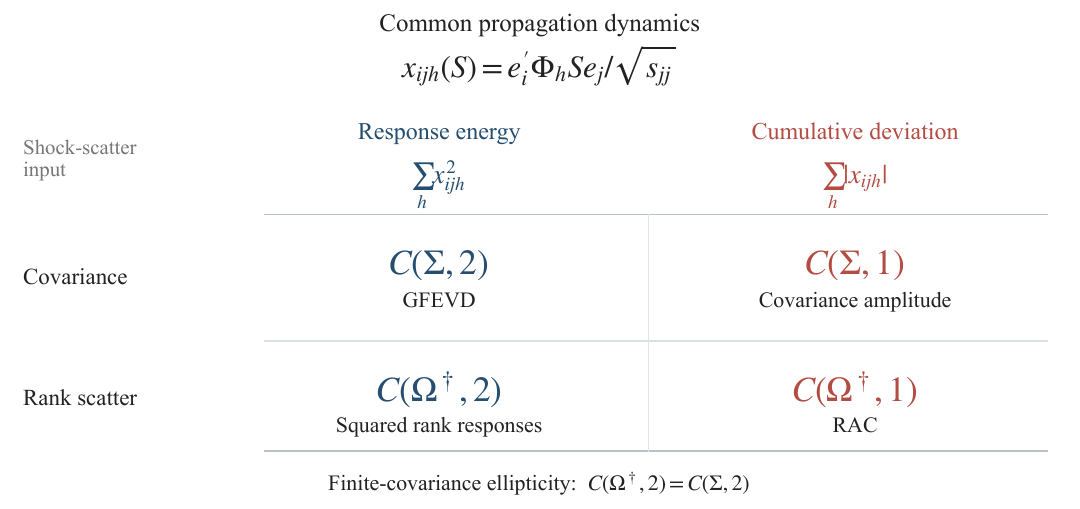}
\caption{\textbf{Two measurement choices under a common propagation system.} Rows change the shock-scatter input; columns change response-energy versus cumulative-deviation attribution. Horizontal comparisons isolate what the network allocates; vertical comparisons isolate the input used to generate its responses. Covariance requires second moments; rank scatter uses quantiles and concordance. GFEVD and RAC occupy opposite corners, and the crossed allocations separate the two changes.}
\label{fig-framework}
\end{figure}

Rank-amplitude connectedness (RAC) supplies a cumulative-deviation target without finite first or second moments. Marginal quantile scales and rank dependence define its shock-scatter input; absolute responses define its attribution rule. Under finite-covariance ellipticity, rank scatter and covariance give the same squared-response allocation. This calibration separates replacement of the shock input from replacement of the accounting rule. Figure~\ref{fig-framework} displays the two choices: rows change the input and columns change attribution, with propagation fixed. Horizontal comparisons also isolate attribution outside ellipticity. Under a fixed external/own cumulative amplitude ratio, an exact identity shows how source breadth and longer relative external duration lower the squared external share.

In forty stocks, holding fitted dynamics and rank scatter fixed while changing attribution replaces 3.14 names in the daily Top 10 and moves recipient ranks by 5.05 places on average. The cumulative-deviation list selects recipients whose external responses persist longer relative to their own adjustment. Thus the change in names corresponds to a change in the response profiles brought under review.

In five futures contracts, the difference changes whether responses are attributed mainly to the contract itself or to other markets. Under the same common-input comparison, mean external share rises from 28.0\% to 55.9\%, and 70.9\% of recipient--dates switch from own to external dominance. Corn illustrates the effect of external responses lasting longer than its own response; E-mini futures illustrate the combined contribution of several external sources with similar durations. Persistent or widely distributed responses can dominate the cumulative account even when own responses dominate the squared account.

Both findings recur in public data beyond these panels. In 96 global banks, common-input Top 24 lists differ by 7.22 names on average. In the 88 industrial sectors of \citet{brunner2023estimating}, fixing fitted dynamics and rank scatter yields own-to-external reversals in 91.6\% of sector--month observations over a three-month horizon. Broad contributions across sectors reproduce the source-breadth mechanism in production growth.

Our direct benchmark is the GFEVD underlying Diebold--Yilmaz connectedness. Refinements of this benchmark address how shocks are identified and how their contributions are combined. \citet{lanne2016generalized} construct generalized shares that sum to one; \citet{lastrapes2021joint} use joint conditional forecasts to measure aggregate spillovers. \citet{buchwalter2026clustered} allow dependence within clusters while orthogonalizing shocks across them. \citet{caloia2019normalization} show that normalization changes net-spillover signs and magnitudes. These results make clear that choices within variance decomposition affect network assessments. Our comparison reaches beyond these choices: for the generalized responses of \citet{pesaran1998generalized}, opposite source priorities arise with both responses and row normalization fixed. The difference lies in what each source is credited with. We also establish when normalized covariance shares have a deterministic limit under heavy tails. Summing to one does not ensure that a reported network estimates a fixed allocation.

An important extension of this benchmark is the quantile-connectedness framework of \citet{ando2022quantile}. Separating common from idiosyncratic variation, they show that small credit-risk shocks remain largely local, whereas shocks in both tails spread strongly through the sovereign--financial system. Their result changes the diagnosis of external dependence by changing the conditional state. We show that the diagnosis can also reverse within a fixed response system: own responses can dominate the energy account while external responses dominate cumulative deviations. The distinction is between state dependence and attribution dependence. At each quantile, their FEVD still allocates forecast-error variance using a finite covariance input. Our analysis establishes a cumulative-deviation population target beyond finite covariance and proves that attribution reversals also arise within stable Gaussian systems.

Related work models nonlinear tail interdependence through TENET \citep{hardle2016tenet}, traces distributional dynamics with quantile impulse responses \citep{gao2026quantilevar}, and distinguishes good- and bad-volatility spillovers \citep{barunik2016asymmetric}. These approaches refine the state or sign of the responses being examined.

Financial assessment also depends on how its objectives are implemented. \citet{glasserman2025fair} show why applying a common stress-test model need not treat heterogeneous banks fairly. \citet{ghamami2022collateralized} show how collateral access and payment timing shape defaults and clearing. Their objects are projected losses and clearing payments; ours is the source attribution of a given dynamic response. The paper makes two foundational contributions. First, it establishes when a computed variance network has a fixed population basis and constructs an estimable cumulative-deviation target beyond finite covariance. Second, it characterizes the two accounting objects and derives exact conditions for channel and own/external reversals within the domain of valid Gaussian variance attribution. The financial evidence shows which recipients and external dependencies these choices bring into view.

Section~\ref{sec-main-method} establishes the covariance boundary. Section~\ref{sec-rac} characterizes the two accounting objects and their calibrated comparison. Section~\ref{sec-mc} connects estimation to priority recovery. Section~\ref{sec-empirical} examines recipient and source priorities across stocks, banks, currencies and futures. The Electronic Companion contains proofs, full simulation designs and supplementary empirical evidence.

\section{Population foundations of variance connectedness}\label{sec-main-method}
A row of shares always sums to one, but this arithmetic does not ensure that larger samples recover a fixed network. We first state the joint covariance requirement, then determine when normalization cancels extreme fluctuations and how estimation behaves when covariance exists.

\Needspace{6\baselineskip}
\subsection{The joint covariance requirement}\label{sec-estimand}
We follow a shock from its source to each recipient over $H$ horizons. For $K$ variables, $\Phi_h(u)$ describes propagation at horizon $h$ and evaluation date $u$, with $\Phi_0(u)=I$. The horizons run from impact through $H-1$. In a VAR(1) frozen at $u$, these matrices are powers of its coefficient matrix; finite-order VARs use the usual moving-average recursion.\footnote{Variance-based connectedness has been extended to larger systems \citep{demirer2018estimating}, time-varying dynamics \citep{antonakakis2020dynamic} and frequency decompositions \citep{barunik2018measuring}; \citet{diebold2023past} review its development.}

For a positive-diagonal scatter matrix $S$, define the generalized response of recipient $i$ to a shock from source $j$ as
\[
x_{ijh}(S;u)=\frac{e_i^\top\Phi_h(u)S e_j}{\sqrt{s_{jj}}}.
\]
Here the coordinate vectors $e_i$ and $e_j$ select the recipient and source. The vector $S e_j/\sqrt{s_{jj}}$ fixes the shock direction and scale. Covariance gives the GFEVD construction \citep{pesaran1998generalized}; replacing scatter changes the shock input while keeping propagation fixed.

A source's share compares its accumulated responses with those from all sources. For response power $r>0$, define
\begin{equation}
C_{ij}(S,r;u)=\frac{\sum_{h=0}^{H-1}|x_{ijh}(S;u)|^r}
{\sum_{k=1}^{K}\sum_{h=0}^{H-1}|x_{ikh}(S;u)|^r}.
\label{eq-main-allocation}
\end{equation}
Each row sums to one, with recipients indexing rows and sources indexing columns. Generalized shocks may be correlated.

Covariance and squared responses recover row-normalized GFEVD, $C(\Sigma,2)$: the forecast-error variance term in the unnormalized decomposition cancels within each row. This covariance construction supplies the conventional variance interpretation; a general scatter input defines response-energy shares.

The full $K$-node covariance input exists only when every innovation coordinate has a finite second moment. If $E\varepsilon_j^2=\infty$ for even one node $j$, no finite $K\times K$ covariance matrix exists, regardless of the moments of the remaining nodes.

Positive changes of measurement units leave the allocation and its rankings unchanged. Multiplying variable $i$ by $b_i>0$ and transforming dynamics and scatter consistently multiplies row $i$'s responses by $b_i$; the factor $b_i^r$ cancels from its shares.

\subsection{Deterministic limits under infinite variance}\label{sec-moment-boundary}
The same extreme observations enter a share's numerator and denominator. Their effects may cancel in the ratio even when covariance has no finite limit. Let $\nu$ denote the regular-variation index of the innovation norm; its square has index $\nu/2$. 
\begin{proposition}[When normalized covariance shares have a deterministic limit]\label{prop:fixed-network}
Fix $K\ge2$, finite $H\ge1$ and known propagation matrices with $\Phi_0=I$. For $0<\nu<2$, suppose positive-diagonal sample covariance inputs satisfy $c_T\widehat\Sigma_T\Rightarrow S_\nu$ for a deterministic tail-dependent normalization $c_T>0$ (specified in Electronic Companion Section~EC.1.3), where $S_\nu$ is a positive-semidefinite $\nu/2$-stable random matrix with strictly positive diagonal almost surely. Then
\[
C(\widehat\Sigma_T,2)\Rightarrow C(S_\nu,2).
\]
For any deterministic matrix $C_0$,
\[
C(\widehat\Sigma_T,2)\xrightarrow{p}C_0
\quad\Longleftrightarrow\quad
\Pr\{C(S_\nu,2)=C_0\}=1.
\]
If the support of $S_\nu$ contains a nonempty open subset of the positive-definite cone, every individual limiting share is nondegenerate. No individual sample share then converges in probability to a deterministic constant.
\end{proposition}

Under the stated support condition, knowing propagation exactly and collecting more data do not recover fixed source shares: the limiting allocation remains random. Electronic Companion Section~EC.1.3 proves the result and gives a condition on tail directions ensuring nondegeneracy.

Exact cancellation occurs in a system driven by one common scalar shock. For innovations $Z_t v$, $A=0$, $H=1$ and every $v_j\ne0$, all sample shares equal $1/K$ even when the iid continuous scalar $Z_t$ has infinite variance. Fluctuations along this single scatter direction cancel in every share. The deterministic ratios have no finite-population-variance interpretation.

Table~\ref{tbl-moment-boundary} separates covariance existence from the sampling behavior of the shares. At $\nu=2,4$, moment existence depends on the slowly varying tail factor; the corresponding Student moments diverge.

\begin{table}[!htbp]\centering\small\renewcommand{\arraystretch}{1.05}\setlength{\tabcolsep}{5pt}
\caption{Moment regimes for covariance-based connectedness.}\label{tbl-moment-boundary}
\begin{tabular}{p{.12\linewidth}p{.30\linewidth}p{.48\linewidth}}\toprule
Tail index & Covariance input & Implication for connectedness\\\midrule
$0<\nu<2$ & At least one second moment is infinite & Sample shares remain computable without a finite-covariance attribution; their limits depend on response geometry.\\[3pt]
$2<\nu<4$ & Covariance exists; a quadratic input has infinite variance & Consistency can hold, but a noncancelling share can inherit a slower, non-Gaussian covariance limit.\\[3pt]
$\nu>4$ & Fourth moments exist & The fourth-moment requirement is satisfied; conventional covariance CLTs also require the stated dependence and regularity conditions.\\\bottomrule
\end{tabular}
\par\smallskip\begin{minipage}{.97\linewidth}\footnotesize\raggedright\textit{Notes.} The table concerns covariance-based squared allocation $C(\widehat\Sigma_T,2)$. Sampling statements use the dependence, identification, limit and allocation-map conditions in the Electronic Companion.\end{minipage}
\end{table}

\subsection{Nonstandard precision with finite covariance}\label{sec-share-precision}
For $2<\nu<4$, the population variance network exists. The question is how quickly sample shares approach it. Here the scaling applies to covariance estimation error; Proposition~\ref{prop:fixed-network} instead rescales the covariance level. Suppose $a_T(\widehat\Sigma_T-\Sigma)\Rightarrow S_\nu$, where $a_T\to\infty$ is the valid covariance normalization and $S_\nu$ is a $\nu/2$-stable matrix. With fixed propagation and positive diagonal covariance,
\[
a_T\{C_{ij}(\widehat\Sigma_T,2)-C_{ij}(\Sigma,2)\}
\Rightarrow DC_{ij,\Sigma}[S_\nu].
\]
The derivative measures how a small covariance change affects the share. With full-dimensional stable fluctuations and positive-definite covariance outside an algebraic exceptional set, every individual share inherits the non-Gaussian first-order limit. If the derivative vanishes on the fluctuation support, the leading term cancels. Electronic Companion Section~EC.1.3 gives the normalization and support conditions. Estimated propagation contributes its joint estimation error \citep{davisresnick1986,davis1998sample,resnick2007heavy}.

A finite-covariance target, a deterministic limit of normalized sample shares, and conventional sampling precision are separate properties. Tail behavior and response geometry determine which properties the reported network possesses.\footnote{Sensitivity to the underlying distribution also arises in inventory pooling \citep{bimpikis2016inventory}, portfolio risk \citep{jin2021tail,lim2011conditional}, divergence-based uncertainty sets \citep{kruse2021toolkit} and tail-aware state-price-density estimation \citep{li2026tail}.}

\section{Attribution dependence: response energy and cumulative deviation}\label{sec-rac}
For the same response system, energy and cumulative deviation can imply opposite source priorities. To isolate this accounting choice, we first define rank scatter and calibrate it to the covariance benchmark.

\subsection{Population definition and elliptical calibration}\label{sec-rank-definition}
RAC uses marginal scales and rank dependence to define its shock input. MAD supplies scale without tail moments; Kendall's rank correlation measures dependence through bounded pairwise signs. The resulting scatter defines the population estimand $C(\Omega^\dagger,1)$ in Equation~\eqref{eq-main-allocation}. Collect the marginal scales in
\begin{align*}
D&=\operatorname{diag}(s_i),\\
s_i&=1.4826\,\operatorname{Med}|\varepsilon_i-\operatorname{Med}(\varepsilon_i)|,
\end{align*}
Write $\tau_{ij}$ for Kendall's rank correlation between innovation coordinates $i$ and $j$. Set $\widetilde R_{ij}=\sin(\pi\tau_{ij}/2)$ and let $R^\dagger$ be the Frobenius projection of $\widetilde R$ onto positive-semidefinite correlation matrices \citep{lindskog2003kendall,higham2002computing}. The scatter is
\[
\Omega^\dagger=DR^\dagger D.
\]
The sine transformation calibrates concordance to correlation under ellipticity. The nearest-correlation projection gives a positive-semidefinite scatter beyond that domain.

With continuous marginals, unique medians and positive MADs, the scatter is well defined beyond ellipticity. The own-source impact response is $\sqrt{\omega^\dagger_{ii}}>0$, so every row denominator is positive. The sums contain finitely many population responses. Both $r=1$ and $r=2$ therefore give finite shares without requiring first or second innovation moments.

Under continuous finite-covariance ellipticity, rank scatter is proportional to covariance \citep{fang1990symmetric}. Common scale cancels:
\[
C(\Omega^\dagger,2)=C(\Sigma,2).
\]
Outside ellipticity, $\Omega^\dagger$ retains its quantile-and-concordance definition but need not be proportional to covariance. Holding this input fixed continues to isolate attribution.

\subsection{Characterizing the accounting objects}\label{sec-choices}
The accounting choice determines whether surveillance accumulates deviations or gives greater weight per unit to larger departures. Both objectives add penalties across periods and treat signs symmetrically. Preserving a same-sign split treats a departure of size $a+b$ in one period like departures of sizes $a$ and $b$ in two periods. Quadratic scaling instead makes doubling a departure quadruple its penalty. Theorem~\ref{thm:response-accounting} characterizes the path score under each principle.\footnote{Other characterizations connect regulatory preferences to systemic-risk measures \citep{chen2013axiomatic} and risk measures to diversification indices \citep{han2025diversification}. Here the axioms concern how a given response path is scored.}

\begin{theorem}[Characterization of response accounting]\label{thm:response-accounting}
Let $V$ be a nonnegative score on finite response paths, invariant to changing any response's sign, with $V(0)=0$ and $V(1)=1$ for single-period paths. Suppose scores add across successive periods: $V(x,y)=V(x)+V(y)$, where $(x,y)$ concatenates two paths.

\textit{Cumulative deviation.} The score preserves same-sign splits, $V(a+b)=V(a,b)$ for all $a,b\ge0$, if and only if
\[
V(x)=\sum_h|x_h|.
\]
\textit{Response energy.} The score scales quadratically, $V(tx)=t^2V(x)$ for every $t\ge0$, if and only if
\[
V(x)=\sum_h x_h^2.
\]
Normalizing these scores across separately evaluated source paths gives $C(S,1)$ and $C(S,2)$, respectively.
\end{theorem}

In the log-realized-variance application, each response is the log ratio of perturbed and baseline variance paths. A surveillance objective charging the same penalty per unit of absolute log-variance displacement at each horizon uses the cumulative path loss. Each standardized source path is evaluated separately before shares are formed. Electronic Companion Section~EC.2.1 proves the characterization.

\subsection{Channel-level priority reversals}\label{sec-attribution-priorities}
For source $j$, let $a_j=\sum_h|x_{ijh}|$ denote cumulative amplitude and $b_j=\sum_hx_{ijh}^2$ response energy. For a channel with a nonzero response, effective duration is $d_j=a_j^2/b_j\in[1,H]$. With $p_{jh}=|x_{ijh}|/a_j$, the identity $d_j=1/\sum_h p_{jh}^2$ counts effective horizons: one for a single-period response and $H$ for equal absolute responses throughout. Holding cumulative amplitude fixed, a shorter effective duration raises response energy $b_j=a_j^2/d_j$.

\begin{proposition}[Attribution and channel priorities under common inputs]\label{prop:amplitude-concentration}
Fix a recipient, a common scatter input and a finite response horizon. For every active channel,
\begin{equation}
\frac{C_{ij}(S,2)}{C_{ij}(S,1)}
=\frac{b_j/a_j}{\sum_k b_k/\sum_k a_k}
=\frac{a_j/d_j}{\sum_k b_k/\sum_k a_k}.
\label{eq-main-attribution-tilt}
\end{equation}
The two allocation rows coincide if and only if $b_j/a_j$ is constant across active channels.
For two channels with $a_j>a_k>0$, amplitude allocation gives priority to $j$ and squared allocation gives priority to $k$ exactly when
\begin{equation}
1<\frac{a_j}{a_k}<\sqrt{\frac{d_j}{d_k}}.
\label{eq:duration-reversal}
\end{equation}
For any $H\ge3$, these opposite priorities occur on an open set of stable VAR systems with Gaussian innovations.
\end{proposition}

Equation~\eqref{eq-main-attribution-tilt} shows how the accounting choice changes the network: squaring increases a source's share when its energy per unit cumulative deviation exceeds the recipient-wide average. By Equation~\eqref{eq:duration-reversal}, a channel with greater cumulative amplitude loses squared priority when its longer duration outweighs that amplitude advantage.

\subsection{Opposite priorities with Gaussian innovations}\label{sec-gaussian}
\begin{figure}[!htbp]
\centering
\includegraphics[width=\linewidth]{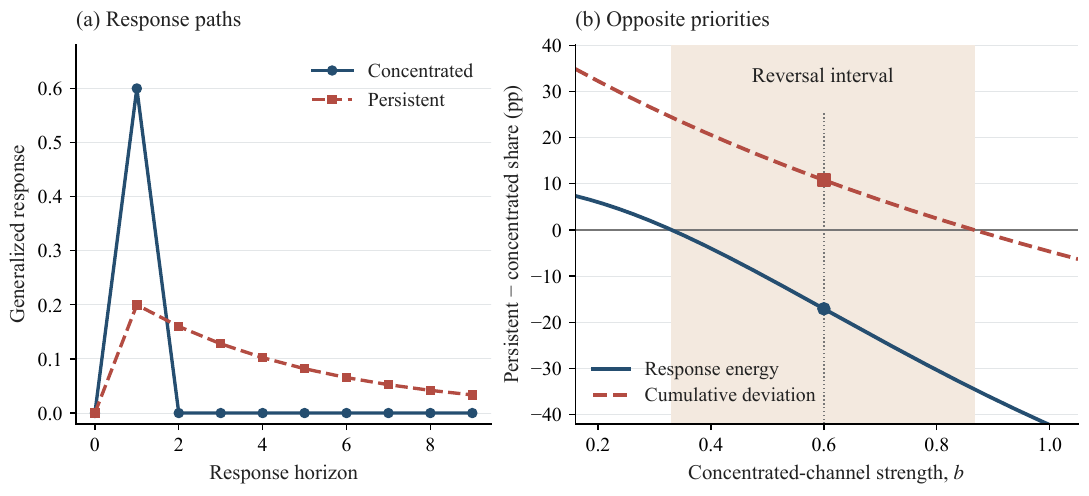}
\caption{\textbf{Aggregation alone reverses channel priorities over an interval.} A Gaussian VAR fixes $K=3$, innovation covariance $I_3$, $H=10$, and $A=\left(\begin{smallmatrix}0&b&0.2\\0&0&0\\0&0&0.8\end{smallmatrix}\right)$. (a) Response paths at $b=0.60$: source 2 produces a concentrated response; source 3 produces a smaller, persistent response at the same recipient. (b) Exact population share gaps $C_{13}-C_{12}$ as $b$ varies. Positive values favor the persistent channel. Opposite signs throughout $0.3303<b<0.8658$ identify the reversal interval. At each $b$, dynamics and scatter are identical across aggregation rules.}
\label{fig-main-allocation-mechanism}
\end{figure}

Figure~\ref{fig-main-allocation-mechanism} evaluates both allocations at known population parameters in a stable Gaussian system. GFEVD has its full variance interpretation and elliptical calibration is exact. One source produces a concentrated response and the other a smaller response that persists. With dynamics, innovation covariance and horizon fixed, GFEVD favors the concentrated channel and RAC the persistent channel.

When the concentrated response has magnitude $b=0.60$, RAC assigns 10.8 percentage points more to the persistent channel; GFEVD assigns 17.1 points more to the concentrated channel. The opposite ordering holds throughout $0.3303<b<0.8658$.

The strict reversal condition applies to any two active sources in a larger network. Moving from individual channels to all external sources adds breadth to the comparison: cumulative amplitude can be spread across sources as well as across horizons.

\subsection{Source breadth and relative duration}\label{sec-breadth-design}
The own/external diagnosis compares the combined external responses with the recipient's own response. Let $A_E=\sum_{j\ne i}a_j$ and $A_O=a_i$ denote external and own cumulative amplitudes for recipient $i$. Measure the effective number of external sources and their duration by
\[
n_E=\frac{A_E^2}{\sum_{j\ne i}a_j^2},\qquad
d_E=\frac{\sum_{j\ne i}a_j^2}{\sum_{j\ne i,h}x_{ijh}^2},
\]
with $d_O=a_i^2/\sum_hx_{iih}^2$. For positive external amplitude, $1\le n_E\le K-1$ and $1\le d_E,d_O\le H$. The effective source count $n_E$ equals the number of active external sources when their cumulative amplitudes are equal.

The product $n_Ed_E=A_E^2/\sum_{j\ne i,h}x_{ijh}^2$ measures effective breadth across source--horizon pairs: it is the inverse sum of squared shares of external absolute response. If $F_r=\sum_{j\ne i}C_{ij}(S,r)$, then
\begin{equation}
\frac{F_1}{1-F_1}=\frac{A_E}{A_O},\qquad
\frac{F_2}{1-F_2}=\left(\frac{A_E}{A_O}\right)^2\frac{d_O}{n_Ed_E}.
\label{eq-external-response-odds}
\end{equation}
An externally dominated amplitude allocation and an own-dominated squared allocation occur precisely when
\begin{equation}
F_1>\tfrac12>F_2
\quad\Longleftrightarrow\quad
1<\frac{A_E}{A_O}<\sqrt{\frac{n_Ed_E}{d_O}}.
\label{eq-external-dominance}
\end{equation}
Equation~\eqref{eq-external-dominance} concerns aggregate own/external dominance; Equation~\eqref{eq:duration-reversal} ranks individual channels. The aggregate reversal requires $n_Ed_E>d_O$. A single external source can satisfy it through longer duration; several sources add breadth.

\begin{figure}[!htbp]\centering
\includegraphics[width=\linewidth]{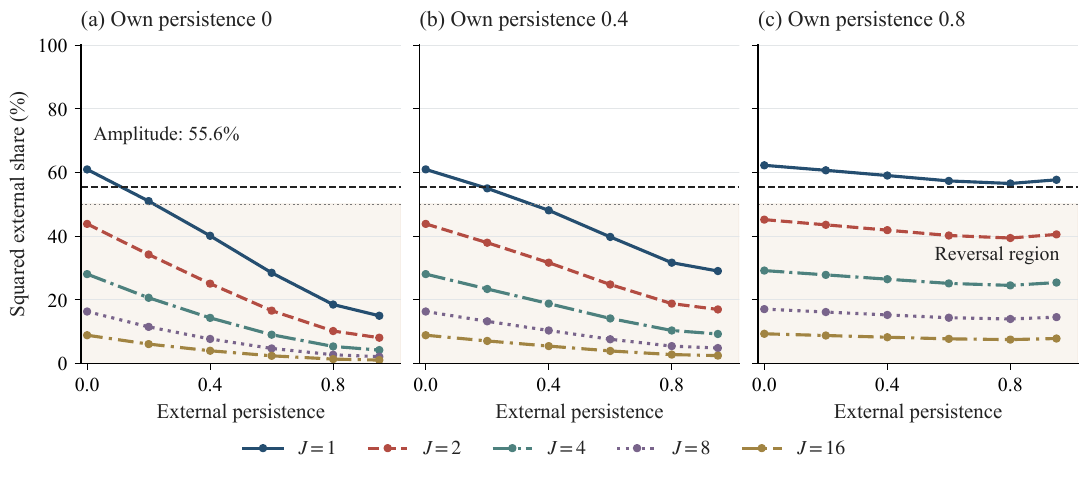}
\caption{\textbf{At a fixed amplitude ratio, breadth and duration can reverse source dominance.} All cells use $H=10$ and external amplitude share $1.25/2.25=55.6\%$. Curves give squared external percentages for each source count over six external-persistence values; panels vary own persistence. The dotted line marks 50\%; the dashed black line is the common 55.6\% amplitude share. Squared shares below 50\% identify reversals.}\label{fig-suite-geometry}
\end{figure}

Figure~\ref{fig-suite-geometry} isolates breadth and duration in a Gaussian star system with $J$ equally strong external channels. Let $r=A_E/A_O$, and write $D_e$ and $D_o$ for the common external effective duration and the own effective duration, respectively. The shares are
\[
F_1=\frac{r}{1+r},\qquad
F_2=\frac{r^2}{r^2+JD_e/D_o}.
\]
At a fixed external/own amplitude ratio, broader sources and longer relative duration reduce the external energy share while leaving cumulative-deviation attribution unchanged.

Figure~\ref{fig-suite-geometry} fixes $H=10$ and $r=1.25$, giving an external amplitude share of 55.6\% in every cell. At zero own and external persistence, splitting one external channel into four lowers the squared external share from 61.0\% to 28.1\%. With one source, increasing its persistence from zero to 0.8 lowers that share to 18.5\%. The three panels vary own persistence: a longer own response changes the relative-duration threshold, even with the same cumulative external/own ratio.

Electronic Companion Section~EC.2.5 gives the star construction and complete parameter grid.

The same breadth--duration adjustment governs recipient rankings. Write $r_i=A_{E,i}/A_{O,i}$ and $g_i=n_{E,i}d_{E,i}/d_{O,i}$. Amplitude ranks recipients by $r_i$, while squared allocation ranks them by $r_i^2/g_i$. For $r_i>r_k>0$, the rankings oppose exactly when
\[
1<\frac{r_i}{r_k}<\sqrt{\frac{g_i}{g_k}}.
\]

\FloatBarrier
\section{Estimation and recovery of network priorities}\label{sec-mc}
To recover monitoring priorities, an estimated network must select the right recipients and classify their responses as own- or externally dominated. Small share errors can change names near the list cutoff while leaving the list's mean external share almost unchanged. We give conditions for recovering the list and the classification, and bound the external share forgone when the list differs.

\subsection{The surveillance allocation problem}\label{sec-summaries}
Incoming strength $f_i=FROM_i=100(1-C_{ii})$ is the recipient's external share. It measures relative dependence: row normalization removes total response magnitude.

Suppose monitoring capacity permits reviewing a fixed number of recipients with the greatest relative external dependence. For a rank cutoff $1\le q<K$, the list $I_q$ of the largest $f_i$ maximizes mean external share, $V_f(I)=q^{-1}\sum_{i\in I}f_i$. Its objective loss from choosing another list $J$ is
\begin{equation}
L_f(J)=V_f(I_q)-V_f(J).
\label{eq-monitoring-loss}
\end{equation}
Response-energy accounting selects the list with the greatest mean external variance share; cumulative-deviation accounting selects the list with the greatest mean external share of absolute response. List overlap records whether the same recipients are selected. The loss $L_f$ records the external share forgone under the chosen objective.

The relevant accuracy is relative to the gap at the list cutoff and the distance from the own/external boundary (Proposition~\ref{prop-monitoring-resolution}).

\begin{proposition}[Recovery of network priorities]\label{prop-monitoring-resolution}
For recipient $i$, let $x_{ijh}$ and $\widehat x_{ijh}$ be target and fitted responses with positive cumulative totals. Put
\[
\begin{aligned}a_{ij}&=\sum_h|x_{ijh}|,\qquad A_i=\sum_j a_{ij},\\
\delta_i&=\sum_{j,h}|\widehat x_{ijh}-x_{ijh}|,\end{aligned}
\]
and define $\widehat a_{ij},\widehat A_i$ analogously. Let $f_i=100\sum_{j\ne i}a_{ij}/A_i$ and let $\widehat f_i$ be its fitted counterpart. Then
\begin{equation}
\begin{aligned}|\widehat f_i-f_i|&\le
\min\left\{100,\frac{100\delta_i}{\max(A_i,\widehat A_i)}\right\}
\\&\le \epsilon_i:=\min\left\{100,\frac{100\delta_i}{\widehat A_i}\right\}.
\end{aligned}\label{eq-response-score-bound}
\end{equation}
Let $I_q$ and $\widehat I_q$ maximize the target and fitted sums over lists of size $q$. Put $\epsilon=\max_i\epsilon_i$ and $d=|I_q\setminus\widehat I_q|$.
\Needspace{4\baselineskip}
Consequently:
\begin{enumerate}
\item A target gap $f_{(q)}-f_{(q+1)}>2\epsilon$ gives $\widehat I_q=I_q$. Alternatively, the same condition on the fitted cutoff gap certifies that its selected list is target-optimal.
\item If $|\widehat f_i-50|>\epsilon_i$, fitted and target own/external classifications agree.
\item The selected list satisfies
\[
0\le \frac1q\left(\sum_{i\in I_q}f_i-\sum_{i\in\widehat I_q}f_i\right)
\le \frac{2d\epsilon}{q}\le2\epsilon.
\]
\end{enumerate}
\end{proposition}

The deterministic bound includes dynamic and scatter estimation error in $\delta_i$. The companion proves the proposition and evaluates $\delta_i$ using known target paths in complete-estimation simulations.

Applying standard connectedness summaries to $C(S,r)$, outgoing strength identifies sources contributing across recipients, and total connectedness averages incoming strength:
\begin{align*}
TO_j&=100\sum_{i\ne j}C_{ij},\\
FROM_i&=100\sum_{j\ne i}C_{ij},\\
NET_i&=TO_i-FROM_i,\\
TCI&=K^{-1}\sum_i FROM_i.
\end{align*}
Spearman correlation compares rankings across all nodes. Displayed lists contain ten stocks or two smaller-panel instruments. Let $E_V$ and $E_R$ contain the $q_E$ largest external links under GFEVD and RAC, respectively. Their link replacement rate is $1-|E_V\cap E_R|/q_E$. 

\subsection{Estimating the responses that determine priorities}\label{sec-estimation}
Under the conditions below, the estimator consistently recovers recipient shares in a fixed stationary model or at an interior rescaled date of a time-varying model. A positive population gap at the cutoff then gives exact Top-$q$ recovery with probability tending to one. Recipients whose population shares lie away from the own/external boundary are consistently classified. Even with a tied cutoff, the selected list's mean external share converges to the population optimum.

\paragraph*{Construction.}
The current window estimates propagation by self-weighted local-linear Huber regression. Self-weighting limits extreme-lag influence, and the bounded Huber score controls innovation influence. Residual MADs and projected sine--Kendall dependence estimate the rank-scatter input.\footnote{The construction connects bounded-score regression \citep{huber1981robust} with heavy-tail autoregressive estimation \citep{davis1992mestimation,ling2005selfweighted} and local modeling \citep{gao2024tvp}. Related robust methods address multivariate scatter \citep{maronna1976robust,tyler1987distribution}, scale \citep{rousseeuw1993alternatives} and high-breakdown covariance \citep{rousseeuw1999fastmcd}. \citet{mikosch2013heavy} analyze how heavy-tailed noise enters finite-sample OLS estimates.} Earlier observations guide coefficient regularization and correlation shrinkage. The resulting coefficient adjustment, the fraction of high-leverage observations deleted and the amount of scatter shrinkage all vanish asymptotically. Electronic Companion Sections~EC.3.2--EC.3.3 specify the criterion, tuning and complete algorithm.

\paragraph*{Conditions and guarantee.}
Dynamics are smooth and uniformly contracting. Innovations are stationary ergodic with positive conditional densities and a positive fractional moment, which may be below one. The distributions defining residual medians and MADs must cross their target quantiles uniquely; ellipticity is not required.

Identification requires the location targeted by the conditional Huber score to be the same across past histories at each scale. A free intercept absorbs this common location, leaving the lag coefficients and translation-invariant scatter unchanged. This condition covers asymmetric iid innovations and conditionally symmetric innovations with persistent scale variation.

For fixed dimension, lag order and response horizon, these conditions give
\[
\max_{i,j}\left|\widehat C_{ij}(\widehat\Omega^\dagger,1;u)
                   -C_{ij}(\Omega^\dagger(u),1;u)\right|\xrightarrow{p}0.
\]
The guarantee includes pilot-scale estimation, coordinate calibration, slope-penalty selection and historical regularization. Electronic Companion Section~EC.3.5 gives the local conditions and proof; Proposition~EC.5 establishes consistency of the complete estimator. Proposition~\ref{prop-monitoring-resolution} converts response consistency into the recovery conclusions above.

\subsection{Population stability across moment regimes}
The three-node experiment isolates scatter estimation by holding dynamics and squared aggregation fixed. Seven innovation laws span tail indices $.8,1.5,2.1,3,4.1,8$ and Gaussian shocks, with 300 independent samples per law and nested sizes $n=125,500,2000,8000$. The companion specifies the population inputs.

\begin{figure}[!htbp]\centering
\includegraphics[width=\linewidth]{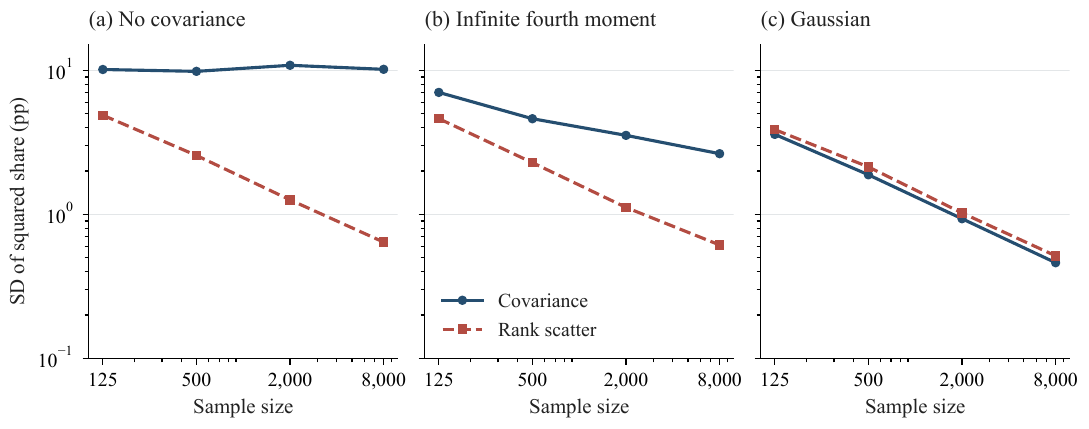}
\caption{\textbf{Growing samples need not stabilize covariance shares.} Dispersion of $C_{12}$ under the same squared-response allocation across 300 independent innovation samples, with known dynamics and $H=10$. Sample sizes are nested within each path; both axes are logarithmic. Panels use Student tail indices 1.5 and 3 and a Gaussian benchmark. Covariance shares in the first panel are sample ratios without a population variance interpretation.}\label{fig-main-boundary}
\end{figure}

At tail index 1.5, increasing the sample from $n=125$ to $n=8000$ leaves the standard deviation of a sample squared share at 10.16 and 10.19 percentage points. Replacing covariance with rank scatter, with squared aggregation unchanged, reduces dispersion from 4.89 to 0.64 points (Figure~\ref{fig-main-boundary}). In the rank-one control, normalization instead cancels the covariance fluctuations and the shares remain constant at every sample size.

The forty-variable factor and directed-block designs measure recovery when both dynamics and scatter are estimated. Stationary and time-varying dynamics combine with four elliptical Student laws, Gaussian innovations and independent coordinates with unequal tail indices, giving 24 cells with 100 independent paths each. Each path uses its own initial 500 observations for calibration and tuning, followed by estimation at $m=125,250,500$ and $H=10$.

Recovery is evaluated against each allocation's own population target using
\[
\mathcal L_K=K^{-1}\sum_i\sum_{j\ne i}|\widehat C_{ij}-C_{ij}|.
\]
Multiplication by 100 expresses off-diagonal allocation error per recipient in percentage points. Elliptical calibration gives rank-squared and covariance-squared allocations a common target where covariance exists. With the same fitted dynamics, rank scatter reduces mean allocation error for the squared-response target by 17.58 points at $\nu=3$, with paired 95\% Monte Carlo interval $[16.44,18.72]$. Table~\ref{tbl-suite-squared} reports both networks and the finite-fourth-moment and Gaussian comparisons. Its common-dynamics columns isolate scatter estimation; the local-linear OLS column compares complete fitting procedures.

\begin{table}[!htbp]\centering\small\renewcommand{\arraystretch}{1.05}\setlength{\tabcolsep}{4pt}
\caption{Recovering a common squared-response target, $m=250$.}\label{tbl-suite-squared}
\begin{tabular}{llrrr}\toprule
& & \multicolumn{2}{c}{Common RAC dynamics} & Local-linear\\ \cmidrule(lr){3-4} Network & Innovation law & Covariance & Rank scatter & OLS\\\midrule

Factor & 3 & 36.78 & 21.81 & 67.74\\
Factor & 8 & 20.85 & 19.98 & 50.92\\
Factor & Gaussian & 17.54 & 19.07 & 48.02\\
\addlinespace[3pt]
Blocks & 3 & 43.73 & 23.54 & 96.40\\
Blocks & 8 & 20.76 & 21.76 & 84.03\\
Blocks & Gaussian & 16.96 & 20.40 & 79.85\\
\bottomrule\end{tabular}
\par\smallskip\begin{minipage}{.98\linewidth}\footnotesize\raggedright\textit{Notes.} Entries are mean off-diagonal allocation errors per recipient in percentage points. Stationary and time-varying cells have equal weight. Finite-covariance elliptical calibration gives all columns the same population target. OLS uses sample covariance and the same local-linear design and current window. Each network--law row contains 200 independent paths.\end{minipage}
\end{table}

\subsection{Recovering recipient lists and source dominance}\label{sec-simulation-monitoring}
\begin{figure}[!htbp]\centering
\includegraphics[width=\linewidth,height=.70\textheight,keepaspectratio]{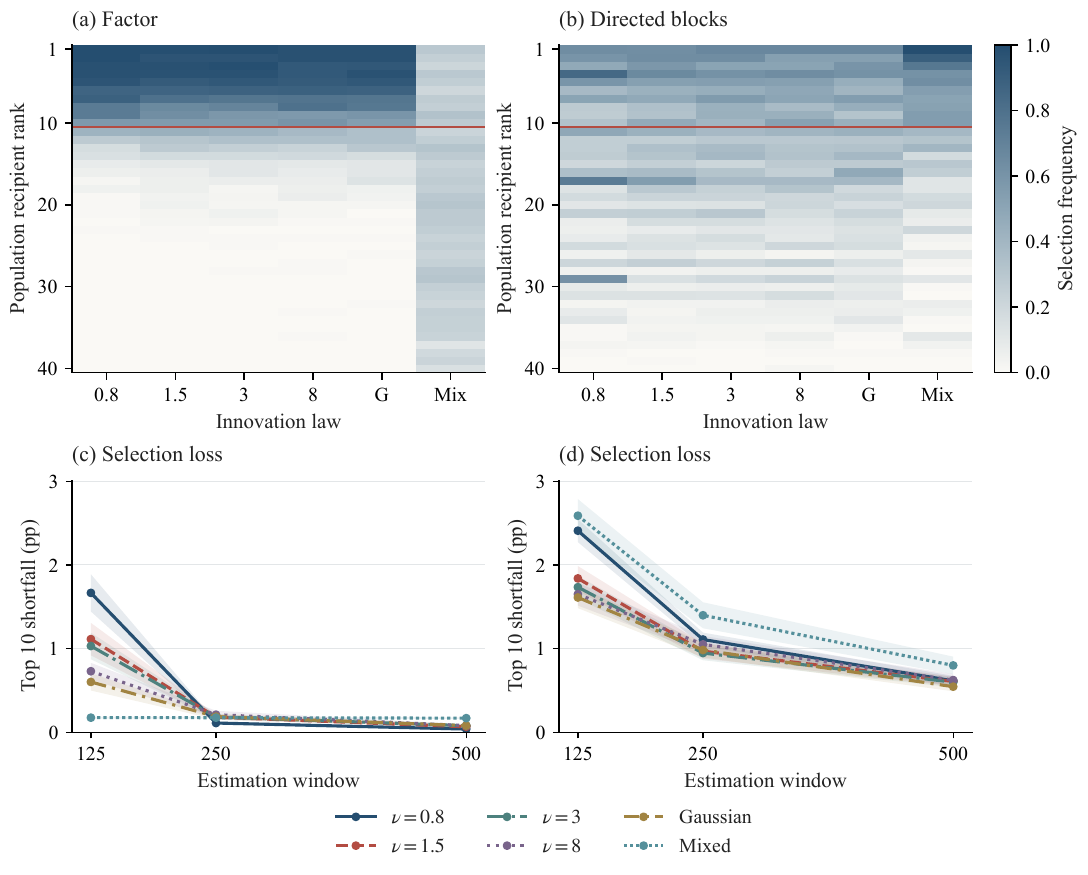}
\caption{\textbf{Recovering recipient lists and their population external share.} (a,b) Top 10 selection frequencies in time-varying systems at $m=250$. Rows order all forty nodes by population amplitude receiving strength within each column; the horizontal line marks the population Top 10. G denotes Gaussian innovations and Mix independent Student coordinates with unequal tail indices. (c,d) Mean population-strength shortfall of the selected ten recipients as the estimation window changes. Shaded bands are $\pm1.96$ Monte Carlo standard errors.}\label{fig-suite-monitoring}
\end{figure}

Figure~\ref{fig-suite-monitoring} distinguishes recovering the population Top 10 from attaining its mean external share. The upper panels show how often each recipient is selected; the lower panels report the external share forgone by the fitted list. Both use estimated dynamics and scatter, with 100 paths in each network--innovation cell.

At $m=250$, Top 10 population-strength shortfall ranges from 0.10 to 1.51 percentage points across the 24 cells, and mean recipient-share error ranges from 2.24 to 5.98 points. The independent-margin factor design tests list recovery under weak population separation: all forty population scores span 0.667 points and the cutoff gap is 0.000612 points. Selected-list shortfall is 0.17--0.18 points, against 0.192 for random selection; Top 10 overlap is 26.3--27.8\%, against 25\% for a random ten-name list. Near-tied recipients make names difficult to distinguish, while alternative lists can attain almost identical objective values. The companion pairs recovery with population separation and reports classification by source category and distance from the dominance boundary.

At $m=250$, own/external classification accuracy ranges from 98.0\% to 100.0\% across the 24 cells. Each elliptical design has externally dominated population recipients; each independent-margin design has own-dominated recipients. The extended designs vary both categories within a network and their distance from the boundary.

A further 26 designs span null, factor, block, dense, signed, star and weakly correlated systems. Across twelve paths per design, recipient-share error at the final endpoint averages 4.44 points and Top 10 shortfall 0.75 points. In the Gaussian star designs, the focal external recipient is recovered in eight of twelve paths at population share $.6$ and all twelve at $.8$.

\FloatBarrier
\section{Financial monitoring: recipients and sources}\label{sec-empirical}
The stock evidence follows changes in recipient lists to the response profiles behind them. The futures evidence asks whether a contract's responses are predominantly own or external, and which sources call for joint assessment. Fitted dynamics remain the same in every comparison. GFEVD--RAC comparisons change both scatter and attribution; holding rank scatter fixed isolates attribution, without requiring ellipticity.
\subsection{Data and common TVP design}\label{sec-data}
We model the logarithm of RV5, the realized-variance measure constructed from five-minute returns, in VOLARE panels of forty stocks, five currency pairs and five futures contracts \citep{cipollini2026volare}. Following \citet{diebold2014network}, connectedness describes volatility propagation.\footnote{This object differs from other systemic-risk summaries: principal components and Granger-causality networks \citep{billio2012econometric}, system tail risk conditional on institutional distress \citep{adrian2016covar}, and expected capital shortfall in a severe market decline \citep{brownlees2017srisk}. \citet{benoit2017risks} distinguish mechanism-based analyses from broad market-based measures.} Realized variance links intraday price variation to daily volatility, and log-volatility dynamics accommodate persistence \citep{andersen2003realized,corsi2009simple}. Stocks run from January 2, 2015 to September 30, 2026. Currency pairs AUD/USD, EUR/USD, GBP/USD, USD/CAD and USD/JPY begin on September 25, 2009; futures C, CL, ES, GC and NG begin on September 28, 2009. Both smaller panels end on September 30, 2026.

Each market is fitted separately, so external shares refer to the other assets within that panel. The first 500 positions provide fixed median/MAD coordinate calibration and select the slope-penalty multiplier: $\gamma=10$ for stocks and futures and $\gamma=3$ for currencies. The current window contains $m=250$ original market-date positions. Historical templates use up to 1,500 observations before that window. Fits retain originally adjacent complete lag--response pairs and use $H=10$ response horizons. The estimates cover 2,452 stock dates, 3,914 currency dates and 3,900 futures dates; every primary fitted coefficient matrix has spectral radius below one.

\subsection{Network-level tail diagnostics}\label{sec-tail-diagnostics}
\begin{figure}[!htbp]\centering
\includegraphics[width=\linewidth,height=.63\textheight,keepaspectratio]{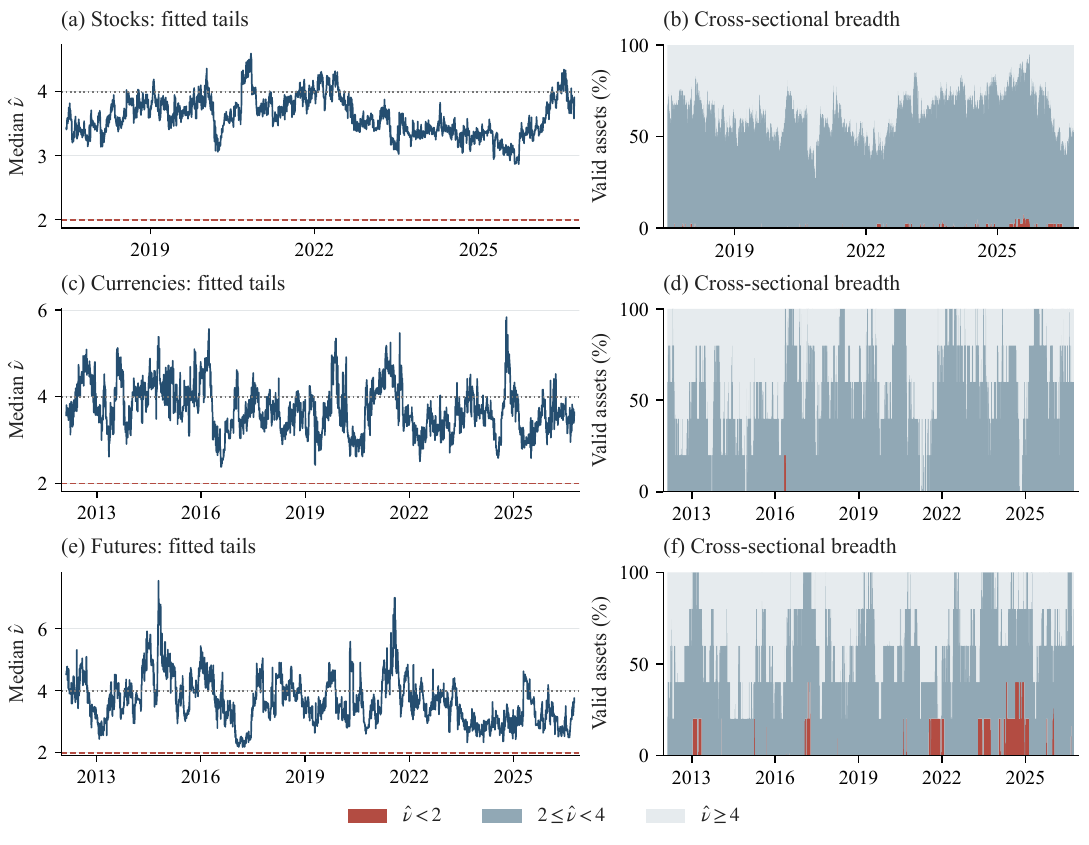}
\caption{\textbf{Fitted innovation tails and their breadth.} Left: daily cross-asset medians of bias-corrected Hill estimates, using trailing $W=120$ positions and a 10\% upper tail. Right: percentages of valid asset estimates below two, between two and four, and at least four. Every estimate uses its own trailing history.}\label{fig-dynamic-tails}
\end{figure}

The tail diagnostics describe innovations to the fitted log-RV5 system. We estimate bias-corrected Hill indices from absolute endpoint-residual deviations around their trailing median, using $W=120$ original positions and a 10\% upper-tail fraction \citep{hill1975simple,hill2010dependent,caeiro2005direct}. Figure~\ref{fig-dynamic-tails} reports the daily cross-asset median and the cross-sectional breadth of each tail regime. Under regular variation, tail indices two and four mark the thresholds for covariance existence and finite variance of quadratic inputs \citep{resnick2007heavy}. 

A single infinite innovation second moment removes the finite covariance input for the full network. The empirical counterpart of this joint requirement is the number of estimated indices below two,
\[
B_t^{(2)}=\sum_{j=1}^{K}\mathbf1\{\widehat\nu_{jt}<2\},\qquad
I_t^{(2)}=\mathbf1\{B_t^{(2)}>0\}.
\]
Table~\ref{tbl-network-tail} reports network dates rather than pooling asset--dates. At least one estimate falls below two on 18.9\% of stock dates and 19.0\% of futures dates, compared with 0.5\% of currency dates. In stocks, the corresponding asset--date frequency is 0.6\%: infrequent node-level signals can recur across a large network. Conditional on a signal, the stock network has 1.17 flagged nodes on average.
\begin{table}[!htbp]\centering\small\renewcommand{\arraystretch}{1.12}\setlength{\tabcolsep}{5pt}
\caption{A node-level tail signal is a network-level covariance diagnostic.}\label{tbl-network-tail}
\begin{tabular}{lrrrrr}\toprule
Market & $K$ & Dates & Asset--dates (\%) & Network dates (\%) & Nodes if flagged\\\midrule
Stocks & 40 & 2,333 & 0.6 & 18.9 & 1.17\\
Currencies & 5 & 3,795 & 0.1 & 0.5 & 1.00\\
Futures & 5 & 3,781 & 4.3 & 19.0 & 1.15\\\bottomrule\end{tabular}
\par\smallskip\begin{minipage}{.98\linewidth}\footnotesize\raggedright\textit{Notes.} The fourth and fifth columns report percentages with $\widehat\nu<2$: across asset--dates and across network dates with at least one such node. The last column is $E[B_t^{(2)}\mid B_t^{(2)}>0]$ as a sample average. All nodes have valid estimates on every included date. Bias-corrected Hill uses $W=120$ and a 10\% upper tail.\end{minipage}
\end{table}

For every primary point flag, the companion's 95\% block-bootstrap upper bound for the tail index reaches or exceeds two. Across the covered stock event neighborhoods, mean daily counts below four range from 19.21 to 28.32 out of forty.

\subsection{Stocks: which external dependencies receive attention}\label{sec-stock-priorities}
\begin{figure}[!htbp]\centering
\includegraphics[width=\linewidth,height=.72\textheight,keepaspectratio]{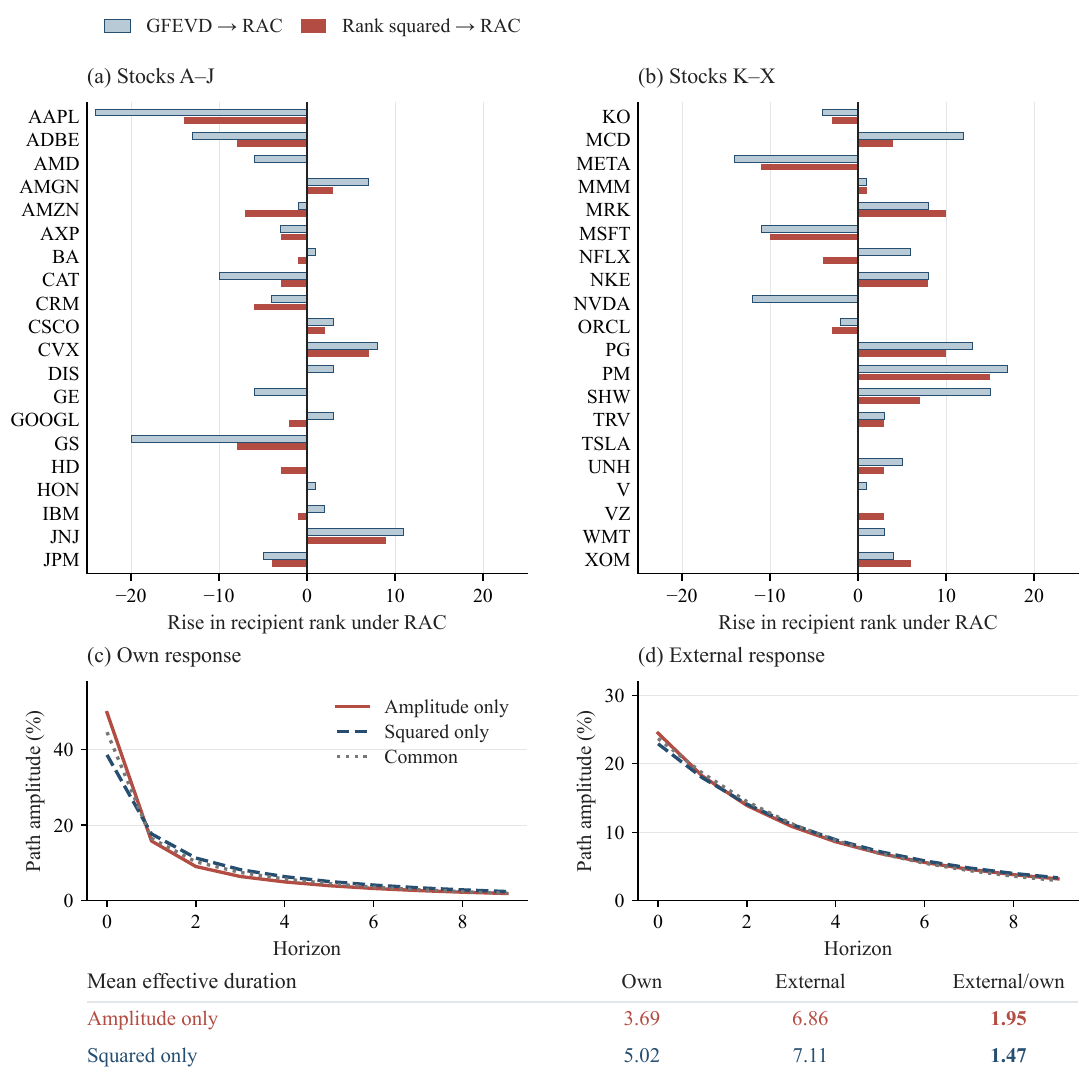}
\caption{\textbf{Stock review lists select different response profiles.} (a,b) All forty stocks in ticker order. Bars give the comparison rank minus the RAC rank, using full-sample mean incoming strength. Positive values indicate a rise under RAC. (c,d) Own and external absolute response paths for daily Top 10 groups under common rank scatter. Each recipient--date path is normalized separately to unit cumulative amplitude before averaging. The duration readout gives means of recipient--date effective durations and their external/own ratios.}\label{fig-stock-monitoring-priorities}
\end{figure}

Table~\ref{tbl-stock-top10} ranks recipients by full-sample mean external share. In the GFEVD--RAC comparison, CVX, MCD, PG and SHW enter the Top 10 in place of AAPL, AMZN, META and MSFT. The rank-squared column holds scatter fixed relative to RAC: SHW rises from 9th under squared rank responses to 2nd under amplitude. Figure~\ref{fig-stock-monitoring-priorities}a--b locates the full-sample rank changes across every stock.
\begin{table}[!htbp]\centering\small\renewcommand{\arraystretch}{1.05}\setlength{\tabcolsep}{4pt}
\caption{Stock recipients ranked by mean external share.}\label{tbl-stock-top10}
\begin{tabular}{r l r l r l r}\toprule
& \multicolumn{2}{c}{GFEVD} & \multicolumn{2}{c}{Rank squared} & \multicolumn{2}{c}{RAC}\\\cmidrule(lr){2-3}\cmidrule(lr){4-5}\cmidrule(lr){6-7} Rank & Stock & Share (\%) & Stock & Share (\%) & Stock & Share (\%)\\\midrule
1 & MSFT & 88.69 & V & 88.35 & V & 95.69\\
2 & V & 88.37 & \textbf{MSFT} & 88.10 & SHW & 95.62\\
3 & AXP & 88.30 & AXP & 88.07 & CVX & 95.61\\
4 & JPM & 88.28 & \textbf{AMZN} & 87.85 & CSCO & 95.55\\
5 & META & 88.13 & JPM & 87.45 & \textbf{XOM} & 95.54\\
6 & AAPL & 87.71 & CSCO & 87.33 & AXP & 95.53\\
7 & CSCO & 87.67 & HON & 87.27 & HON & 95.50\\
8 & HON & 87.58 & \textbf{META} & 87.16 & \textbf{PG} & 95.48\\
9 & XOM & 87.50 & SHW & 87.13 & JPM & 95.45\\
10 & AMZN & 87.39 & CVX & 86.98 & \textbf{MCD} & 95.43\\\bottomrule\end{tabular}
\par\smallskip\begin{minipage}{.98\linewidth}\footnotesize\raggedright\textit{Notes.} Ranks order full-sample mean external percentages over all 2,452 dates. Rank squared is $C(\Omega^\dagger,2)$; RAC is $C(\Omega^\dagger,1)$. Fitted dynamics are common. Bold names appear in only one of the two common-rank-scatter lists.\end{minipage}
\end{table}

\FloatBarrier

For daily lists under common rank scatter, changing attribution moves recipient ranks by 5.05 places on average across all forty stocks and replaces 3.14 names in the Top 10.

Recipients selected by amplitude but not by the squared rule adjust more quickly to their own shocks. Their mean own-response duration is 3.69 horizons, compared with 5.02 for squared-only recipients; external responses occupy 6.86 and 7.11 horizons, respectively. Relative to their own adjustment, external responses therefore persist longer in the amplitude-only group: the mean external/own duration ratio is 1.95 versus 1.47. This ordering holds in all five periods (Table~\ref{tbl-stock-monitoring-periods}).

Under common rank scatter, the daily lists share 68.6\% of their names, with a time-average HAC interval of [65.3\%,71.8\%].\footnote{Under the companion's mean and limit conditions, Bartlett HAC intervals summarize time-average fitted comparisons. They account for overlapping-window dependence using lag 250, with sensitivity lags 125 and 500 \citep{newey1987simple}.} Overlap ranges from 62.4\% to 72.7\% across five consecutive periods.

At a ten-name review budget, the squared-selected list forgoes 0.087 percentage points of mean external share under the amplitude objective; the amplitude-selected list forgoes 0.479 points under the squared objective.

\begin{table}[!htbp]\centering\small\renewcommand{\arraystretch}{1.05}\setlength{\tabcolsep}{4pt}
\caption{Stock recipient selection across five consecutive periods.}\label{tbl-stock-monitoring-periods}
\begin{tabular}{lrrrrr}\toprule
& & \multicolumn{2}{c}{Daily Top 10 lists} & \multicolumn{2}{c}{External/own duration}\\\cmidrule(lr){3-4}\cmidrule(lr){5-6} Period & Dates & Overlap (\%) & Replaced & Amplitude only & Squared only\\\midrule
2016--2018 & 491 & 72.7 & 2.73 & 1.93 & 1.48\\
2018--2020 & 491 & 62.4 & 3.76 & 1.77 & 1.31\\
2020--2022 & 490 & 69.4 & 3.06 & 1.77 & 1.38\\
2022--2024 & 490 & 69.4 & 3.06 & 2.00 & 1.55\\
2024--2026 & 490 & 68.7 & 3.13 & 2.30 & 1.64\\\bottomrule\end{tabular}
\par\smallskip\begin{minipage}{.98\linewidth}\footnotesize\raggedright\textit{Notes.} Both allocations use common rank scatter. The last two columns give mean external/own effective-duration ratios for recipients selected exclusively by that rule in the daily Top 10. Periods partition every evaluation date into nearly equal consecutive blocks.\end{minipage}
\end{table}

\paragraph{External evidence from global banks.}
The recipient-list differences extend to the public panel of 96 banks in \citet{demirer2018estimating}. Applying the same estimator at 103 monthly endpoints, GFEVD and RAC select Top 24 lists covering one quarter of the panel and differing by 8.10 names on average. Holding dynamics and rank scatter fixed still replaces 7.22 names, with a time-average HAC interval of [6.29,8.16]. Source classification also changes: 14.5\% of bank--month observations are own-dominated under squared responses and externally dominated under amplitude. See Electronic Companion Section~EC.5.3 for the bank comparisons.

\subsection{Futures: reversing the dominant source}\label{sec-external-monitoring}
The five futures contracts are mainly own-driven under the squared account and mainly externally driven under the cumulative account. With dynamics and rank scatter fixed, mean external share rises from 28.0\% to 55.9\%, and all five contracts cross the 50\% boundary in their mean shares.

\begin{figure}[!htbp]\centering
\includegraphics[width=\linewidth]{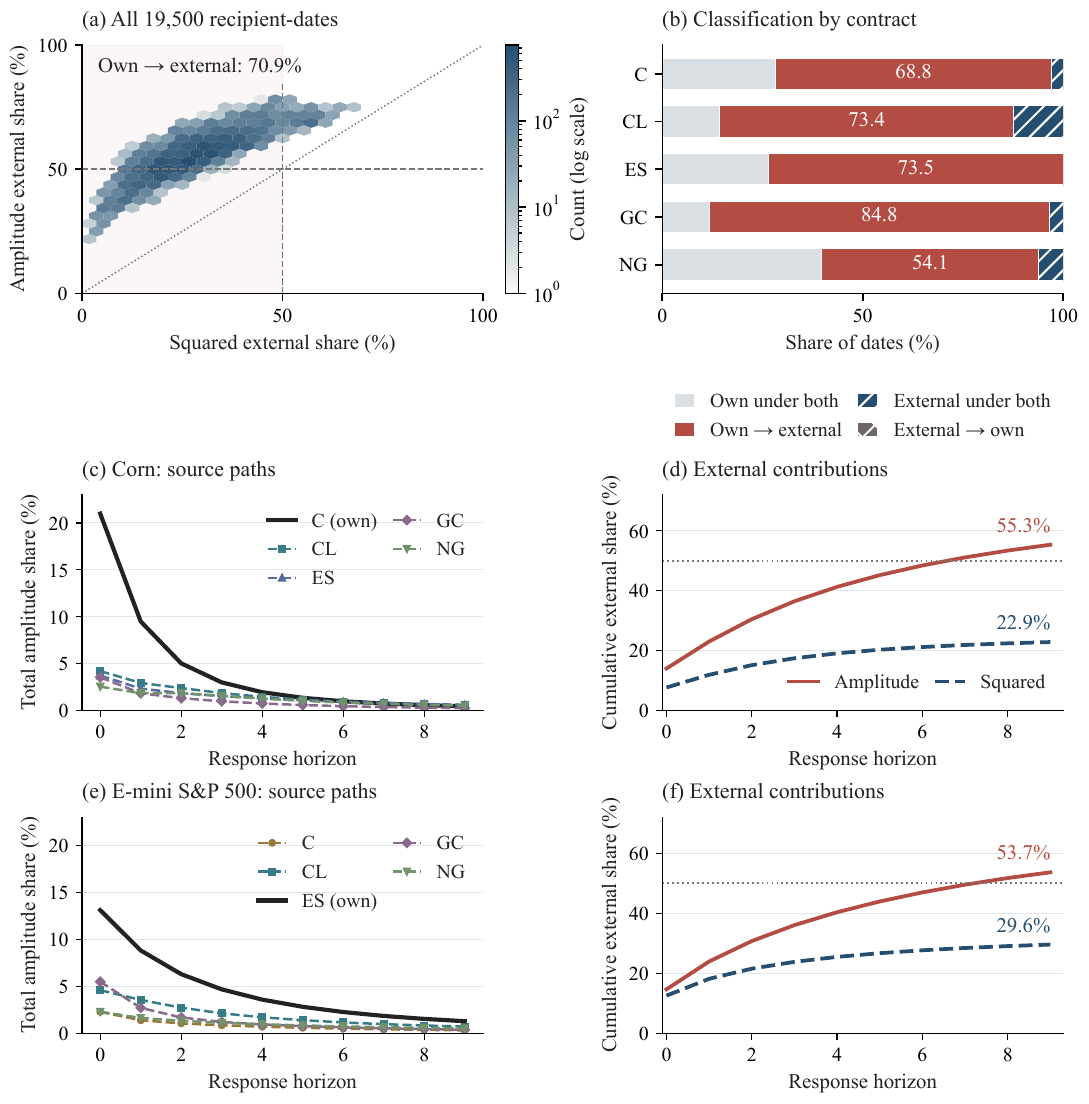}
\caption{\textbf{The same response system gives opposite dominant-source diagnoses.} Dynamics and rank scatter are fixed. (a) Joint external shares for all 19,500 recipient--dates; hexagon colors give counts on a logarithmic scale. (b) Exhaustive source classifications by contract. (c,e) Source contributions by horizon, normalized by total cumulative amplitude for each recipient--date. (d,f) Cumulative external contributions under each rule, divided by their respective full source--horizon totals. Normalization precedes averaging. Dotted references mark 50\%.}\label{fig-futures-monitoring}
\end{figure}

The switch occurs on 13,828 of 19,500 recipient--dates (70.9\%): $F_1>1/2>F_2$. Both shares lie more than five percentage points beyond the dominance boundary on their respective sides in 46.4\% of all observations. Figure~\ref{fig-futures-monitoring}a--b reports the full joint distribution and each contract's classification frequencies. For surveillance of cumulative deviations, these contracts call for joint review of external sources even when the squared account assigns the majority to their own response.

The responses after impact help explain the change. At impact, reversals occur on 41.6\% of recipient--dates; 51.8\% of the ten-horizon reversals are absent at impact. Corn illustrates the duration channel: external responses occupy 6.32 effective horizons, against 3.60 for its own response. External share rises from 22.9\% to 55.3\%, directing cumulative-deviation surveillance toward those persistent external responses.

E-mini futures illustrate source breadth. External and own durations are similar, at 6.49 and 6.35 horizons. Jointly, however, external channels account for 53.7\% of amplitude versus 29.6\% of squared response. The source diagnosis changes through their combined contribution.

The counterfactual decomposition assigns ten-horizon reversals to four cases: breadth alone (59.7\%), relative duration alone (0.5\%), either factor sufficient on its own (22.5\%), and both required together (17.2\%). These comparisons hold each measured external/own amplitude ratio fixed in Equation~\eqref{eq-external-dominance}. Table~EC.16 reports the contract-level decomposition. Table~\ref{tbl-monitoring-margins} gives reversal frequencies at four opposed-margin thresholds.
\begin{table}[!htbp]\centering\small\renewcommand{\arraystretch}{1.05}\setlength{\tabcolsep}{4pt}
\caption{Futures source-dominance reversals at increasing opposed margins.}\label{tbl-monitoring-margins}
\begin{tabular}{lr*{4}{>{\raggedleft\arraybackslash}p{.09\linewidth}}}\toprule
& & \multicolumn{4}{c}{Reversal frequency (\%)}\\\cmidrule(lr){3-6} Contract & Recipient--dates & $>0$ & $>2.5$ & $>5$ & $>10$\\\midrule
All & 19500 & 70.9 & 59.5 & 46.4 & 18.3\\
C & 3900 & 68.8 & 56.0 & 43.9 & 23.2\\
CL & 3900 & 73.4 & 64.1 & 51.7 & 17.7\\
ES & 3900 & 73.5 & 60.0 & 38.5 & 5.9\\
GC & 3900 & 84.8 & 78.9 & 71.3 & 33.3\\
NG & 3900 & 54.1 & 38.4 & 26.3 & 11.4\\\bottomrule\end{tabular}
\par\smallskip\begin{minipage}{.98\linewidth}\footnotesize\raggedright\textit{Notes.} Entries are percentages of all recipient--dates with $\min(100F_1-50,50-100F_2)$ above the column threshold, in percentage points. $F_1,F_2$ are external amplitude and squared shares under common rank scatter.\end{minipage}
\end{table}

\subsection{From recipient priorities to network structure}\label{sec-network}
In the GFEVD--RAC comparison, which changes scatter and attribution together, rank correlation is higher for transmitters than recipients across the three VOLARE panels and the global banks (Table~\ref{tbl-current-ranks}). Mean TCIs under GFEVD and RAC, respectively, are 85.17\% and 95.10\% in stocks, 60.31\% and 71.45\% in currencies, and 31.12\% and 55.92\% in futures.

\begin{table}[!htbp]\centering\small\renewcommand{\arraystretch}{1.05}\setlength{\tabcolsep}{4pt}
\caption{Transmitter and recipient agreement between GFEVD and RAC.}\label{tbl-current-ranks}
\begin{tabular}{lrrrrr}\toprule
& & \multicolumn{2}{c}{Rank correlation} & \multicolumn{2}{c}{List overlap (\%)}\\\cmidrule(lr){3-4}\cmidrule(lr){5-6} Panel & Budget & Transmitters & Recipients & Transmitters & Recipients\\\midrule
Stocks (daily) & 10 & 0.824 & 0.714 & 70.5 & 60.4\\
Banks (monthly) & 24 & 0.909 & 0.812 & 81.2 & 66.3\\
Currencies (daily) & 2 & 0.746 & 0.712 & 73.7 & 76.5\\
Futures (daily) & 2 & 0.836 & 0.534 & 86.4 & 65.9\\\bottomrule\end{tabular}
\par\smallskip\begin{minipage}{.98\linewidth}\footnotesize\raggedright\textit{Notes.} Spearman correlations use all assets; overlaps are percentages of the selected sets shared by the two allocation rules. Budgets are ten stocks, 24 banks and two assets in the five-node panels. The VOLARE panels use all 10,266 daily evaluation dates; banks use 103 monthly endpoints. Each comparison holds fitted dynamics fixed.\end{minipage}
\end{table}

Figure~\ref{fig-full-market-networks} places the stock recipient changes in the complete network. MSFT remains the largest net transmitter under both allocations, even as the recipient lists in Table~\ref{tbl-stock-top10} change. XOM$\to$CVX remains the largest individual external share, at 7.34\% under GFEVD and 3.37\% under RAC.

\begin{figure}[p]\centering
\includegraphics[width=\linewidth,height=.83\textheight,keepaspectratio]{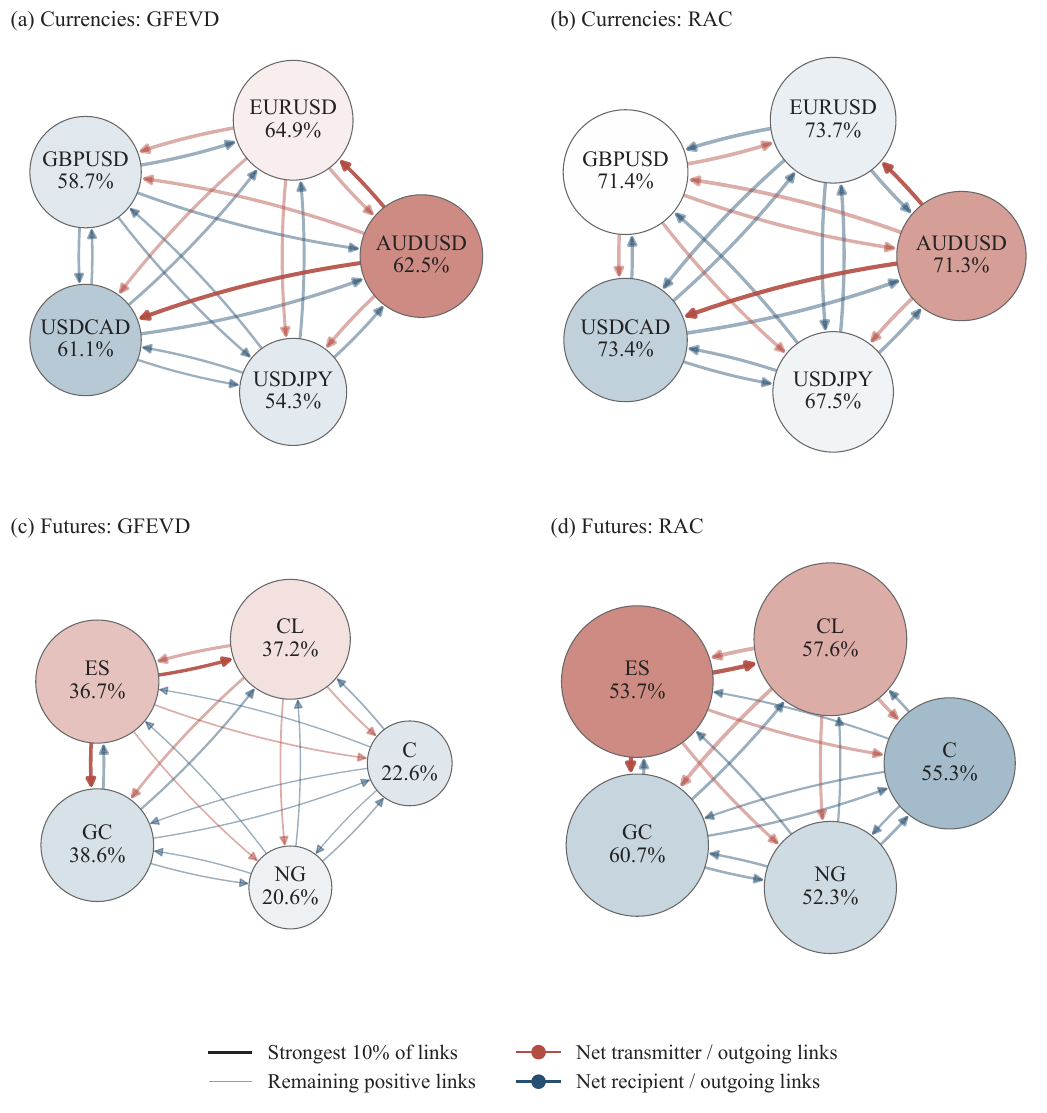}
\caption{\textbf{Complete currency and futures networks under GFEVD and RAC.} Top: currencies; bottom: futures. Each panel uses the full-sample mean allocation. Labels give each recipient's external share; node area represents incoming plus outgoing strength. Blue denotes a net recipient and red a net transmitter. Arrows follow the source's net-role color and point toward the recipient; width increases with the allocation share, and the strongest 10\% are darker. Within each market, both panels use the same node positions and the same scales for size, color and arrow width. Self-links are omitted. Futures are corn (C), crude oil (CL), E-mini S\&P 500 (ES), gold (GC) and natural gas (NG).}\label{fig-full-small-networks}
\end{figure}

The strongest bilateral channels can stay the same while the combined external share crosses the dominance boundary. In the full-sample mean networks of Figure~\ref{fig-full-small-networks}, all five currency pairs remain externally dominated, and AUD/USD$\to$USD/CAD and AUD/USD$\to$EUR/USD remain the strongest channels under both rules. Futures retain ES$\to$CL and ES$\to$GC as their strongest channels, yet every contract switches from mean own dominance under GFEVD to mean external dominance under RAC. The aggregate diagnosis credits all external sources jointly.

\paragraph{External evidence from industrial production.}
The public replication panel of \citet{brunner2023estimating} carries the comparison from volatility to production growth across 88 industrial sectors. At $H=3$, common fitted dynamics and rank scatter yield own-to-external reversals in 5,965 of 6,512 sector--month observations (91.6\%). External amplitude spans 67.63 effective source channels on average, whereas own and external effective durations are 1.53 and 1.94 horizons. These many external contributions jointly exceed the own response under the cumulative account, extending the source-breadth mechanism to production-growth networks. Electronic Companion Section~EC.6 reports the full comparison.

\begin{figure}[H]\centering
\includegraphics[width=\linewidth]{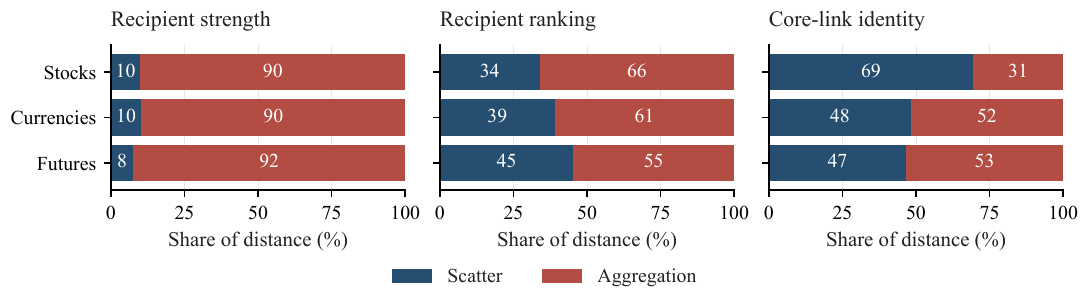}
\caption{\textbf{Scatter and aggregation contributions to allocation differences.} Shapley contributions to incoming strength, incoming ranks and the strongest 10\% of external channels are summed over daily networks before normalization. GFEVD is the reference. The companion also reports the RAC reference.}\label{fig-aggregation-direct}
\end{figure}

\paragraph{Scatter and aggregation.}\label{sec-application-support}\label{sec-sensitivity}\label{sec-average-uncertainty}
Shapley allocation separates scatter and aggregation by averaging their contributions over both orders of change \citep{shorrocks2013decomposition}. Aggregation accounts for 90.0\% of stock incoming-strength distance and 66.1\% of incoming-rank distance under the GFEVD reference; scatter accounts for 69.4\% of stock core-link distance. Figure~\ref{fig-aggregation-direct} extends the four-way comparison to all three markets.

\Needspace{.5\textheight}

\paragraph{Sensitivity.}
On common evaluation dates, changing the current window to 125 or 500 positions and the historical cap to 500 or 2,000 gives 2.86--4.24 stock Top 10 replacements. Futures source-dominance reversals range from 70.6\% to 73.3\%. Every specification retains the same RAC construction and initial tuning. Table~EC.19 reports the matched-date sensitivity comparisons.

Adding a second lag, with the same estimation procedure and all original evaluation dates, gives 3.01 stock Top 10 replacements and 68.3\% futures source-dominance reversals (Table~EC.20).

Stock recipient rankings diverge further at longer horizons: incoming-rank correlation is 0.760 at $H=5$, 0.714 at $H=10$ and 0.684 at $H=20$. Table~EC.18 reports the horizon comparisons across markets.

\FloatBarrier
\begin{figure}[p]\centering
\includegraphics[page=1,height=.88\textheight,width=\linewidth,keepaspectratio]{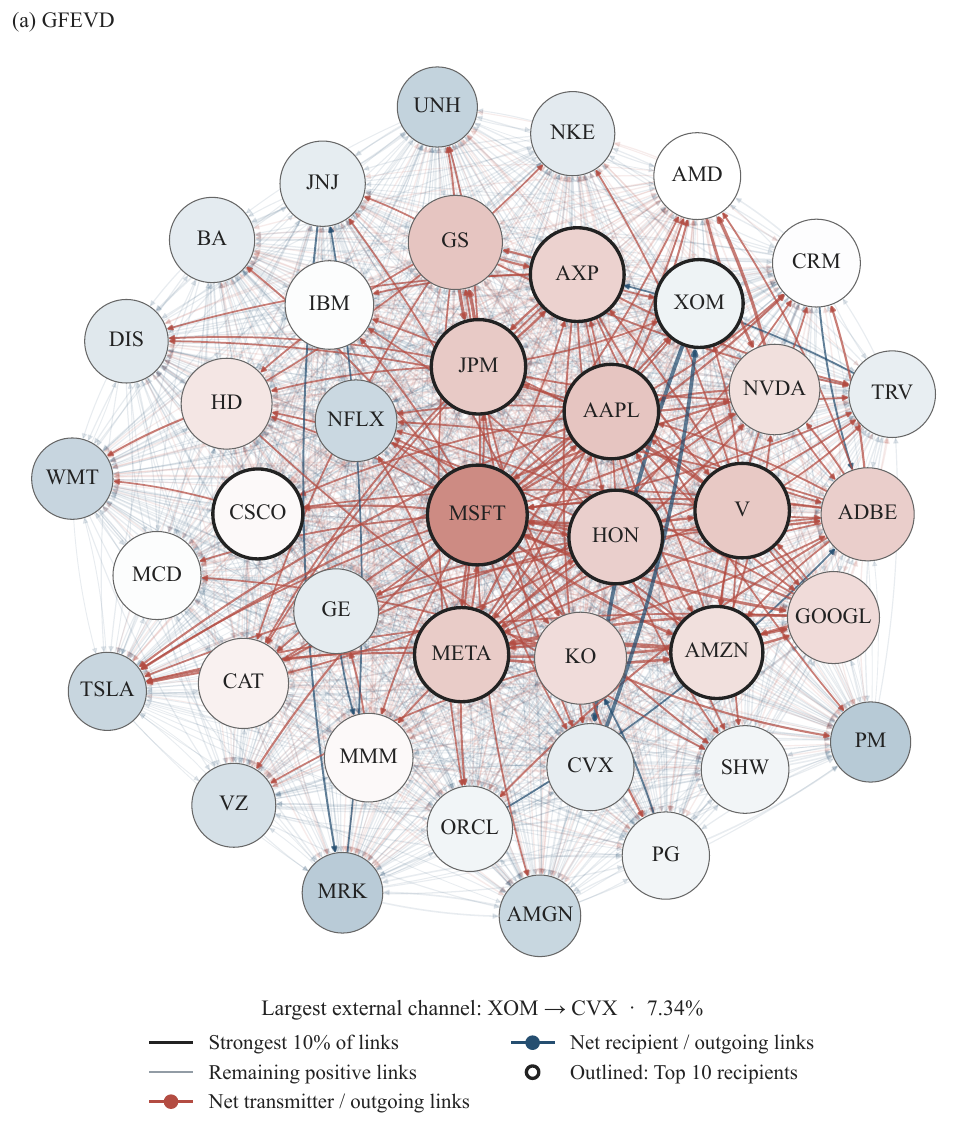}
\caption{\textbf{Complete stock networks under GFEVD and RAC.} Each panel uses the full-sample mean allocation. Node area represents incoming plus outgoing strength; blue denotes a net recipient and red a net transmitter. Arrows run from source to recipient; their color follows the source's net role and width increases with the allocation share. The strongest 10\% of external channels are darker; outlined nodes are the Top 10 recipients ranked by mean external share. Both panels use the same node positions and the same scales for size, color and arrow width. Self-links are omitted.}\label{fig-full-market-networks}\end{figure}

\begin{figure}[p]\ContinuedFloat\centering
\includegraphics[page=2,height=.88\textheight,width=\linewidth,keepaspectratio]{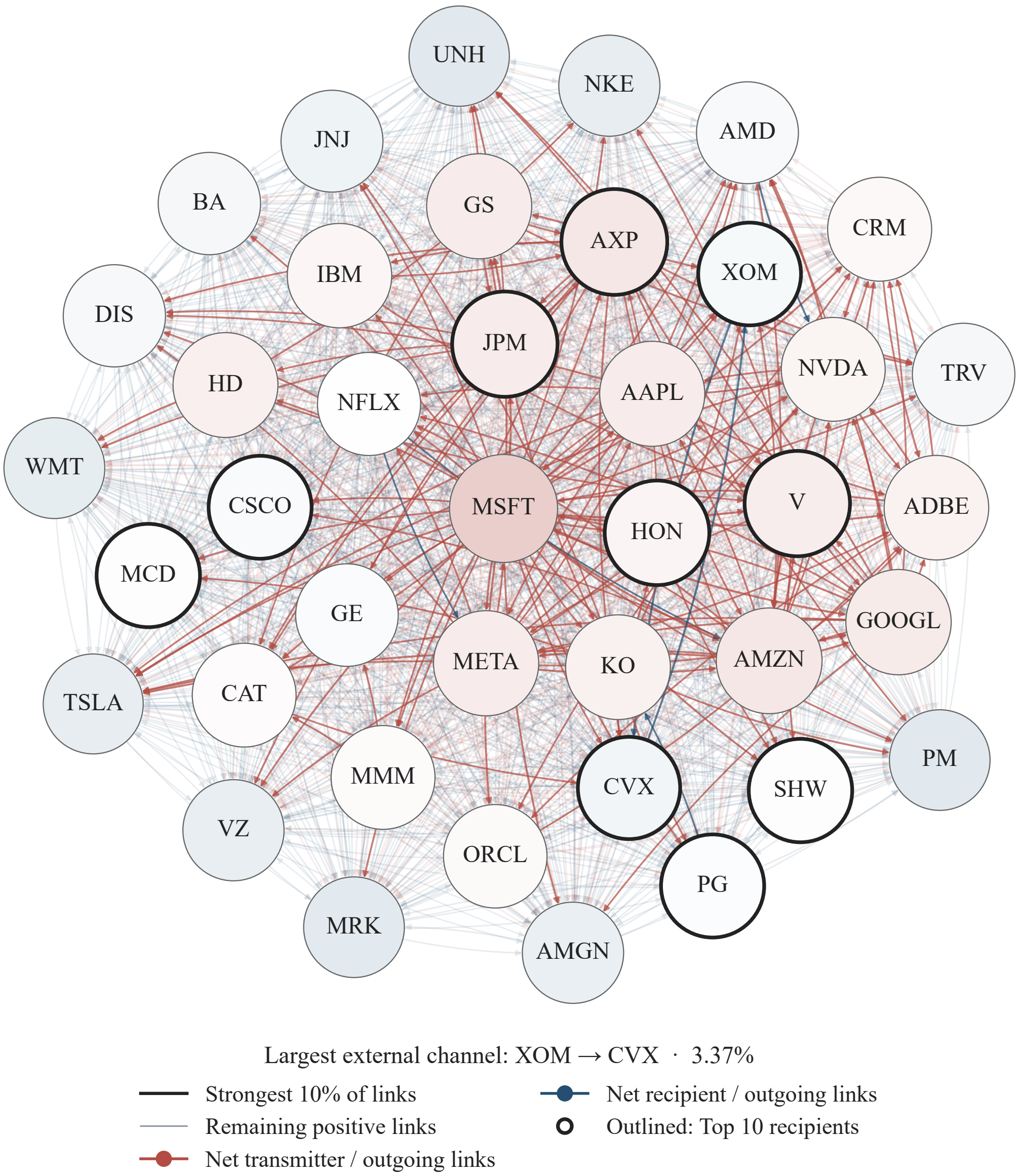}
\caption[]{\textbf{Complete stock networks under GFEVD and RAC (continued).} RAC uses the same node positions and the same scales for size, color and arrow width as GFEVD on the preceding page.}
\end{figure}

\FloatBarrier
\section{Conclusion}\label{sec-discussion}
Propagation alone does not define a financial connectedness network. The attribution rule determines which feature of its responses receives credit: even in a fully known Gaussian system, the source responsible for the larger cumulative departure can receive the smaller variance share. Energy and cumulative deviation are distinct accounting objects, with source breadth and relative duration determining when they give opposite priorities.

Before those priorities can be estimated, the allocation must have a population foundation. If even one innovation has an infinite second moment, the full network lacks a finite covariance input. Normalization yields a fixed limit only when response geometry cancels the extreme fluctuations. Quantile scales and rank dependence define a cumulative-deviation target without requiring finite first or second moments, with consistent recovery under the stated local conditions.

The stock and futures results distinguish two monitoring choices: whom to examine and whether external sources warrant joint review. Agreement on the strongest bilateral channels leaves the second question open, because collective external dependence includes all sources. Attribution is part of the network's economic definition: the monitoring objective helps determine which dependencies the network brings into view.

\Needspace{10\baselineskip}
\section*{Data availability statement}
The realized-variance data are from VOLARE \citep{cipollini2026volare}, available at \texttt{volare.unime.it} subject to its access conditions. The global-bank panel of \citet{demirer2018estimating} is publicly available from the \href{http://qed.econ.queensu.ca/jae/2018-v33.1/demirer-et-al/}{JAE Data Archive}. The industrial-production panel is available in the public replication package of \citet{brunner2023estimating}. Monte Carlo experiments use simulated data. Replication materials contain estimation and analysis code, deterministic seed rules, simulation outputs and data hashes.

\section*{Disclosure statement}
The authors declare no competing interests.

\FloatBarrier

\clearpage
\appendix
\setcounter{section}{0}
\renewcommand{\thesection}{EC.\arabic{section}}
\renewcommand{\theHsection}{ec.\arabic{section}}
\setcounter{equation}{0}
\renewcommand{\theequation}{EC.\arabic{equation}}
\renewcommand{\theHequation}{ec.\arabic{equation}}
\setcounter{theorem}{0}
\renewcommand{\thetheorem}{EC.\arabic{theorem}}
\renewcommand{\theHtheorem}{ec.\arabic{theorem}}
\setcounter{lemma}{0}
\renewcommand{\thelemma}{EC.\arabic{lemma}}
\renewcommand{\theHlemma}{ec.\arabic{lemma}}
\setcounter{proposition}{0}
\renewcommand{\theproposition}{EC.\arabic{proposition}}
\renewcommand{\theHproposition}{ec.\arabic{proposition}}
\setcounter{corollary}{0}
\renewcommand{\thecorollary}{EC.\arabic{corollary}}
\renewcommand{\theHcorollary}{ec.\arabic{corollary}}
\setcounter{definition}{0}
\renewcommand{\thedefinition}{EC.\arabic{definition}}
\renewcommand{\theHdefinition}{ec.\arabic{definition}}
\setcounter{assumption}{0}
\renewcommand{\theassumption}{EC.\arabic{assumption}}
\renewcommand{\theHassumption}{ec.\arabic{assumption}}
\setcounter{remark}{0}
\renewcommand{\theremark}{EC.\arabic{remark}}
\renewcommand{\theHremark}{ec.\arabic{remark}}
\setcounter{example}{0}
\renewcommand{\theexample}{EC.\arabic{example}}
\renewcommand{\theHexample}{ec.\arabic{example}}
\setcounter{figure}{0}
\renewcommand{\thefigure}{EC.\arabic{figure}}
\renewcommand{\theHfigure}{ec.\arabic{figure}}
\setcounter{table}{0}
\renewcommand{\thetable}{EC.\arabic{table}}
\renewcommand{\theHtable}{ec.\arabic{table}}

\begin{center}\Large\bfseries Electronic Companion\end{center}
\medskip
Section~\ref{app-foundations} establishes the population boundary and calibration. Section~\ref{app-attribution} characterizes the accounting rules and their priority reversals. Section~\ref{app-estimation} specifies the estimator and proves consistency, from fitted responses to network priorities. Section~\ref{app-finite-sample} gives the simulation designs and population targets; Sections~\ref{app-monitoring}--\ref{app-external-industry} provide data construction and supplementary financial and industrial evidence.

\section{Population networks and the covariance boundary}\label{app-foundations}\label{app-formal-results}
\subsection{Conditional-median dynamics}\label{sec-model}\label{sec-target}
Let $y_t\in\mathbb R^K$ satisfy
\begin{equation}
y_t=\sum_{\ell=1}^p A_\ell y_{t-\ell}+\varepsilon_t,\qquad
\operatorname{Med}(\varepsilon_{it}\mid\mathcal F_{t-1})=0,\quad i=1,\ldots,K.
\label{eq:model}
\end{equation}
The specification has zero intercept; including and transforming an intercept gives positive affine invariance. Write $X_t=(y_{t-1}',\ldots,y_{t-p}')'$ and $\beta_i'=(e_i'A_1,\ldots,e_i'A_p)$. The matrices $\Phi_0=I$ and
\[
\Phi_h=\sum_{\ell=1}^{p}A_\ell\Phi_{h-\ell},\qquad \Phi_h=0\ (h<0),
\]
describe the pathwise response to an innovation perturbation in the linear additive system, with subsequent innovations held fixed. The conditional-median restriction specifies the one-step dynamics.

\subsection{Response allocations}\label{sec-dy}
Fix a response horizon $H$. For a shock scatter $S$ with positive diagonal and a response aggregation rule $\rho$, write $x_{ijh}(S)=e_i'\Phi_h S e_j/\sqrt{s_{jj}}$. The coordinate vectors select the recipient and source. The resulting source share is
\begin{equation}
C^{(H,\rho)}_{ij}(S)=\frac{\sum_{h=0}^{H-1}\rho\!\left(x_{ijh}(S)\right)}{\sum_{k=1}^K\sum_{h=0}^{H-1}\rho\!\left(x_{ikh}(S)\right)}.
\label{eq:csa}
\end{equation}
Scatter $S$ fixes the generalized shock direction and scale, while $\rho_r(x)=|x|^r$ aggregates the resulting responses. Every fixed $r>0$ applies to finite population responses without an additional moment condition. The benchmarks are covariance--squared connectedness based on generalized forecast error variance decomposition (GFEVD) and rank-amplitude connectedness (RAC).

The GFEVD shares are the $C(\Sigma,\rho_2)$ instance of generalized response allocation \citep{pesaran1998generalized}. With $\Sigma=\operatorname{Cov}(\varepsilon_t)$ of positive diagonal, the generalized variance share and its row normalization are
\begin{align}
\theta^g_{ij}(H)&=
\frac{\sigma_{jj}^{-1}\sum_{h=0}^{H-1}(e_i'\Phi_h\Sigma e_j)^2}
{\sum_{h=0}^{H-1}e_i'\Phi_h\Sigma\Phi_h'e_i},\nonumber\\
C^{\mathrm{GFEVD}}_{ij}(H)&=
\frac{\theta^g_{ij}(H)}{\sum_k\theta^g_{ik}(H)}
=\frac{\sum_{h=0}^{H-1}|x_{ijh}(\Sigma)|^2}
{\sum_k\sum_{h=0}^{H-1}|x_{ikh}(\Sigma)|^2}.
\label{eq:dy}
\end{align}
The mean-based forecast interpretation requires the corresponding moments and dynamic specification. Generalized shocks may be correlated; structural attribution requires additional identification.

\subsection{The covariance boundary}\label{sec-validity}
The covariance requirement applies jointly to all $K$ innovation coordinates. A finite covariance matrix exists exactly when every coordinate has a finite second moment. Indeed,
\[
E\|\varepsilon_t\|_2^2=\sum_{j=1}^K E\varepsilon_{jt}^2,
\qquad
E|\varepsilon_{it}\varepsilon_{jt}|
\le (E\varepsilon_{it}^2\,E\varepsilon_{jt}^2)^{1/2}.
\]
Thus finite marginal second moments give finite means and cross-products, whereas one infinite marginal second moment excludes a finite covariance matrix for the full network. For heterogeneous regularly varying marginals, $\min_j\nu_j>2$ is sufficient for this moment requirement; any $\nu_j<2$ violates it. At index two, the slowly varying tail factor decides the boundary moment. Convergence of normalized sample shares is a separate property of the allocation map.

\begin{proposition}[Moment restrictions and their scope]\label{prop:moment}
Suppose the innovation norm is regularly varying with index $\nu>0$ and has a nontrivial tail direction.
For $\nu<2$, at least one second moment is infinite and the full covariance input for \eqref{eq:dy} is unavailable.
For $\nu>2$, second moments are finite.
For $2<\nu<4$, covariance exists but at least one diagonal quadratic term has infinite variance, precluding the usual iid finite-variance covariance CLT.
For $\nu>4$, fourth moments are finite; covariance and dynamic-estimator CLTs still require their own dependence and identification conditions.
At $\nu=2$ or $4$, the slowly varying tail factor determines moment existence. For Student-$t$ innovations, second and fourth moments exist exactly when $\nu>2$ and $\nu>4$, respectively.
\end{proposition}

\subsubsection{From moments to normalized shares}\label{app-moments}
Let $R=\|\varepsilon\|$ have survival function
$P(R>x)=x^{-\nu}L(x)$ with slowly varying $L$. For $q>0$,
\[
E R^q=q\int_0^\infty x^{q-1}P(R>x)\,dx.
\]
Tail integration gives finiteness for $q<\nu$ and divergence for $q>\nu$. At $q=\nu$, the relevant integral is proportional to
$\int^\infty L(x)\,dx/x$; regular variation alone does not decide it. Since $R^2=\sum_i\varepsilon_i^2$, divergence of $E R^2$ makes at least one coordinate second moment infinite. The same conclusion for fourth moments follows from equivalence of finite-dimensional norms. This proves Proposition~\ref{prop:moment}. Student-$t$ tails have a nonzero asymptotically constant slowly varying factor, so the boundary moments diverge.

Normalization determines whether covariance fluctuations survive in a share. With iid finite-fourth-moment innovations, covariance entries satisfy a standard multivariate CLT; dependent innovations require the corresponding dependence CLT. For $2<\nu<4$, at least one diagonal quadratic input has infinite variance. Proposition~\ref{prop:dy-share-rate} transfers an explicit stable-input limit through the allocation derivative under a noncancellation condition.

For iid $0<\nu<2$, choose the tail quantile $b_T$ by
$T P(R>b_T)\to1$. Under multivariate domain-of-attraction conditions, point-process convergence of extremes implies
\[
b_T^{-2}\sum_t\varepsilon_t\varepsilon_t'\Rightarrow S_\nu,
\]
a positive-semidefinite stable random matrix with index $\nu/2$ \citep{resnick2007heavy,davis1998sample}. For the sample covariance, this gives the main paper's scaling $c_T=T/b_T^2$. The reciprocal scale $b_T^2/T$ combines the power $T^{2/\nu-1}$ with a slowly varying factor. Centering by the sample mean has negligible effect on this normalized iid second-product limit.

Because $C(cS,2)=C(S,2)$ for $c>0$, a continuous-mapping argument gives $C(\widehat\Sigma_T,2)\Rightarrow C(S_\nu,2)$ where the required diagonals and denominators are positive. When fitted dynamics are used, their joint convergence must also be supplied. Nondegeneracy of this ratio limit requires $C(\cdot,2)$ to be nonconstant on the limiting support.

\begin{proposition}[A deterministic ratio without covariance]\label{prop:det-ratio}
Let $\varepsilon_t=Z_t v$, where $Z_t$ are iid symmetric continuous regularly varying scalars of index below two, and all $v_j$ are nonzero. With $A=0$ and $H=1$, all covariance-based sample directional shares equal $1/K$, although the population covariance does not exist.
\end{proposition}
\begin{proof}
The sample covariance is $q_T vv'$ for a positive scalar $q_T$. Therefore
\[
|x_{ij0}|^2=(q_T v_i v_j)^2/(q_T v_j^2)=q_T v_i^2,
\] independent of $j$. Row normalization gives $1/K$. The full covariance has infinite diagonal moments.
\end{proof}
Proposition~\ref{prop:ratio-support} gives support conditions for nondegenerate share limits and covariance-rate fluctuations.

With finite covariance, the target is fixed. The scaling below magnifies covariance estimation error around its finite limit; the earlier scaling applies to the covariance level itself.

\begin{proposition}[GFEVD share rate and cancellation]\label{prop:dy-share-rate}
Let $2<\nu<4$ and let the mean-zero innovations be iid regularly varying with index $\nu$. Assume the centered sample covariance satisfies $T^{(\nu-2)/\nu}(\widehat\Sigma_T-\Sigma)\Rightarrow S_\nu$ for a nondegenerate $\nu/2$-stable random element $S_\nu$ of the space of symmetric matrices \citep{davisresnick1986,resnick2007heavy}. The centered stable limit may be indefinite. The displayed power normalization is an additional assumption; general regular variation can require a slowly varying normalization. Fix $H$ and a stable propagation $\{\Phi_h\}$, and write the GFEVD share as $C_{ij}(\Sigma)=N_{ij}(\Sigma)/D_i(\Sigma)$ with
\[
N_{ij}(\Sigma)=\sum_{h=0}^{H-1}(e_i'\Phi_h\Sigma e_j)^2/\sigma_{jj},\qquad D_i=\sum_k N_{ik}>0.
\] Assume all $\sigma_{jj}>0$. The map $\Sigma\mapsto C_{ij}(\Sigma)$ is differentiable at $\Sigma$. By the delta method applied to the stable covariance limit,
\[
T^{(\nu-2)/\nu}\{C_{ij}(\widehat\Sigma_T)-C_{ij}(\Sigma)\}
\;\Rightarrow\;
\frac{\big\langle \nabla N_{ij}-C_{ij}\textstyle\sum_k\nabla N_{ik},\;S_\nu\big\rangle}{D_i},
\]
where $\nabla N_{ij}$ is the gradient of $N_{ij}$ in $\Sigma$ and $\langle\cdot,\cdot\rangle$ is the Frobenius trace inner product. Therefore:
\begin{enumerate}
\item (Noncancelling stable limit) If $\big\langle \nabla N_{ij}-C_{ij}\sum_k\nabla N_{ik},\,S_\nu\big\rangle$ is nondegenerate, the GFEVD share inherits the covariance estimator's stable rate $T^{-(\nu-2)/\nu}$ with a non-Gaussian limit.
\item (Cancellation) If $\langle\nabla N_{ij}-C_{ij}\sum_k\nabla N_{ik},S_\nu\rangle=0$ almost surely, the leading stable term vanishes and $C_{ij}(\widehat\Sigma_T)-C_{ij}(\Sigma)=o_p(T^{-(\nu-2)/\nu})$, so the share is estimated strictly faster than the covariance input.
\end{enumerate}
The cancellation condition is a joint restriction on propagation geometry and the support of $S_\nu$: the share derivative must vanish along the fluctuating covariance directions. Cancellation removes the first-order term at the specified covariance input. The rank-one geometry of Proposition~\ref{prop:det-ratio} also gives exact cancellation with finite covariance if the scalar innovation instead has index between two and four. This proposition holds propagation fixed; using estimated propagation requires its error to be negligible at the displayed covariance rate.
\end{proposition}

\subsubsection{Deterministic share limits and surviving fluctuations}
The allocation map must be constant on the limiting scatter support for normalized shares to approach a deterministic network. Nonconstancy identifies the fluctuations that survive normalization.

Let $\mathbb S^K$ be the space of real symmetric $K\times K$ matrices, with Frobenius inner product, and let $\mathbb S_{++}^K$ be its positive-definite cone. Fix $K\ge2$, finite $H\ge1$ and finite propagation matrices $\Phi_0=I,\Phi_1,\ldots,\Phi_{H-1}$. Write $f_{ij}(S)=C^{(H,2)}_{ij}(S)$. As in Proposition~\ref{prop:dy-share-rate}, propagation is fixed.

\paragraph*{Proof of the deterministic-limit proposition in the main paper.}
\begin{proof}
Take positive-diagonal sample covariance inputs with $c_T\widehat\Sigma_T\Rightarrow S_\nu$, where $c_T>0$ is deterministic and $S_\nu$ is positive semidefinite with strictly positive diagonal almost surely. On the positive-semidefinite cone with positive diagonal, every response is continuous in scatter. The own impact term is $x_{ii0}(S)=\sqrt{s_{ii}}>0$, so the squared-response row denominator is strictly positive. The allocation map is therefore continuous on this domain and invariant to positive common scaling. Hence
\[
C(\widehat\Sigma_T,2)=C(c_T\widehat\Sigma_T,2)
\Rightarrow C(S_\nu,2).
\]
If $C(S_\nu,2)=C_0$ almost surely, convergence in distribution to the constant matrix $C_0$ is convergence in probability. Conversely, convergence in probability to $C_0$ implies convergence in distribution to that constant; uniqueness of the weak limit gives $C(S_\nu,2)=C_0$ almost surely. This proves the equivalence.

Continuity also makes almost-sure constancy equivalent to constancy on the limiting support within this domain. If the map differed from $C_0$ at a support point in the domain, a neighborhood of that point would have positive probability and keep the map separated from $C_0$, contradicting almost-sure equality. The converse follows because $S_\nu$ belongs to its support and to the domain almost surely. Lemma~\ref{lem:share-nonconstant} and the first part of Proposition~\ref{prop:ratio-support} below establish nondegeneracy of every individual share under the stated open-support condition.
\end{proof}

\begin{lemma}[Nonconstancy of individual shares]\label{lem:share-nonconstant}
Every $f_{ij}$ is real analytic and nonconstant on $\mathbb S_{++}^K$. For each fixed propagation system there is a proper algebraic exceptional set $\mathcal E$, of Lebesgue measure zero in $\mathbb S^K$, such that
\[
\nabla f_{ij}(S)\ne0\quad\text{for all }i,j
\quad\text{whenever }S\in\mathbb S_{++}^K\setminus\mathcal E.
\]
Consequently no individual share is constant on a nonempty open subset of $\mathbb S_{++}^K$.
\end{lemma}
\begin{proof}
\textit{Nonconstancy.} Put $a_{ij}=\sum_{h=0}^{H-1}(e_i'\Phi_h e_j)^2$, so $a_{ii}\ge1$. On a diagonal matrix $S=\operatorname{diag}(d_1,\ldots,d_K)$,
\[
f_{ij}(S)=\frac{a_{ij}d_j}{\sum_k a_{ik}d_k}.
\]
For $i\ne j$ with $a_{ij}>0$, varying $d_j$ changes this ratio because its denominator also contains $a_{ii}d_i>0$. If $a_{ij}=0$, the ratio is zero on diagonal inputs, whereas at $S=I+\rho(e_i e_j'+e_j e_i')$, $0<|\rho|<1$, its impact numerator is positive. Thus every off-diagonal share is nonconstant. For $i=j$, a positive $a_{ik}$ with $k\ne i$ makes the diagonal ratio vary with $d_k$. Otherwise $f_{ii}=1$ on diagonal inputs, while the same off-diagonal perturbation makes another channel positive and hence $f_{ii}<1$.

\textit{Analyticity.} Each numerator $N_{ij}(S)$ is a sum of quadratic polynomials divided by positive $s_{jj}$. Multiplying every numerator and its row denominator by $\prod_k s_{kk}$ gives polynomials $P_{ij}$ and $Q_i$, with $Q_i>0$ on $\mathbb S_{++}^K$; the own impact term ensures this positivity. Hence $f_{ij}=P_{ij}/Q_i$ is analytic there.

\textit{The exceptional set.} Use the independent upper-triangular entries as coordinates $z$. Since $f_{ij}$ is nonconstant on the connected cone, at least one polynomial
$Q_i\partial_{z_r}P_{ij}-P_{ij}\partial_{z_r}Q_i$ is nonzero. Its zero set contains the critical points of $f_{ij}$. A nonzero polynomial has a Lebesgue-null zero set: this follows by induction on dimension, viewing it as a polynomial in the last coordinate and applying Fubini outside the null set where all coefficient polynomials vanish. Taking the finite union over $i,j$ gives a proper algebraic null set $\mathcal E$. Constancy on an open subset would put that subset inside $\mathcal E$, which is impossible.
\end{proof}

\begin{proposition}[Sufficient conditions for covariance-rate transfer]\label{prop:ratio-support}
Maintain the respective covariance-input convergence assumptions in the preceding subsection.
\begin{enumerate}
\item For $0<\nu<2$, suppose the required diagonal and denominator conditions hold almost surely at $S_\nu$, and its topological support contains a nonempty open subset of $\mathbb S_{++}^K$. Then every individual limit $f_{ij}(S_\nu)$ is nondegenerate. Thus no corresponding sample share converges in probability to a deterministic constant.
\item For $2<\nu<4$, suppose the stable matrix limit is full-dimensional: no nonzero $G\in\mathbb S^K$ makes $\langle G,S_\nu\rangle$ deterministic. If $\Sigma\in\mathbb S_{++}^K\setminus\mathcal E$, every individual share has a nondegenerate stable first-order limit at the covariance rate in Proposition~\ref{prop:dy-share-rate}.
\end{enumerate}
\end{proposition}
\begin{proof}
For the first part, if $f_{ij}(S_\nu)=c$ almost surely, continuity forces $f_{ij}=c$ at every support point in its continuity domain. Otherwise a neighborhood of a support point with $f_{ij}\ne c$ would have positive probability, contradicting almost-sure equality. Constancy on the specified open subset contradicts Lemma~\ref{lem:share-nonconstant}. Convergence in probability to a constant would also imply a degenerate distributional limit.

For the second part, the derivative in Proposition~\ref{prop:dy-share-rate} is precisely $\langle\nabla f_{ij}(\Sigma),S_\nu\rangle$. Its coefficient matrix is nonzero outside $\mathcal E$, so full dimensionality excludes cancellation. A nondegenerate linear image of an $\alpha$-stable law, $\alpha=\nu/2\in(1,2)$, remains non-Gaussian $\alpha$-stable. The stated covariance-rate delta method applies.
\end{proof}

\begin{corollary}[A rank-one tail-direction criterion]\label{cor:quadratic-tail-directions}
Let $\Lambda$ be the finite angular measure of a rank-one matrix-stable input: its L\'evy jumps have form $r vv'$, $r>0$, $v$ on the unit sphere, with intensity
$r^{-1-\alpha}\,dr\,\Lambda(dv)$, $\alpha=\nu/2$. Suppose
\begin{equation}
\operatorname{span}\{vv':v\in\operatorname{supp}\Lambda\}=\mathbb S^K.
\label{eq:quadratic-tail-span}
\end{equation}
For $\alpha<1$, take the positive stable input with no deterministic drift; it is positive definite almost surely and its support contains a positive-definite open set. For $1<\alpha<2$, allow the centered stable input; it is full-dimensional. Thus \eqref{eq:quadratic-tail-span} supplies the support conditions in Proposition~\ref{prop:ratio-support}, retaining the covariance exception $\mathcal E$ in its second part.
\end{corollary}
\begin{proof}
Select finitely many directions $v_1,\ldots,v_m$ whose outer products span $\mathbb S^K$. They also span $\mathbb R^K$. The linear map
$L(w)=\sum_{\ell=1}^m w_\ell v_\ell v_\ell'$ is surjective and therefore open; $L((0,\infty)^m)$ is a nonempty open set consisting of positive-definite matrices.

For $\alpha<1$, the stable input is the convergent sum of its positive L\'evy jumps. Every $L(w)$ with $w_\ell>0$ belongs to its support. Indeed, fix an error tolerance and a small cutoff $\delta<\min_\ell w_\ell$. Prescribe jumps in disjoint neighborhoods of the finitely many $w_\ell v_\ell v_\ell'$ and no other jumps larger than $\delta$. The finite-intensity Poisson process gives this event positive probability; directions with identical outer products can first be merged. The independent remainder has expected trace
\[
E\operatorname{tr}(S_{\le\delta})
=\frac{\Lambda(\mathbb S^{K-1})}{1-\alpha}\delta^{1-\alpha},
\]
which can be made arbitrarily small. Markov's inequality then gives positive probability that its norm is below the remaining tolerance. This proves the support claim.

For almost-sure positive definiteness, choose $K$ linearly independent supported directions and small angular neighborhoods such that every choice of one direction from each neighborhood remains independent. Each neighborhood has infinite jump intensity near zero and hence contains a jump almost surely. Their positive outer-product sum is positive definite, and all other jumps are positive semidefinite.

For $1<\alpha<2$, the characteristic-function modulus of a projection is
\[
\left|E\exp\{it\langle G,S_\nu\rangle\}\right|
=\exp\!\left\{-c_\alpha|t|^\alpha
\int|v'Gv|^\alpha\,\Lambda(dv)\right\},
\]
where $c_\alpha=\int_0^\infty(1-\cos r)r^{-1-\alpha}\,dr\in(0,\infty)$. This follows by taking the real part of the L\'evy exponent and substituting $r|t v'Gv|$. A zero integral forces $v'Gv=0$ on the angular support by continuity, and \eqref{eq:quadratic-tail-span} then forces $G=0$. Every nonzero projection is therefore nondegenerate.
\end{proof}

For example, the $K(K+1)/2$ directions $e_i$ and $(e_i+e_j)/\sqrt2$, $i<j$, span all symmetric matrices through their outer products. Full angular support also satisfies \eqref{eq:quadratic-tail-span}. Independent coordinate-axis tails span the diagonal subspace. With diagonal propagation, all diagonal scatter inputs then yield $C(S,2)=I$, despite noncollinear factors.

The criterion is sufficient, with genericity defined over covariance coordinates; a fixed innovation law can lie in the exceptional set. The conclusions apply to individual shares. Row sums are deterministic, and other aggregates may also cancel the fluctuations.

\subsection{Rank scatter and population existence}\label{sec-method}\label{sec-scatter}
Population medians are lower $0.5$ quantiles. For an innovation law $P$, let $m_i$ be its marginal median,
$d_i=\operatorname{Med}(|\varepsilon_i-m_i|)$, and $s_i=1.4826d_i$.
Let $\varepsilon'$ be an independent draw from $P$. Define
\[
\tau_{ij}=E\!\left[\operatorname{sign}
\{(\varepsilon_i-\varepsilon_i')(\varepsilon_j-\varepsilon_j')\}\right],
\quad \tau_{ii}=1,\quad
\widetilde R_{ij}=\sin(\pi\tau_{ij}/2).
\]
The matrix $\widetilde R$ need not be positive semidefinite outside an elliptical model. With $D=\operatorname{diag}(s_i)$ and
$\mathcal C=\{R=R':R\succeq0,\operatorname{diag}(R)=\mathbf1\}$, set
\begin{equation}
R^\dagger=\arg\min_{R\in\mathcal C}\|R-\widetilde R\|_F,\qquad
\Omega^\dagger=D R^\dagger D.
\label{eq:scatter}
\end{equation}
The Frobenius nearest-correlation problem has a unique, nonexpansive solution \citep{higham2002computing}. Projection therefore supplies a well-defined scatter outside ellipticity and carries convergence of the concordance matrix to convergence of scatter.

For $H\ge1$ and $r\in\{1,2\}$, use the comparison notation
\begin{equation}
C^{(H,r)}_{ij}(S)=
\frac{\sum_{h=0}^{H-1}|x_{ijh}(S)|^r}
{\sum_{k=1}^K\sum_{h=0}^{H-1}|x_{ikh}(S)|^r},
\qquad \mathrm{RAC}_{ij}(H)=C^{(H,1)}_{ij}(\Omega^\dagger).
\label{eq:mrc}
\end{equation}
The index $i$ denotes recipient and $j$ source. In percentage points, outbound transmission is
$TO_j=100\sum_{i\ne j}C_{ij}$, the hub is $\arg\max_j TO_j$, and
$TCI=100K^{-1}\sum_{i\ne j}C_{ij}$. RAC assigns each source its relative contribution to cumulative absolute rank-scatter response.

\begin{theorem}[Existence and calibration of the rank-scatter implementation]\label{thm:exist}
With finite dynamic coefficients and continuous innovation marginals having unique medians and finite strictly positive MADs, RAC is finite and has positive row denominators without a moment condition. In particular, its denominator is at least
$|x_{ii0}(\Omega^\dagger)|=s_i>0$.
If innovations are continuous elliptical with positive shape matrix $\Omega^*$, the Kendall identity and common radial MAD scaling imply
$\Omega^\dagger=\kappa_F\Omega^*$ for a positive constant $\kappa_F$.
If covariance also exists, $\Sigma=c_F\Omega^*$, and
\[
C^{(H,2)}(\Omega^\dagger)=C^{(H,2)}(\Sigma)=C^{\mathrm{GFEVD}}(H).
\]
\end{theorem}

\begin{remark}[The RAC target under heavy tails]\label{rem:rac-target}
The direction $\Omega^\dagger e_j/\sqrt{\omega^\dagger_{jj}}$ combines marginal MAD scales and projected Kendall dependence. Covariance uses expected products of centered innovations; MAD measures central spread through a quantile of absolute deviations, and Kendall dependence measures pairwise concordance.
\end{remark}

\subsubsection{Projection}
The correlation set
$\mathcal C=\{R=R':R\succeq0,\operatorname{diag}(R)=\mathbf1\}$
is nonempty, closed, convex, and bounded. The bound $|R_{ij}|\le1$ follows from positive semidefinite $2\times2$ principal minors. The Frobenius projection $\Pi_{\mathcal C}$ therefore exists and is unique.

For $P_A=\Pi_{\mathcal C}(A)$ and $P_B=\Pi_{\mathcal C}(B)$, the projection variational inequalities are
$\langle A-P_A,P_B-P_A\rangle\le0$ and
$\langle B-P_B,P_A-P_B\rangle\le0$. Adding gives
\[
\|P_A-P_B\|_F^2\le
\langle P_A-P_B,A-B\rangle,
\]
hence $\|P_A-P_B\|_F\le\|A-B\|_F$.
Nonexpansiveness transfers consistency and rates of the input, including when the limiting correlation lies on the boundary of $\mathcal C$.

The algorithm follows \citet{higham2002computing}: subtract the Dykstra correction, project onto positive semidefinite symmetric matrices, update the correction, and replace the diagonal by ones. Diagonal replacement is the Frobenius projection onto the unit-diagonal affine set. Convergence is checked both by relative iterate change and minimum eigenvalue.

\subsubsection{Existence}
Set $D=\operatorname{diag}(1.4826d_i)$,
$\widetilde R_{ij}=\sin(\pi\tau_{ij}/2)$ and
$\Omega^\dagger=D\Pi_{\mathcal C}(\widetilde R)D$.
The Kendall kernel is bounded, so $\tau_{ij}$ is a finite expectation. The projection is finite and has unit diagonal, giving
$\omega^\dagger_{ii}=(1.4826d_i)^2>0$.
All finite-horizon propagation matrices and responses are finite. The own impact response is
\[
x_{ii0}(\Omega^\dagger)
=\omega^\dagger_{ii}/\sqrt{\omega^\dagger_{ii}}
=1.4826d_i>0.
\]
The own impact alone keeps the denominator positive, including when all cross-responses vanish. This proves existence without a moment or full-rank-scatter requirement.

\subsection{Elliptical calibration}
For a continuous elliptical law with positive shape $\Omega^*$,
\[
\sin(\pi\tau_{ij}/2)
=\frac{\omega^*_{ij}}{\sqrt{\omega^*_{ii}\omega^*_{jj}}}
\]
by \citet{lindskog2003kendall}. The transformed concordances already form a correlation matrix, so projection leaves them unchanged. Elliptical marginals also share a standardized law, giving MADs of the form
$d_i=d_F\sqrt{\omega^*_{ii}}$ for a common $d_F>0$. Consequently
\[
\Omega^\dagger=(1.4826d_F)^2\Omega^*.
\]
When covariance exists, $\Sigma=c_F\Omega^*$ for a radial-law constant $c_F>0$. Writing $\Sigma=c\Omega^\dagger$ gives
$x_{ijh}(\Sigma)=\sqrt c\,x_{ijh}(\Omega^\dagger)$. Row normalization cancels this common factor for every power allocation. In particular
$C(\Omega^\dagger,2)=C(\Sigma,2)$.
For Student-$t_\nu$ shape $\Omega^*$, $c_F=\nu/(\nu-2)$ when $\nu>2$; the constant depends on the radial law.

Outside ellipticity, $\Omega^\dagger$ remains defined but need not be proportional to covariance. Common-scatter comparisons isolate aggregation for any innovation law: the accounting characterization and breadth--duration identities concern fixed response paths.

\subsection{Unit invariance}\label{app-axioms}

For a positive diagonal unit transformation $B$, dynamic coefficients transform as $A_\ell^*=BA_\ell B^{-1}$, hence $\Phi_h^*=B\Phi_h B^{-1}$.
MAD scales transform by $s_i^*=B_{ii}s_i$, Kendall ranks are unchanged, and
$\Omega^{\dagger *}=B\Omega^\dagger B$. Also
$\Sigma^*=B\Sigma B$ where covariance exists.
For either scatter,
\[
q_j^*=\frac{S^*e_j}{\sqrt{s_{jj}^*}}=Bq_j,
\qquad
x_{ijh}^*=B_{ii}x_{ijh}.
\]
The factor $B_{ii}^r$ cancels in recipient-row normalization. Translation equivariance requires the correspondingly transformed intercept.

If marginal changes preserve MAD and the copula while dynamics are fixed, they preserve $\Omega^\dagger$ and RAC. Common covariance rescaling also preserves GFEVD.

\FloatBarrier

\section{Attribution dependence and external dependence}\label{app-attribution}\label{app-existence}\label{app-ratios}

\subsection{Cumulative-response accounting}\label{app-cumulative-accounting}
Compare two paths of the same linear system, starting from the same initial conditions and receiving the same subsequent innovations. At impact, add $S e_j/\sqrt{s_{jj}}$ to one path. With propagation matrices, scatter $S$ and horizon $H$ fixed, the difference at recipient $i$ and horizon $h$ is $x_{ijh}(S)$. For $y_i=\log RV_i$, it is the log ratio of the two modeled variance paths. Summing its absolute magnitude gives $a_{ij}=\sum_{h=0}^{H-1}|x_{ijh}(S)|$, the accumulated separation over the specified horizon. Evaluating each source perturbation separately and normalizing the resulting totals gives the recipient's source shares.

The two accounting objectives attach different penalties to the same path separation. A constant penalty per unit of absolute separation gives cumulative-deviation loss; a quadratic penalty gives response-energy loss. Both use the same rule at every date and treat signs symmetrically. Cumulative accounting preserves same-sign splits: spreading a fixed total departure across periods leaves its score unchanged. Under quadratic accounting, doubling a departure quadruples its score. The characterization is defined on finite response paths, so its accounting distinction is independent of the particular dynamic model that generates them.

\paragraph*{Proof of the response-accounting characterization in the main paper.}
\begin{proof}
Let $V$ be a nonnegative score on all finite real response paths. It is invariant to changing any response's sign, additive under concatenation, and satisfies $V(0)=0$, $V(1)=1$ on single-period paths. Define $g(a)=V(a)$ for $a\ge0$. Repeated concatenation and sign invariance give
\[
V(x)=\sum_h g(|x_h|),\qquad g(0)=0,\quad g(1)=1.
\]
For cumulative accounting, preservation of same-sign splits gives $g(a+b)=g(a)+g(b)$ for $a,b\ge0$. Nonnegativity makes $g$ nondecreasing: if $b\ge a$, then $g(b)=g(a)+g(b-a)\ge g(a)$. Additivity and normalization imply $g(q)=q$ for every nonnegative rational $q$. Rational sequences approaching any real $a\ge0$ from below and above then give $g(a)=a$. Thus $V(x)=\sum_h|x_h|$ without an additional continuity assumption.

For quadratic accounting, apply $V(tx)=t^2V(x)$ to the single-period path $x=1$. This gives $g(t)=t^2$ for every $t\ge0$, hence $V(x)=\sum_hx_h^2$. Conversely, the absolute and squared sums each satisfy the common properties and their respective split or scaling condition.

At a recipient, evaluate each standardized source path separately. Normalizing its score by the sum of source scores gives $C^{(H,1)}_{ij}(S)$ in the cumulative case and $C^{(H,2)}_{ij}(S)$ in the quadratic case. The own-impact response ensures a positive denominator for positive-diagonal scatter and $\Phi_0=I$.
\end{proof}

\subsection{Exact comparison with squared allocation}\label{sec-ratio}
For common scatter $S$, let $a_{ij}=\sum_h|x_{ijh}(S)|$ and
$b_{ij}=\sum_h|x_{ijh}(S)|^2$. For active channels $a_{ij}>0$,
\begin{equation}
\frac{C^{(H,2)}_{ij}(S)}{C^{(H,1)}_{ij}(S)}
=\frac{b_{ij}/a_{ij}}{\sum_k b_{ik}/\sum_k a_{ik}}.
\label{eq:ratio}
\end{equation}
In Equation~\eqref{eq:ratio}, the numerator is the channel's energy per unit cumulative deviation; the denominator averages this quantity across the recipient's sources, weighted by cumulative amplitude. Squaring raises the share exactly when that ratio exceeds the average. If the two rows agree, each active-channel ratio must equal the row average; conversely, a common ratio makes the two normalizations identical. Channels with zero amplitude have zero energy and zero share under both accounts.

\paragraph*{Proof of the attribution and channel-priority proposition in the main paper.}
For each active channel, set $d_{ij}=a_{ij}^2/b_{ij}$. The statement compares two channels with $a_{ij}>a_{ik}>0$ and their normalized absolute and squared shares.

\begin{proof}
For $u_h=|x_{ijh}|\ge0$, positivity of the cross-products and Cauchy--Schwarz give
\[
\sum_hu_h^2\le\Big(\sum_hu_h\Big)^2\le H\sum_hu_h^2.
\]
Thus $1\le d_{ij}\le H$, and $b_{ij}=a_{ij}^2/d_{ij}$. Within a recipient row the normalization is common across sources. Therefore, conditional on $a_{ij}>a_{ik}>0$, squared allocation reverses the absolute ordering if and only if
$a_{ij}^2/d_{ij}<a_{ik}^2/d_{ik}$, which is the stated channel-priority condition.

\textit{Gaussian construction.} Take $H\ge3$, $0<q<1$, $a,b>0$, $S=I_3$, and
\[
A=\begin{pmatrix}0&b&a\\0&0&0\\0&0&q\end{pmatrix},\qquad
G_1=\sum_{\ell=0}^{H-2}q^\ell,\quad G_2=\sum_{\ell=0}^{H-2}q^{2\ell}.
\]
The eigenvalues are $0,0,q$, so it is stable. At $S=I_3$, both cross-responses are zero at impact. Channel $2\to1$ equals $b$ at horizon one and zero thereafter, whereas channel $3\to1$ equals $aq^{h-1}$ for $h=1,\ldots,H-1$. Hence their absolute totals are $b$ and $aG_1$, and their squared totals are $b^2$ and $a^2G_2$. Thus RAC ranks $3\to1$ above $2\to1$ and squared allocation does the reverse precisely when $a\sqrt{G_2}<b<aG_1$. Since $H-1\ge2$ and $q>0$,
\[
G_1^2-G_2=2\sum_{0\le r<s\le H-2}q^{r+s}>0,
\]
so the interval is nonempty.

Finite-horizon responses and row shares are continuous in $(A,S)$ wherever the scatter diagonal is positive; own impact responses keep their row denominators positive. Stability and positive definiteness are open conditions. The strict ordering margins therefore persist on an open neighborhood. For Gaussian innovations with any scatter in this neighborhood, population rank scatter is a common positive multiple of covariance, so the same reversal separates RAC from the calibrated squared specialization.

\end{proof}

\subsection{Impact inequality}\label{app-impact-inequality}
For positive-semidefinite $S$ with positive diagonal, at impact $\Phi_0=I$, so
$x_{ij0}=s_{ij}/\sqrt{s_{jj}}$.
Let
$\delta_{ij}=|s_{ij}|/\sqrt{s_{ii}s_{jj}}=|\rho_{ij}(S)|$.
Then $|x_{ij0}|=\sqrt{s_{ii}}\delta_{ij}$, and
\[
C^{(1,1)}_{ij}=\frac{\delta_{ij}}{\sum_k\delta_{ik}},\qquad
C^{(1,2)}_{ij}=\frac{\delta_{ij}^2}{\sum_k\delta_{ik}^2}.
\]
Writing $\bar\delta_i=\sum_k\delta_{ik}^2/\sum_k\delta_{ik}$ gives $C^{(1,2)}_{ij}/C^{(1,1)}_{ij}=\delta_{ij}/\bar\delta_i$ for $\delta_{ij}>0$.
For $\delta_{ij}=0$, both shares are zero and their ratio is undefined.
The diagonal is $\delta_{ii}=1$, and all entries lie in $[0,1]$.
Thus
$\sum_k\delta_{ik}^2\le\sum_k\delta_{ik}$,
which implies $C^{(1,2)}_{ii}\ge C^{(1,1)}_{ii}$.
Subtracting own shares from one shows that the total external share under absolute aggregation is at least as large as under squared aggregation.

For $H>1$, $b_j$ sums the squared responses at each horizon. Paths $(1,1)$ and $(1,0)$, for example, both have $b_j/a_j=1$ and receive the same shares under the two rules. The common-scatter ratio compares GFEVD and RAC under elliptical calibration, or isolates aggregation when scatter is held fixed.

\subsection{Source breadth and relative duration}\label{app-breadth-duration}
The channel comparison extends to own versus external attribution by accounting for how cumulative amplitude is divided among external sources. Fix a recipient and common scatter. Write $a_j=\sum_h|x_{ijh}|$, $b_j=\sum_h x_{ijh}^2$, and separate own source $i$ from external sources $j\ne i$. Let
\[
A_E=\sum_{j\ne i}a_j,\quad A_O=a_i,\quad
B_E=\sum_{j\ne i}b_j,\quad B_O=b_i,
\]
and, for positive external response, define
\[
n_E=\frac{A_E^2}{\sum_{j\ne i}a_j^2},\qquad
d_E=\frac{\sum_{j\ne i}a_j^2}{B_E},\qquad d_O=\frac{A_O^2}{B_O}.
\]
Cauchy--Schwarz gives $1\le n_E\le K-1$ and $1\le d_E,d_O\le H$. In particular, $d_E$ is the external squared-response-weighted mean of the channel durations.

The product $n_Ed_E=A_E^2/B_E$ is the effective breadth across external source--horizon pairs: normalizing their absolute responses by $A_E$ makes their squared shares sum to $B_E/A_E^2$. Since $B_E=A_E^2/(n_Ed_E)$ and $B_O=A_O^2/d_O$, putting $r=A_E/A_O$ and $\delta=d_E/d_O$ yields
\[
F_1=\frac{r}{1+r},\qquad F_2=\frac{r^2}{r^2+n_E\delta}.
\]
Consequently, $F_1>1/2>F_2$ exactly when $1<r<\sqrt{n_E\delta}$.

For two recipients $i,k$ with positive external amplitudes, define $r_i=A_{E,i}/A_{O,i}$ and $g_i=n_{E,i}d_{E,i}/d_{O,i}$, with analogous definitions for $k$. Their squared-share difference has the sign of $r_i^2g_k-r_k^2g_i$, whereas their amplitude-share difference has the sign of $r_i-r_k$. Thus, conditional on $r_i>r_k$, opposite recipient priorities occur exactly when $1<r_i/r_k<\sqrt{g_i/g_k}$. 

To separate the two factors, hold $r$ fixed and evaluate $F_2$ at $(n_E,\delta)$, $(1,\delta)$, $(n_E,1)$ and $(1,1)$. These algebraic comparisons vary breadth and relative duration at the same measured amplitude ratio. Among actual reversals, breadth alone is sufficient when $r^2<n_E$, and relative duration alone is sufficient when $r^2<\delta$. This partitions observations into breadth only, duration only, either factor individually, and both jointly required.

\subsection{Breadth and duration at a fixed amplitude ratio}\label{app-suite-geometry}
A Gaussian star system isolates breadth and duration at a fixed external/own cumulative amplitude ratio. Consider a system with covariance $I_{J+1}$, one recipient of persistence $a$ and $J$ external sources, each with persistence $q$ and coefficient $b$ into the recipient. The own response is $a^h$; each external response at $h\ge1$ is
\[
b\sum_{\ell=0}^{h-1}a^\ell q^{h-1-\ell}.
\]
We set $b$ so that the total external/own amplitude ratio is $r$. Equal external channels give $n_E=J$.

For horizon $H$, define
\[
G_O=\sum_{h=0}^{H-1}a^h,\qquad
G_E=\sum_{h=1}^{H-1}\sum_{\ell=0}^{h-1}a^\ell q^{h-1-\ell}.
\]
The coefficient is $b=rG_O/(JG_E)$ for the positive $G_E$ in the grid. This fixes the cumulative external/own ratio at $r$ before computing either allocation.

The grid crosses $H=5,10,20$, $J=1,2,4,8,16$, $q=0,.2,.4,.6,.8,.95$, $a=0,.4,.8$, and $r=.25,.5,.75,1,1.25,1.5,2,3$, giving 2,160 stable systems. Alternating channel signs leaves the absolute and squared allocations unchanged; at $H=1$ all dynamic external responses vanish. An equal-response rank-one scatter supplies a control in which both allocations coincide.

\FloatBarrier

\section{Recovering priorities from estimated responses}\label{app-estimation}
Sections~\ref{app-local-construction}--\ref{app-complete-construction} specify the estimator. The consistency argument starts with residual scatter in a stationary system (Section~\ref{app-stationary-foundation}), extends to local dynamics (Section~\ref{app-local-consistency}), and incorporates the implementation adjustments in Proposition~\ref{prop-regularized-target}.

\subsection{Recovery of network priorities}\label{app-monitoring-resolution}
The proof first bounds recipient-share error by the error in fitted response magnitudes, including both dynamics and scatter. Comparing that bound with the list's cutoff gap and the own/external boundary gives the recovery conditions. Summing score errors over replaced names bounds the selected list's population-strength shortfall.

\paragraph*{Proof of the priority-recovery proposition in the main paper.}
Use $a_{ij}=\sum_h|x_{ijh}|$, $A_i=\sum_j a_{ij}$ and $\delta_i=\sum_{j,h}|\widehat x_{ijh}-x_{ijh}|$, with fitted counterparts defined in the same way. Scores are external shares in percent, $f_i=100\sum_{j\ne i}a_{ij}/A_i$. Write $\epsilon_i=\min\{100,100\delta_i/\widehat A_i\}$, $\epsilon=\max_i\epsilon_i$, and $d=|I_q\setminus\widehat I_q|$.
\begin{proof}
For nonnegative vectors $a,b$ with positive sums $A,B$, write $d_1=\|a-b\|_1$. If $A\ge B$, then
\[
\left\|\frac aA-\frac bB\right\|_1
\le\frac{d_1+|A-B|}{A}\le\frac{2d_1}{A}.
\]
Interchanging the vectors covers $B>A$. A subset's probability difference is at most half the $\ell_1$ distance between probability vectors. Moreover $\|\widehat a_i-a_i\|_1\le\delta_i$ by the reverse triangle inequality. Apply the subset bound to the external sources and multiply by 100 to obtain the response-to-score bound in the main paper. The classification assertion follows immediately.

For selection, common list members cancel. Since $\widehat I_q$ maximizes the fitted sum,
\[
\sum_{i\in I_q} f_i-\sum_{i\in\widehat I_q} f_i
\le\sum_{i\in I_q\setminus\widehat I_q}(f_i-\widehat f_i)
  +\sum_{i\in\widehat I_q\setminus I_q}(\widehat f_i-f_i)
\le2d\epsilon.
\]
Every target-selected score exceeds every omitted score when its cutoff gap exceeds $2\epsilon$; the same argument with target and fit interchanged proves the fitted-gap criterion.
\end{proof}

\paragraph*{Consistent recovery.}
At fixed dimension and horizon, the dynamic and rank-scatter consistency results below give convergence of every response and hence $\delta_i\xrightarrow{p}0$. Positive own-impact totals imply $\epsilon=\max_i\epsilon_i\xrightarrow{p}0$. At a fixed stationary model or a fixed interior rescaled date satisfying the estimation conditions, a positive population cutoff gap eventually exceeds $2\epsilon$ with probability tending to one, yielding exact Top-$q$ recovery. Nonzero population distances from 50 give consistent own/external classification. The objective-loss bound tends to zero even at a tied cutoff.

Complete-estimation simulations use known target paths to evaluate $\delta_i$; Table~\ref{tbl-suite-classification} reports classification accuracy.

\subsection{Local dynamics and initial-history tuning}\label{app-local-construction}\label{app-final-estimation}
The current window estimates coefficient levels and local time slopes. Use the local-linear design $d_s=(1,x_s',z_s,z_sx_s')'$, time coordinate $z_s=(s-t)/m$, and self-weight $w_s=(1+\|x_s\|)^{-1}$. A local-constant self-weighted Huber pilot supplies the fixed positive residual MAD $s_j$. For each equation, minimize
\[
\frac{\sum_s w_s s_j^2\rho_{1.345}((y_{js}-d_s'b)/s_j)}{\sum_s w_s}
+\frac{\lambda_n}{2}b'Pb,\qquad
P=\operatorname{diag}(0_{K+1},1_{K+1}),\quad
\lambda_n=\gamma\sqrt{(K+1)/n}.
\]
Here $n$ counts complete adjacent pairs. The penalty acts on the time intercept and time slopes, leaving coefficient levels unpenalized. Median/MAD coordinate calibration preserves positive affine equivariance, and the fit imposes no spectral constraint.

\subsubsection{Initial-history tuning}
The candidate grid is $\gamma\in\{0,.01,.03,.1,.3,1,3,10\}$. Within the initial 500 positions, two 250-position training windows end at zero-indexed positions 349 and 449; the next 50 positions form each validation block. Each fold obtains its coordinate calibration strictly before its training window.

One-step validation predictions use observed lags and the extrapolated fitted line. Pooling all equations, the criterion compares weighted Huber losses with those of the common local-constant pilot. Subtracting the common pilot loss preserves candidate ordering and gives a bounded contrast on bounded coefficient sets.

Choose the strength minimizing pooled validation loss; ties within $10^{-12}$ favor the larger candidate. The selected empirical values are 10 for stocks, 3 for currencies and 10 for futures. Subsequent fits use the same selected strength and the initial-500 coordinate calibration.

\subsection{From local dynamics to the RAC estimator}\label{app-complete-construction}\label{app-regularized-implementation}
Historical observations select coefficient directions; the current window estimates coefficients and residual scatter.

\paragraph*{Historical directions.}
A template specifies which lag directions enter a recipient's regression. Each retains the recipient's own lag and may add the equal-weight aggregate, rank-factor directions, or $q$ individual lags. Distinct sparse sizes are $q\in\{1,2,4,8,16,K-1\}$, capped at $K-1$; the own-lag template adds no direction. For endpoint $t$, the historical block ends before the current 250 response positions and contains at most 1,500 observations. Only originally adjacent complete pairs enter. Coordinate calibration and the slope penalty remain fixed at their initial-500 values.

Historical self-weighted Huber fits score the templates. Let $n_h$ count the usable historical pairs and $s_{i,h}$ denote the historical pilot's scaled residual MAD. Factor directions are eigenvectors of the historical lag correlation whose eigenvalues exceed $1+2\pi\sqrt{K/n_h}/3+\pi^2K/(9n_h)$, rescaled by lag MADs. Sparse directions follow the full historical pilot's absolute coefficients.

\paragraph*{Template selection.}
With historical effective sample size $n_e=(\sum w_s)^2/\sum w_s^2$, template criterion $\ell_{ij}$ equals the weighted Huber loss contrast against the historical pilot plus $\operatorname{tr}(H_{ij}^{-1}J_{ij})$, where $H_{ij}$ is the weighted empirical Huber Hessian and $J_{ij}$ the sum of outer products of the weighted, sample-normalized scores. Define
\[
L_{ij}\propto\exp\{-n_e\ell_{ij}/s_{i,h}^{2}-b_j\},\qquad
\widehat p=\arg\max_{p\in\Delta}\left\{\sum_i\log\Big(\sum_jL_{ij}p_j\Big)+\sum_j a_j\log p_j\right\}.
\]
The search cost is $b_j=\log {K-1\choose q}$ for a sparse template and zero otherwise. Pseudo-counts $a_j$ allocate $.4,.1,.1$ to the own, aggregate and factor templates and the remaining $.4$ equally across sparse sizes. For each recipient, we select the template maximizing $L_{ij}\widehat p_j$ and refit its directions on the current pairs.

\paragraph*{Current coefficients and residuals.}
For the VAR(1) implementation, let $n$ denote the number of usable current pairs and let $\widehat A$ be the local coefficient-level estimate, expressed in the calibrated coordinates. Let $B$ be the coefficient matrix obtained by refitting the historically selected directions on the current pairs. For recipient $i$, set
\[
\widetilde A_{ij}=\widehat A_{ij}
+\operatorname{clip}\{B_{ij}-\widehat A_{ij},-c_n\widehat s_i,c_n\widehat s_i\},
\qquad c_n=(\log n/n)^{1/4},
\]
where $\widehat s_i$ is the fixed pilot scale. Transform the coefficients back to observed units and form aligned residuals $y_s-\widetilde A y_{s-1}$. A common intercept translation leaves MADs and Kendall coefficients unchanged.

\paragraph*{Residual scatter.}
Discard current residual rows with the largest calibrated lag norms above the median plus two scaled MADs, up to a fraction
$q_n=\min\{.3,2\sqrt{\log n/n}\}$.
Use the retained residual rows to form the projected sine--Kendall correlation $U$ and marginal MAD matrix $D$. The covariance benchmark uses all aligned current residuals.

\paragraph*{Shrinkage target.}
A target $T=I+\rho(gg'-I)$ pools signed concordance. The sign vector $g$ is selected from the all-positive vector and signs of the leading and smallest eigenvectors fitted to the first half of the historical residuals. The second half selects the largest absolute pooled signal divided by its bounded-kernel dispersion scale. We activate the target when this ratio exceeds the tuning threshold 2.5; otherwise $T=I$. Historical halves use the same leverage rule, and $\rho$ is the sine transform of mean signed Kendall concordance in the retained current residuals.

\paragraph*{Shrinkage intensity.}
Write $N$ for the retained count, $g_{s,ij}$ for the mean sign-product kernel over all partners of row $s$, and $\widehat\tau_{ij}=N^{-1}\sum_sg_{s,ij}$. The shrinkage scale is
\[
V_N=\sum_{i\ne j}\left\{\frac\pi2\cos(\pi\widehat\tau_{ij}/2)\right\}^{\!2}
\frac{4(N-1)}{N(N-2)^2}\sum_s(g_{s,ij}-\widehat\tau_{ij})^2.
\]
With $d=\|U-T\|_F$, set $\alpha=(1-V_N/d^2)_+$, taking $\alpha=0$ when $d=0$, and define
\[
\widetilde R=\Pi_{\mathcal C}\{\alpha U+(1-\alpha)T\},\qquad
\widetilde\Omega=D\widetilde R D,
\]
where $\mathcal C$ is the set of positive-semidefinite correlation matrices.

The coefficient adjustment is bounded by $c_n\widehat s_i$, the deleted fraction by $q_n$, and the correlation perturbation by $\sqrt{V_N}$. These bounds vanish under the consistency conditions below, preserving the population target.

\subsection{Stationary foundation for residual-scatter consistency}\label{app-stationary-foundation}\label{ec-sec-estimation}\label{sec-consistency}
The stationary proof proceeds from coefficient consistency to residual quantiles and Kendall averages. It uses conditional symmetry; Section~\ref{app-local-consistency} extends the argument to local dynamics under ergodicity and a common conditional Huber location, covering asymmetric iid innovations and conditionally symmetric scale dependence. For each equation, estimate $\beta_i$ by minimizing
\begin{equation}
\sum_{t=1}^{T} w(X_t)\,\rho_c\!\left(\frac{y_{it}-b'X_t}{\hat s_i}\right),
\label{eq:huber}
\end{equation}
where $\rho_c(u)=u^2/2$ for $|u|\le c$ and $\rho_c(u)=c|u|-c^2/2$ for $|u|>c$ is the Huber loss with tuning constant $c=1.345$, and $\hat s_i$ is a positive estimated MAD pilot scale held fixed during coefficient optimization. Under conditional symmetry the Huber score is centered at the conditional median for every positive scale. For the moment-free consistency result use $w(x)=(1+\|x\|)^{-1}$.

Assume a stationary ergodic sequence, a unique conditional median and conditionally symmetric innovations. Identification additionally requires the conditional Huber risk to be uniquely minimized at zero uniformly over $s\in[\underline s,\overline s]$: for almost every past, every $s$ in that interval, and every $a\ne0$,
\[
E\!\left[\rho_c\!\left(\frac{\varepsilon_{it}-a}{s}\right)-\rho_c\!\left(\frac{\varepsilon_{it}}{s}\right)\middle|\mathcal F_{t-1}\right]>0.
\]
The pilot scale must lie in a fixed compact interval $[\underline s,\overline s]\subset(0,\infty)$ with probability tending to one, and the coefficient returned by the solver must minimize the Huber criterion at that scale (up to $o_p(1)$ in normalized objective value). Together with
$P(v'X_t\ne0)>0$ for every nonzero $v$, these conditions identify the coefficient.

\begin{theorem}[Feasible consistency of rank-scatter connectedness]\label{thm:consistency}
Under these estimation conditions and continuous innovation marginals having uniquely crossed medians and MADs, the self-weighted Huber M-estimator is consistent. Sample MADs and Kendall ranks of its residuals converge to their innovation functionals. Hence, for fixed $K,p$ and finite $H_{\max}$,
\[
\max_{i,j,\;1\le H\le H_{\max}}
|\widehat{\mathrm{RAC}}_{ij}(H)-\mathrm{RAC}_{ij}(H)|\xrightarrow{p}0.
\]
A population transmission hub with a strictly positive gap from every competitor is recovered with probability tending to one; the same conclusion holds for edge rankings with positive gaps.
\end{theorem}
\label{app-consistency}
\subsubsection{Coefficient identification and uniform convergence}
The bounded difference between two Huber losses supplies the population criterion, including when the individual losses have infinite expectation. For an equation write $\delta=b-\beta_i$ and
\[
\ell_t(\delta;s)=w(X_t)\left\{\rho_c\!\left(\frac{\varepsilon_{it}-\delta'X_t}{s}\right)-\rho_c\!\left(\frac{\varepsilon_{it}}{s}\right)\right\},
\qquad w(x)=(1+\|x\|)^{-1},
\]
where $\rho_c$ is the Huber loss with tuning constant $c=1.345$ and $s\in[\underline s,\overline s]$. The Huber derivative $\psi_c(u)=\min(|u|,c)\operatorname{sgn}(u)$ is bounded by $c$, so
$|\ell_t(\delta;s)|\le c\|\delta\|w(X_t)\|X_t\|/\underline s\le c\|\delta\|/\underline s$.
The contrast is Lipschitz in $\delta$ uniformly in $s$.

It is also Lipschitz in $s$ on this interval: writing $q(u)=u\psi_c(u)$, which is $2c$-Lipschitz, gives
$|\partial_s\ell_t(\delta;s)|\le 2c\|\delta\|w(X_t)\|X_t\|/\underline s^2\le2c\|\delta\|/\underline s^2$
where the derivative exists, and the bound extends across the Huber knots by absolute continuity.

The ergodic theorem, a finite net in $(\delta,s)$, and these uniform Lipschitz bounds give uniform convergence on compact coefficient sets and the scale interval to
$D_i(\delta;s)=E[\ell_t(\delta;s)]$. This expectation is defined even when
$E|\varepsilon_i|=\infty$, because the contrast itself is bounded.

Conditional symmetry and oddness of $\psi_c$ give
$E[\psi_c(\varepsilon_{it}/s)\mid\mathcal F_{t-1}]=0$ for every fixed $s>0$. The conditional Huber separation condition in Section~\ref{ec-sec-estimation} requires that, for almost every past, every $s\in[\underline s,\overline s]$, and every $a\ne0$, the conditional Huber risk difference
\[
Q_i(a\mid\mathcal F_{t-1};s)=
E[\rho_c((\varepsilon_{it}-a)/s)-\rho_c(\varepsilon_{it}/s)\mid\mathcal F_{t-1}]
\]
is strictly positive. For any fixed $\delta\ne0$, the regressor identification condition gives $P(\delta'X_t\ne0)>0$. Since $\delta'X_t$ is $\mathcal F_{t-1}$-measurable, $Q_i(\delta'X_t\mid\mathcal F_{t-1})$ is nonnegative and is strictly positive whenever $\delta'X_t\ne0$. Therefore, for every $s$ in the allowed interval,
\[
D_i(\delta;s)=E\{w(X_t)Q_i(\delta'X_t\mid\mathcal F_{t-1};s)\}>0,
\]
using $w(X_t)>0$.

For any fixed radius $r>0$, continuity in $(\delta,s)$ and compactness give a strictly positive minimum of $D_i$ on $\|\delta\|=r$, $s\in[\underline s,\overline s]$. Uniform convergence on this set makes the sample loss difference positive there with probability tending to one, uniformly over the fitted scale. Convexity along the segment from zero to any candidate minimizer then confines every approximate minimizer to the radius-$r$ ball. Letting $r$ be arbitrary proves consistency.

Sample minimizers exist whenever the design has full column rank: the positive-weight Huber objective is then coercive. The matrix $E[w(X)^2XX']$ is positive definite by the identification condition, and its ergodic sample counterpart implies that full rank holds with probability tending to one. The argument applies to any weight satisfying $w(x)\|x\|\le C$ for finite $C$, including the specified $w(x)=(1+\|x\|)^{-1}$.

\subsubsection{From coefficient consistency to residual ranks}
Let $\Delta_T=\widehat A-A$ collect coefficient errors, and
$\widehat\varepsilon_t=\varepsilon_t-\Delta_T X_t$.
For any fixed $\eta,M>0$, the following is a deterministic sample inequality:
\begin{equation}
\frac1T\sum_t{\bf1}\{\|\Delta_T X_t\|>\eta\}
\le {\bf1}\{\|\Delta_T\|>\eta/M\}
+\frac1T\sum_t{\bf1}\{\|X_t\|>M\}.
\label{eq-fraction}
\end{equation}
The relevant quantity is the fraction of residuals perturbed beyond a fixed tolerance. Equation~\eqref{eq-fraction} controls that fraction even with heavy-tailed regressors. For fixed $M$, the first term vanishes in probability and the second converges by ergodicity to $P(\|X_t\|>M)$. Sending $M$ to infinity makes the fraction of residuals exceeding that perturbation tolerance $o_p(1)$.

For each marginal CDF,
\[
F^\varepsilon_{T,i}(z-\eta)-b_T(\eta)
\le F^{\widehat\varepsilon}_{T,i}(z)
\le F^\varepsilon_{T,i}(z+\eta)+b_T(\eta),
\]
where $b_T(\eta)$ is the fraction exceeding the perturbation tolerance in \eqref{eq-fraction}.
Ergodicity and continuity give uniform marginal empirical-CDF convergence. Quantile crossing and uniqueness give convergence of the residual median. Applying the same bound to absolute deviations from that consistent median gives MAD convergence.

For Kendall ranks, observations exceeding the perturbation tolerance enter at most a fraction
$2b_T(\eta)+O(T^{-1})$ of unordered pairs. For pairs with both residual perturbations within tolerance, a sign can change only if
$|\varepsilon_{it}-\varepsilon_{is}|\le2\eta$ or
$|\varepsilon_{jt}-\varepsilon_{js}|\le2\eta$.

Let $P_T=T^{-1}\sum_t\delta_{\varepsilon_t}$ be the empirical measure of the innovation vectors. Stationary ergodicity gives $P_T\Rightarrow P$ almost surely, hence $P_T\otimes P_T\Rightarrow P\otimes P$. The all-pairs Kendall statistic is the integral of its bounded kernel against $P_T\otimes P_T$, up to the negligible diagonal/$T/(T-1)$ correction; its discontinuity set has $P\otimes P$ measure zero because the relevant marginals are continuous. Thus it converges to the Kendall functional of two independent draws from $P$, without assuming independence of observations over time.

For fixed $\eta$, the proportion of near-tie pairs is controlled by $P_T\otimes P_T$ on the corresponding closed near-tie sets; by weak convergence its limit superior is bounded by their $P\otimes P$ probabilities, which tend to zero as $\eta\downarrow0$. The perturbation-pair fraction also vanishes because $b_T(\eta)=o_p(1)$. Boundedness of the kernel then yields
$\widehat\tau_{ij}(\widehat\varepsilon)\to_p\tau_{ij}(P)$.

Sine continuity and projection nonexpansiveness now imply
$\widehat\Omega^\dagger\to_p\Omega^\dagger$.

\subsubsection{Network and ranking convergence}
For finite $H_{\max}$, each $\Phi_h$ is a polynomial in the coefficients, making every allocation numerator continuous in coefficients and positive-diagonal scatter. The own impact term provides a strictly positive limiting row denominator. Because the sets of indices and horizons are finite, entrywise continuous mapping yields the finite-grid uniform result.

If $j^*$ has a population transmission gap $\gamma>0$, uniform entrywise consistency makes every estimated transmission within $\gamma/3$ of its target with probability tending to one. The hub is then $j^*$.

\subsection{Local consistency of the dynamic component}\label{app-local-consistency}
\subsubsection{Local approximation and the pilot}
\label{app-final-local-theory}
Fix $K,p,H$ and an interior rescaled date $u\in(0,1)$. Write
\[
Y_{t,T}=a(t/T)+\sum_{\ell=1}^p A_\ell(t/T)Y_{t-\ell,T}+L(t/T)\xi_t.
\]
The innovations $\xi_t$ are stationary ergodic and adapted to a stationary filtration $\mathcal G_t$ containing their past. Their conditional law given $\mathcal G_{t-1}$ has an everywhere positive density, and $E\|\xi_t\|^\delta<\infty$ for some $0<\delta\le1$. Conditional densities make finite frozen lag vectors absolutely continuous with full support. Identification uses the common conditional Huber-location condition below, allowing dependence and nonelliptical innovations.

The coefficient functions are continuously differentiable with bounded derivatives, $L$ is nonsingular, and the companion matrices contract uniformly in a common matrix norm: $\sup_v\|F(v)\|\le r<1$. Use the causal solution, or an initialization whose contribution on the local window vanishes in probability. Extend coefficients constantly outside their time domain when constructing the causal representation. For each frozen date, assume uniquely crossed marginal medians and MADs. These sufficient conditions include innovations with regularly varying indices below one by choosing $\delta$ below that index.

Let $t=\lfloor Tu\rfloor$, $h=m/T\to0$, and $m\to\infty$. The local design uses the original positions $t-m+1,\ldots,t$ and all $p$ original lags. Omitted lag--response tuples may number $o_p(m)$; their time coordinates are retained. Thus $n/m\to_p1$.

An affine coordinate calibration $\mathcal C_T=(\mu_{0,T},D_{0,T})$ has positive diagonal scales and satisfies $\|\mu_{0,T}\|+\|D_{0,T}\|+\|D_{0,T}^{-1}\|=O_p(1)$. This includes median/MAD calibration from an initial fixed history of at least three observations: under the causal construction, this finite vector converges in distribution to its initial frozen counterpart, whose coordinate MADs are positive almost surely. The proof below is uniform on compact sets of calibrations and does not require their independence from subsequent observations.

In calibrated coordinates let $X_s=(1,x_s')'$, $d_s=(X_s',z_sX_s')'$, $z_s=(s-t)/m$, and $w_s=(1+\|x_s\|)^{-1}$. The number of level coordinates is $q=Kp+1$. A local-constant pilot returns $(\widetilde\beta,\widetilde s)$ satisfying its self-weighted Huber score up to $o_p(1)$ and
\[
\widetilde s=1.4826\,\operatorname{Med}_s
 \left|Y_{js}-X_s'\widetilde\beta
 -\operatorname{Med}_v(Y_{jv}-X_v'\widetilde\beta)\right|+o_p(1).
\]
The score here is normalized by $\sum_s w_s$ and uses $\psi_c((Y_{js}-X_s'\widetilde\beta)/\widetilde s)$, with pilot constant $c=1.345$. The coefficient optimization holds this positive scale fixed. The asymptotic results below specify vanishing numerical-error sequences.

\paragraph*{Conditional score centering and the intercept.}
For one equation and a fixed calibration, let $\varepsilon_{j,s}(u)$ be its frozen innovation, $s_{j0}(u)$ its marginal scaled MAD, and $\beta_j(u)$ its structural intercept and lag coefficients. For each $\sigma>0$, assume a deterministic location $\kappa_j(u,\sigma)$ common to past histories such that
\[
E\!\left[\psi_c\{(\varepsilon_{j,s}(u)-\kappa_j(u,\sigma))/\sigma\}
\mid\mathcal G_{s-1}\right]=0\quad\text{almost surely}.
\]
This identifies a common conditional Huber location while allowing the conditional distribution to vary over time. Iid innovations satisfy the condition relative to their natural past even when asymmetric. Conditionally centrally symmetric innovations satisfy it with $\kappa_j=0$. An example is a stationary ergodic innovation sequence formed by a predictable nonsingular matrix scale applied to iid centrally symmetric shocks independent of the past, provided the stated fractional moment exists. Its scale and magnitudes may remain serially dependent.

Taking expectations identifies $\kappa_j(u,\sigma)$ as the marginal Huber location. The score is bounded, so its expectation exists without an innovation mean. It is continuous and strictly decreasing in $\kappa$, with limits $c$ and $-c$; positive density makes the root unique. Dominated convergence and strict monotonicity give continuity in $\sigma$. With $e_0$ selecting the intercept, conditional centering identifies
$\beta_j^H(u,\sigma)=\beta_j(u)+\kappa_j(u,\sigma)e_0$.
Indeed, $X_s$ is measurable with respect to $\mathcal G_{s-1}$, and boundedness of $w(X_s)X_s$ and $\psi_c$ permits iterated expectations without lag means:
\[
\begin{aligned}
&E\!\left[w(X_s)X_s\psi_c\{(\varepsilon_{j,s}(u)-\kappa_j(u,\sigma))/\sigma\}\right]\\
&\quad=E\!\left[w(X_s)X_s E\!\left[\psi_c\{(\varepsilon_{j,s}(u)-\kappa_j(u,\sigma))/\sigma\}\mid\mathcal G_{s-1}\right]\right]=0.
\end{aligned}
\]
Positive conditional density gives strict conditional Huber separation at every nonzero shift. Write $\beta_j^H(u)=\beta_j^H(u,s_{j0}(u))$. Back in observed units, lag coefficients remain $A_\ell(u)$; translating innovations leaves their MADs and Kendall coefficients unchanged.

\begin{lemma}[Local empirical approximation and the pilot]
\label{lem-final-pilot}
Under the preceding conditions, the local empirical law of the observations, original lags and deterministic time coordinate converges in probability to the frozen joint law times the uniform law on $[-1,0]$. The centered, self-weighted Huber contrasts obey a uniform law on compact coefficient, positive-scale and calibration sets. The local-constant pilot converges to $\beta_j^H(u)$, and its residual MAD converges to $s_{j0}(u)$. In particular, the pilot consistently estimates the structural lag coefficients.
\end{lemma}
\begin{proof}
\textit{Local approximation.} Uniform contraction gives a causal linear representation with exponentially decaying coefficients. For $0<\delta\le1$, subadditivity of $\|\cdot\|^\delta$ bounds the contribution of innovations more than $J$ lags old by $Cr^{\delta J}$. Comparing the varying coefficients with coefficients frozen at $u$ gives, uniformly on the local window,
\[
 E\|Y_{s,T}-Y_s(u)\|^\delta
 \le C\{h+T^{-1}\}^\delta.
\]
To obtain this bound, telescope each product of companion matrices; its $j$th-lag difference is bounded by $Cj r^{j-1}(h+j/T)$, whose $\delta$th powers are summable. The intercept and innovation-loading terms satisfy the same summable bound. The argument also controls the finite collection of original lags.

\textit{Frozen empirical law.} The frozen process is a measurable shift-equivariant function of the stationary ergodic innovation sequence, so its observation--lag vector $W_s(u)$ is stationary ergodic. For a bounded continuous $f$, the average over any block of length tending to infinity converges in probability to $Ef(W_0(u))$. This also holds for blocks whose starting positions vary: stationarity gives each block average the same distribution as an initial block of that length. Partitioning $[-1,0]$ into finitely many subintervals and applying this fact to each corresponding block gives, for continuous $g$,
\[
\frac1m\sum_{s=t-m+1}^{t}f(W_s(u))g(z_s)
\xrightarrow{p}Ef(W_0(u))\int_{-1}^{0}g(z)\,dz.
\]
Step-function approximation of $g$ and finite sums of product test functions establish the joint empirical-law limit.

To replace frozen observations by actual observations, restrict them to a compact set and apply the coupling bound and uniform continuity. The fractional-moment bound controls the discarded fraction. Removing $o_p(m)$ tuples changes bounded averages by $o_p(1)$ and leaves the limiting time design unchanged. Ergodicity supplies the required law of large numbers.

\textit{Uniform loss convergence.} For a scale $\sigma$ in a compact positive interval, the loss difference
\[
w_s\sigma^2\{\rho_c((Y_{js}-d_s'b)/\sigma)
                  -\rho_c((Y_{js}-d_s'b_0)/\sigma)\}
\]
has absolute value at most $c\sigma w_s\|d_s\|\|b-b_0\|$. Since $|z_s|\le1$ and $w_s\|X_s\|$ is bounded, the envelope is bounded. The contrast is uniformly Lipschitz in coefficients and scale on these sets; truncation of the raw observations supplies uniform continuity in affine calibrations. Finite nets and the empirical-law assertion therefore give the stated uniform law, including a positive limiting denominator $n^{-1}\sum_s w_s$.

\textit{Pilot tightness.} We establish tightness of the data-dependent pilot coefficients and scale. With probability tending to one, more than three quarters of the local observations have both response and regressors in a common compact set. Consequently the residual MAD at any coefficient $b$ is bounded above by $C(1+\|b\|)$. Suppose a sequence of pilot coefficients diverged, and write $b=Rv$, $R\to\infty$, $\|v\|=1$. A subsequence has $v\to v_0$ and $\widetilde s/R\to a\in[0,C]$. The score in direction $v$ then converges to
\[
-\frac{E\{w(X)(v_0'X)\psi_c(v_0'X/a)\}}{Ew(X)}<0,
\]
where at $a=0$ the numerator is interpreted as $cE[w(X)|v_0'X|]$. Frozen lag regressors have a positive density, so $P(v'X=0)=0$ for each nonzero $v$ (a pure intercept is a nonzero constant). The expression is strictly negative uniformly over the compact sphere and $a\in[0,C]$.

The required convergence is uniform: outside a compact raw-data set the score is bounded by $c w\|X\|$, and on that set its only limiting discontinuity, at $a=0$, is confined to an arbitrarily thin hyperplane neighborhood of vanishing probability. A finite-net argument completes this step. This contradicts the pilot score equation and proves coefficient tightness.

\textit{Pilot consistency.} On compact coefficient sets, frozen residual distributions have uniquely crossed medians and positive MADs. Uniform empirical convergence bounds the pilot scale away from zero and infinity. For each such scale $\sigma$, the conditional score equation identifies $\beta_j^H(u,\sigma)$: positive conditional density makes the conditional Huber contrast strictly positive at every nonzero shift. Nondegenerate regressors give coefficient identification. Continuity and compactness yield uniform separation away from this continuous family of minimizers. Convexity, the uniform contrast law and the approximate pilot score therefore give
$\widetilde\beta-\beta_j^H(u,\widetilde s)=o_p(1)$.

The residual law is consequently the frozen innovation law translated by $-\kappa_j(u,\widetilde s)$, up to a vanishing empirical perturbation. Translation invariance of MAD gives $\widetilde s\to_p s_{j0}(u)$; continuity of $\kappa_j$ then gives $\widetilde\beta\to_p\beta_j^H(u)$. All arguments are uniform on calibration compacts. Tightness of the calibration and its inverse permits exhaustion of these compacts without conditioning on its realized value.
\end{proof}

\begin{theorem}[Consistency of preliminary local dynamics]
\label{prop-final-local-consistency}
Under Lemma~\ref{lem-final-pilot}, let $\widehat b_\gamma$ minimize the preliminary dynamic criterion with $P=\operatorname{diag}(0_q,I_q)$ and $\lambda_n=\gamma\sqrt{q/n}$. Normalized objective errors may be $o_p(1)$. Uniformly over any fixed finite grid $\Gamma\subset[0,\infty)$, endpoint levels converge to $\beta_j^H(u)$ and scaled time slopes converge to zero. After back transformation, the lag-coefficient levels converge to $A_\ell(u)$. This holds for any data-selected $\widehat\gamma\in\Gamma$, including the initial-history validation rule. Fitted residual MADs and Kendall coefficients converge to the endpoint innovation functionals. Coordinate back transformation gives consistent dynamic and rank-scatter inputs in observed units.
\end{theorem}
\begin{proof}
\textit{Coefficients.} In calibrated coordinates put $b_0=(\beta_j^H(u)',0')'$. The limiting centered contrast is the frozen Huber contrast integrated over $z\in[-1,0]$. If $b-b_0=(v_0',v_1')'\ne0$, then $v_0+zv_1$ is nonzero except possibly at a single $z$, or is a nonzero constant vector. Conditional score centering makes the same Huber location apply at every frozen history and every $z$. Positive conditional density and the Huber-location equation therefore identify $b_0$ at scale $s_{j0}(u)$. Continuity gives uniform positive separation on every fixed sphere around it for scales in a sufficiently small neighborhood of $s_{j0}(u)$. Pilot-scale consistency places the fitted scale in that neighborhood with probability tending to one.

The penalty contrast on compact sets is $O(\max_{\gamma\in\Gamma}\lambda_n)=o(1)$. The local uniform law and convexity confine every approximate minimizer to each such sphere, proving consistency. Taking the maximum over the finite grid permits selection without requiring a consistent or uniquely optimal validation choice.

\textit{Residuals and rank scatter.} Compare the fitted residual with the frozen endpoint innovation shifted by $-\kappa_j(u,s_{j0}(u))$. The varying-process approximation vanishes in empirical probability by the coupling bound. The remaining coefficient perturbation satisfies, for any $M,\eta>0$,
\[
\frac1n\sum_s{\bf1}\{|d_s'(\widehat b-b_0)|>\eta\}
\le {\bf1}\{\|\widehat b-b_0\|>\eta/(CM)\}
 +\frac1n\sum_s{\bf1}\{\|X_s\|>M\}.
\]
The local empirical law makes the second term small as $M\to\infty$; consistency handles the first. Thus the residual empirical law converges to the translated endpoint innovation law. Unique quantile crossings give the shifted residual median and the unchanged innovation MAD. The product-empirical-measure argument in Section~\ref{app-consistency} gives convergence to the unchanged Kendall coefficients, including vanishing tied-pair fractions. Sine continuity and nonexpansive correlation projection give convergence of the preliminary rank-scatter input. Uniformity over calibration compacts and back transformation yield the result in original coordinates.
\end{proof}

\subsubsection{Propagation to the allocation target}
\label{app-final-continuity}
On a compact coefficient and scatter set with positive marginal scales and row denominators bounded away from zero, finite-horizon RAC is locally Lipschitz in coefficients, marginal scales and Kendall inputs. For a constant $C_H<\infty$,
\[
\|\widehat C-C\|
\le C_H\left\{\sum_{\ell=1}^p\|\widehat A_\ell-A_\ell\|
 +\|\widehat d-d\|+\|\widehat\tau-\tau\|\right\}.
\]
Finite-order response recursions are polynomials; the sine map and correlation projection are Lipschitz; square roots are Lipschitz away from zero; and $||x|-|y||\le|x-y|$. A quotient bound completes the argument. The bound also holds at zero response entries.

\subsection{Consistency of the complete RAC estimator}
The remaining adjustments affect coefficients, residual selection and correlation shrinkage. The proof bounds each adjustment and carries preliminary consistency through to the RAC allocation.
\begin{proposition}[Consistency of the RAC allocation estimator]\label{prop-regularized-target}
Under the fixed-dimensional local consistency conditions for the preliminary local estimator, the preceding regularized implementation satisfies
\[
C(\widetilde A,\widetilde\Omega,r;H)
\xrightarrow{p}C(A(u),\Omega^\dagger(u),r;H),\qquad r\in\{1,2\},
\]
for fixed $H$. Here the coefficient sequence and deletion fraction obey $c_n,q_n\to0$, pilot scales and the calibration matrix and its inverse are bounded in probability, and numerical optimization and projection errors vanish.
\end{proposition}
\begin{proof}
\textit{Coefficient adjustment.} The coefficient bound gives $\|\widetilde A-\widehat A\|=O_p(c_n)$ in observed units. For any $M,\eta>0$, the fraction of residual discrepancies above $\eta$ is bounded by
\[
\mathbf1\{\|\widetilde A-\widehat A\|>\eta/M\}
+\frac1n\sum_s\mathbf1\{\|y_{s-1}\|>M\}.
\]
The local tightness and empirical-frequency conditions used in the feasible-consistency proof let $n$ tend to infinity and then $M$ tend to infinity. Thus the perturbation affects a vanishing fraction of residuals at any fixed tolerance. The implementation omits the fitted intercept when forming residuals, adding a common translation to the Huber-centered residuals. MADs and Kendall coefficients are unchanged. Local intercept and coefficient variation vanish under the same local conditions. Unique median and MAD crossings and continuous marginals then give consistent residual quantiles and Kendall coefficients without finite first or second lag moments.

\textit{Residual deletion.} Deleting at most $q_nn$ rows changes each empirical distribution and each bounded pairwise-kernel average by $O(q_n)$. Hence $D$ and $U$ retain their population limits.

\textit{Correlation shrinkage.} Bounded kernels give $V_N=O(K^2/N)$ deterministically. Projection is nonexpansive and $\Pi_{\mathcal C}(U)=U$, so
\[
\|\widetilde R-U\|_F
\le(1-\alpha)d
=\min\{d,V_N/d\}\le\sqrt{V_N},
\]
with zero difference at $d=0$. This bound holds for any historically selected target. Consequently $\widetilde\Omega\to_p\Omega^\dagger(u)$. Finite-horizon responses are continuous functions of coefficients and scatter, and the positive own impact keeps each allocation denominator positive. Continuous mapping proves the result.
\end{proof}

\subsection{Numerical implementation}\label{app-implementation}
The nearest-correlation calculation uses eigenvalue tolerance $10^{-10}$.
Current-window scatter and shrinkage use the empirical sign-pair Kendall coefficient. Historical template calculations use tie-adjusted $\tau_b$. Under continuous innovations and vanishing residual perturbations, empirical tied-pair fractions vanish, and both normalizations converge to the population sign kernel.

\FloatBarrier

\section{Population stability and priority recovery}\label{app-finite-sample}

\subsection{Sample-size behavior at the covariance boundary}\label{app-suite-boundary}
Known dynamics and a common squared-response rule isolate the effect of estimating the shock input. We fix
\[
A=\begin{pmatrix}.2&.25&0\\0&.4&.15\\.1&0&.3\end{pmatrix},\qquad
S=\begin{pmatrix}1&.35&.15\\.35&1&.25\\.15&.25&1\end{pmatrix}.
\]
For each $\nu\in\{.8,1.5,2.1,3,4.1,8\}$ and a Gaussian benchmark, 300 independent innovation samples contain 8,000 observations. Nested prefixes use $n=125,500,2000,8000$; dynamics and $H=10$ are known. A rank-one control replaces each innovation vector by its first coordinate times $(1,.7,1.3)'$. The resulting normalized shares equal $1/3$ under either scatter at every sample size.

At the lowest tail indices, covariance-share dispersion remains substantial as the sample expands, while rank-scatter dispersion contracts.

\begin{table}[!htbp]\centering\small\renewcommand{\arraystretch}{1.05}\setlength{\tabcolsep}{5pt}
\caption{Quadratic-share dispersion and common-target recovery with known dynamics.}\label{tbl-suite-boundary}
\begin{tabular}{lrrrrr}\toprule
Innovations & $n$ & Cov. share SD & Rank share SD & Cov. loss & Rank loss\\\midrule

0.8 & 125 & 11.81 & 5.57 & --- & 7.29\\
0.8 & 500 & 11.05 & 2.77 & --- & 3.76\\
0.8 & 2,000 & 11.67 & 1.37 & --- & 1.95\\
0.8 & 8,000 & 11.89 & 0.70 & --- & 0.98\\
\addlinespace[3pt]
1.5 & 125 & 10.16 & 4.89 & --- & 6.49\\
1.5 & 500 & 9.87 & 2.57 & --- & 3.44\\
1.5 & 2,000 & 10.87 & 1.26 & --- & 1.72\\
1.5 & 8,000 & 10.19 & 0.64 & --- & 0.85\\
\addlinespace[3pt]
2.1 & 125 & 8.74 & 4.59 & 13.96 & 6.31\\
2.1 & 500 & 7.66 & 2.40 & 11.73 & 3.41\\
2.1 & 2,000 & 6.65 & 1.27 & 9.91 & 1.73\\
2.1 & 8,000 & 6.50 & 0.58 & 9.03 & 0.84\\
\addlinespace[3pt]
3.0 & 125 & 7.04 & 4.63 & 9.83 & 6.04\\
3.0 & 500 & 4.62 & 2.29 & 6.63 & 3.17\\
3.0 & 2,000 & 3.54 & 1.12 & 4.81 & 1.59\\
3.0 & 8,000 & 2.64 & 0.61 & 3.17 & 0.82\\
\addlinespace[3pt]
4.1 & 125 & 5.51 & 4.28 & 7.64 & 5.89\\
4.1 & 500 & 3.23 & 2.12 & 4.57 & 2.95\\
4.1 & 2,000 & 1.90 & 1.16 & 2.59 & 1.51\\
4.1 & 8,000 & 1.01 & 0.54 & 1.38 & 0.78\\
\addlinespace[3pt]
8.0 & 125 & 4.08 & 3.96 & 5.88 & 5.64\\
8.0 & 500 & 2.26 & 2.22 & 3.21 & 3.01\\
8.0 & 2,000 & 1.09 & 1.06 & 1.66 & 1.52\\
8.0 & 8,000 & 0.58 & 0.53 & 0.83 & 0.75\\
\addlinespace[3pt]
Gaussian & 125 & 3.60 & 3.88 & 5.18 & 5.44\\
Gaussian & 500 & 1.89 & 2.14 & 2.59 & 2.76\\
Gaussian & 2,000 & 0.93 & 1.02 & 1.32 & 1.40\\
Gaussian & 8,000 & 0.46 & 0.52 & 0.62 & 0.66\\
\bottomrule\end{tabular}
\par\smallskip\begin{minipage}{.97\linewidth}\footnotesize\raggedright\textit{Notes.}   Dispersion and loss are in percentage points, computed across 300 independent samples per innovation law. Both methods use squared responses. Loss is the summed absolute off-diagonal share error per recipient. Share SD concerns $C_{12}$. Loss targets the common elliptical squared allocation when covariance exists; below that boundary the rank target remains defined and the population covariance loss is omitted. Every rank-one cancellation control has $C_{12}=1/3$ to numerical precision.\end{minipage}
\end{table}

\subsection{Heterogeneous network designs}\label{app-heterogeneous-designs}
The factor and directed-block systems each contain forty variables, with $S=DRD$ and $D_{ii}=\exp\{-0.3+0.6(i-1)/39\}$.

In the factor network, $R=\operatorname{diag}(1-\ell_i^2)+\ell\ell'$, $\ell_i=0.05+0.8(i-1)/39$; dynamic diagonal entries range from 0.15 to 0.35 and off-diagonal entries are $0.08v_iv_j/39$, with $v_i=0.7+0.6(i-1)/39$.

In the block network, $R=LL'+\operatorname{diag}(1-\operatorname{diag}(LL'))$: all nodes load 0.15 on a common factor, and each consecutive block of ten loads from 0.15 to 0.8 on its own factor. Its dynamic diagonal ranges from 0.15 to 0.45, the first subdiagonal is 0.2, and $A_{11,1}=A_{21,11}=A_{31,21}=0.15$.

Innovations are multivariate Student draws with shape $S$. Population amplitude recipient-strength ranges are 24.87 and 21.57 percentage points, respectively.

\subsection{Recipient and source recovery with estimated dynamics}\label{app-suite-pipeline}
In the factor and directed-block systems of Section~\ref{app-heterogeneous-designs}, share error measures how closely fitted external shares match their population values; list recovery measures whether the same recipients are selected. Stationary and locally varying dynamics combine with elliptical Student laws at $\nu=.8,1.5,3,8$, a Gaussian law, or nonelliptical independent Student coordinates with unequal tail indices. The 24 cells contain 100 independent paths each. Every path has 500 burn-in and 3,000 retained observations. Each final endpoint is evaluated at $m=125,250,500$.

Let $D=\operatorname{diag}(s_i)$ be the endpoint marginal scale, $R_0=D^{-1}S_0D^{-1}$, $v_i$ equally spaced on $[-.5,.5]$, and $u_t=(t-2999)/4000$. In the time-varying design,
\[
A_t=A_0+.06u_tI,\quad D_t=\operatorname{diag}\{s_i\exp(.5u_tv_i)\},\quad
R_t=(1+.25u_t)R_0-.25u_tI.
\]
Elliptical innovations use shape $D_tR_tD_t$ and one common Student radial scale per vector. Student draws are divided by $1.4826t_\nu^{-1}(.75)$; Gaussian draws use the corresponding normal quantile. The independent-margin design uses $R_t=I$, with coordinate tail indices cycling through $(.8,1.5,3,8)$ in node order and marginal normalization at each coordinate's own degrees of freedom. Its endpoint rank scatter is exactly $D^2$, despite unequal tail indices. Stationary designs set $u_t=0$.

\paragraph*{Estimation and comparison.}
Each path uses its own initial 500 observations for coordinate calibration and slope-penalty selection. Every fitted RAC coefficient matrix has spectral radius below one. Historical templates use up to 1,500 pre-window observations. The same fitted dynamics generate all four scatter--aggregation combinations. Local-linear OLS and sample covariance supply the conventional dynamic benchmark for squared-response recovery. Each allocation is assessed against its own population target; common-target comparisons use finite-covariance elliptical squared allocations.

\paragraph*{Recipient priorities.}
For population receiving strength $f_i$, let $I_q$ and $\widehat I_q$ be its population and estimated top-$q$ sets. The selection loss is
\[
L_q=\frac1q\left(\sum_{i\in I_q}f_i-\sum_{i\in\widehat I_q}f_i\right).
\]
The primary display uses $q=10$ for forty recipients. List overlap records the names recovered, while $L_q$ measures the external share forgone. Table~\ref{tbl-suite-pipeline} reports both alongside share error and rank agreement. Tables~\ref{tbl-suite-windows} and \ref{tbl-suite-classification} follow share recovery and source classification as the current window changes.
\par
\begin{table}[!htbp]\centering\small\renewcommand{\arraystretch}{1.05}\setlength{\tabcolsep}{4pt}
\caption{Recipient priorities with estimated dynamics, $m=250$.}\label{tbl-suite-pipeline}
\begin{tabular}{lllrrrr}\toprule
Network & Dynamics & Law & Loss (pp) & Share MAE (pp) & Rank corr. & Top 10 (\%)\\\midrule

Factor & Stationary & 0.8 & 0.10 & 2.35 (0.04) & 0.946 & 87.6\\
Factor & Stationary & 1.5 & 0.17 & 2.24 (0.04) & 0.928 & 84.4\\
Factor & Stationary & 3 & 0.25 & 2.36 (0.04) & 0.911 & 81.1\\
Factor & Stationary & 8 & 0.22 & 2.26 (0.04) & 0.915 & 81.5\\
Factor & Stationary & Gaussian & 0.21 & 2.29 (0.04) & 0.918 & 82.7\\
Factor & Stationary & Mixed & 0.18 & 5.98 (0.23) & 0.059 & 27.8\\
\addlinespace[3pt]
Factor & Time varying & 0.8 & 0.11 & 2.28 (0.04) & 0.946 & 87.3\\
Factor & Time varying & 1.5 & 0.18 & 2.33 (0.04) & 0.925 & 83.4\\
Factor & Time varying & 3 & 0.18 & 2.39 (0.04) & 0.914 & 83.4\\
Factor & Time varying & 8 & 0.21 & 2.31 (0.05) & 0.920 & 83.8\\
Factor & Time varying & Gaussian & 0.19 & 2.25 (0.04) & 0.928 & 83.5\\
Factor & Time varying & Mixed & 0.17 & 5.44 (0.23) & 0.064 & 26.3\\
\addlinespace[3pt]
Blocks & Stationary & 0.8 & 1.14 & 2.85 (0.11) & 0.592 & 49.7\\
Blocks & Stationary & 1.5 & 0.86 & 2.55 (0.08) & 0.662 & 54.4\\
Blocks & Stationary & 3 & 1.01 & 2.56 (0.06) & 0.637 & 51.3\\
Blocks & Stationary & 8 & 0.91 & 2.55 (0.07) & 0.653 & 52.9\\
Blocks & Stationary & Gaussian & 0.96 & 2.58 (0.06) & 0.641 & 52.3\\
Blocks & Stationary & Mixed & 1.51 & 5.35 (0.27) & 0.722 & 63.6\\
\addlinespace[3pt]
Blocks & Time varying & 0.8 & 1.11 & 2.73 (0.09) & 0.604 & 49.2\\
Blocks & Time varying & 1.5 & 0.97 & 2.73 (0.09) & 0.636 & 50.9\\
Blocks & Time varying & 3 & 0.95 & 2.59 (0.07) & 0.634 & 52.5\\
Blocks & Time varying & 8 & 1.05 & 2.67 (0.06) & 0.621 & 52.0\\
Blocks & Time varying & Gaussian & 0.98 & 2.56 (0.06) & 0.640 & 50.3\\
Blocks & Time varying & Mixed & 1.40 & 4.42 (0.17) & 0.741 & 65.7\\
\bottomrule\end{tabular}
\par\smallskip\begin{minipage}{.98\linewidth}\footnotesize\raggedright\textit{Notes.} Every cell contains 100 independent paths. Share MAE and Top 10 population-strength shortfall (Loss) are percentage points; parentheses give the MAE Monte Carlo standard error. Mixed denotes independent Student coordinates with indices cycling through $.8,1.5,3,8$. All quantities are evaluated against the amplitude population target. Read rank correlation and list overlap alongside the population separation in Table~\ref{tbl-identification-geometry}; Share MAE measures error in recipient-share levels.\end{minipage}
\end{table}

\par

The cutoff gaps and score ranges in Table~\ref{tbl-identification-geometry} determine how sharply the population distinguishes the selected recipients.

In the elliptical designs, each population recipient is externally dominated; in the independent-margin designs, each is own dominated. Table~\ref{tbl-suite-classification} therefore measures preservation of a common source category within each design. The extended designs below examine rare external recipients and classifications closer to the boundary.
\par
\begin{table}[!htbp]\centering\small\renewcommand{\arraystretch}{1.05}\setlength{\tabcolsep}{4pt}
\caption{Population separation and source categories in the recovery experiments.}\label{tbl-identification-geometry}
\begin{tabular}{llrrrr}\toprule
Network & Innovations & Range & Cutoff gap & Random loss & External (\%)\\\midrule
Factor & Elliptical & 24.871 & 0.178593 & 5.507 & 100\\
Factor & Independent margins & 0.667 & 0.000612 & 0.192 & 0\\
Blocks & Elliptical & 21.568 & 0.058817 & 3.117 & 100\\
Blocks & Independent margins & 42.868 & 0.190471 & 6.704 & 0\\
\bottomrule\end{tabular}
\par\smallskip\begin{minipage}{.98\linewidth}\footnotesize\raggedright\textit{Notes.} Scores, gaps and losses are percentage points. Cutoff gap is the difference between the tenth and eleventh population recipient scores. Random loss is the Top 10 mean minus the mean over all forty recipients, the expected loss of a uniformly selected ten-name list. Elliptical laws share the same endpoint target, as do stationary and time-varying designs. External gives the fraction of population recipients above 50\%.\end{minipage}
\end{table}

\par

\par
\begin{table}[!htbp]\centering\small\renewcommand{\arraystretch}{1.05}\setlength{\tabcolsep}{4pt}
\caption{Recipient-share recovery across window lengths.}\label{tbl-suite-windows}
\begin{tabular}{lllrrr}\toprule
Network & Dynamics & Law & $m=125$ & $m=250$ & $m=500$\\\midrule

Factor & Stationary & 0.8 & 3.90 (0.05) & 2.35 (0.04) & 1.66 (0.02)\\
Factor & Stationary & 1.5 & 3.28 (0.04) & 2.24 (0.04) & 1.71 (0.03)\\
Factor & Stationary & 3 & 3.25 (0.04) & 2.36 (0.04) & 1.84 (0.03)\\
Factor & Stationary & 8 & 3.10 (0.04) & 2.26 (0.04) & 1.84 (0.04)\\
Factor & Stationary & Gaussian & 3.01 (0.05) & 2.29 (0.04) & 1.89 (0.04)\\
Factor & Stationary & Mixed & 8.80 (0.32) & 5.98 (0.23) & 4.86 (0.22)\\
\addlinespace[3pt]
Factor & Time varying & 0.8 & 3.60 (0.05) & 2.28 (0.04) & 1.66 (0.03)\\
Factor & Time varying & 1.5 & 3.31 (0.05) & 2.33 (0.04) & 1.72 (0.03)\\
Factor & Time varying & 3 & 3.26 (0.05) & 2.39 (0.04) & 1.81 (0.03)\\
Factor & Time varying & 8 & 3.03 (0.05) & 2.31 (0.05) & 1.88 (0.04)\\
Factor & Time varying & Gaussian & 2.85 (0.04) & 2.25 (0.04) & 1.94 (0.04)\\
Factor & Time varying & Mixed & 7.86 (0.28) & 5.44 (0.23) & 4.50 (0.22)\\
\addlinespace[3pt]
Blocks & Stationary & 0.8 & 4.89 (0.15) & 2.85 (0.11) & 2.34 (0.07)\\
Blocks & Stationary & 1.5 & 3.71 (0.10) & 2.55 (0.08) & 2.07 (0.06)\\
Blocks & Stationary & 3 & 3.55 (0.08) & 2.56 (0.06) & 2.13 (0.05)\\
Blocks & Stationary & 8 & 3.66 (0.08) & 2.55 (0.07) & 1.95 (0.05)\\
Blocks & Stationary & Gaussian & 3.69 (0.09) & 2.58 (0.06) & 1.93 (0.05)\\
Blocks & Stationary & Mixed & 7.59 (0.27) & 5.35 (0.27) & 3.16 (0.15)\\
\addlinespace[3pt]
Blocks & Time varying & 0.8 & 4.81 (0.16) & 2.73 (0.09) & 2.21 (0.07)\\
Blocks & Time varying & 1.5 & 3.52 (0.10) & 2.73 (0.09) & 2.16 (0.06)\\
Blocks & Time varying & 3 & 3.56 (0.09) & 2.59 (0.07) & 1.95 (0.05)\\
Blocks & Time varying & 8 & 3.53 (0.08) & 2.67 (0.06) & 1.97 (0.05)\\
Blocks & Time varying & Gaussian & 3.58 (0.08) & 2.56 (0.06) & 1.87 (0.05)\\
Blocks & Time varying & Mixed & 7.44 (0.27) & 4.42 (0.17) & 3.32 (0.17)\\
\bottomrule\end{tabular}
\par\smallskip\begin{minipage}{.98\linewidth}\footnotesize\raggedright\textit{Notes.} Entries are mean absolute recipient-share errors in percentage points, with Monte Carlo standard errors in parentheses. Each cell uses the same 100 paths at all three windows; all tuning is based on the initial history.\end{minipage}
\end{table}

\par

\par
\begin{table}[!htbp]\centering\small\renewcommand{\arraystretch}{1.05}\setlength{\tabcolsep}{4pt}
\caption{Own/external classification across window lengths.}\label{tbl-suite-classification}
\begin{tabular}{lllrrr}\toprule
Network & Dynamics & Law & $m=125$ & $m=250$ & $m=500$\\\midrule

Factor & Stationary & 0.8 & 100.0 & 100.0 & 100.0\\
Factor & Stationary & 1.5 & 100.0 & 100.0 & 100.0\\
Factor & Stationary & 3 & 100.0 & 100.0 & 100.0\\
Factor & Stationary & 8 & 100.0 & 100.0 & 100.0\\
Factor & Stationary & Gaussian & 100.0 & 100.0 & 100.0\\
Factor & Stationary & Mixed & 99.8 & 100.0 & 100.0\\
\addlinespace[3pt]
Factor & Time varying & 0.8 & 100.0 & 100.0 & 100.0\\
Factor & Time varying & 1.5 & 100.0 & 100.0 & 100.0\\
Factor & Time varying & 3 & 100.0 & 100.0 & 100.0\\
Factor & Time varying & 8 & 100.0 & 100.0 & 100.0\\
Factor & Time varying & Gaussian & 100.0 & 100.0 & 100.0\\
Factor & Time varying & Mixed & 99.7 & 100.0 & 100.0\\
\addlinespace[3pt]
Blocks & Stationary & 0.8 & 100.0 & 100.0 & 100.0\\
Blocks & Stationary & 1.5 & 100.0 & 100.0 & 100.0\\
Blocks & Stationary & 3 & 100.0 & 100.0 & 100.0\\
Blocks & Stationary & 8 & 100.0 & 100.0 & 100.0\\
Blocks & Stationary & Gaussian & 100.0 & 100.0 & 100.0\\
Blocks & Stationary & Mixed & 94.9 & 98.0 & 99.2\\
\addlinespace[3pt]
Blocks & Time varying & 0.8 & 100.0 & 100.0 & 100.0\\
Blocks & Time varying & 1.5 & 100.0 & 100.0 & 100.0\\
Blocks & Time varying & 3 & 100.0 & 100.0 & 100.0\\
Blocks & Time varying & 8 & 100.0 & 100.0 & 100.0\\
Blocks & Time varying & Gaussian & 100.0 & 100.0 & 100.0\\
Blocks & Time varying & Mixed & 95.1 & 98.3 & 99.1\\
\bottomrule\end{tabular}
\par\smallskip\begin{minipage}{.98\linewidth}\footnotesize\raggedright\textit{Notes.} Entries are percentages of all recipients whose fitted and population external shares lie on the same side of 50\%. All 100 paths and all forty recipients enter every cell.\end{minipage}
\end{table}

\par

\begin{figure}[!htbp]\centering
\includegraphics[width=\linewidth]{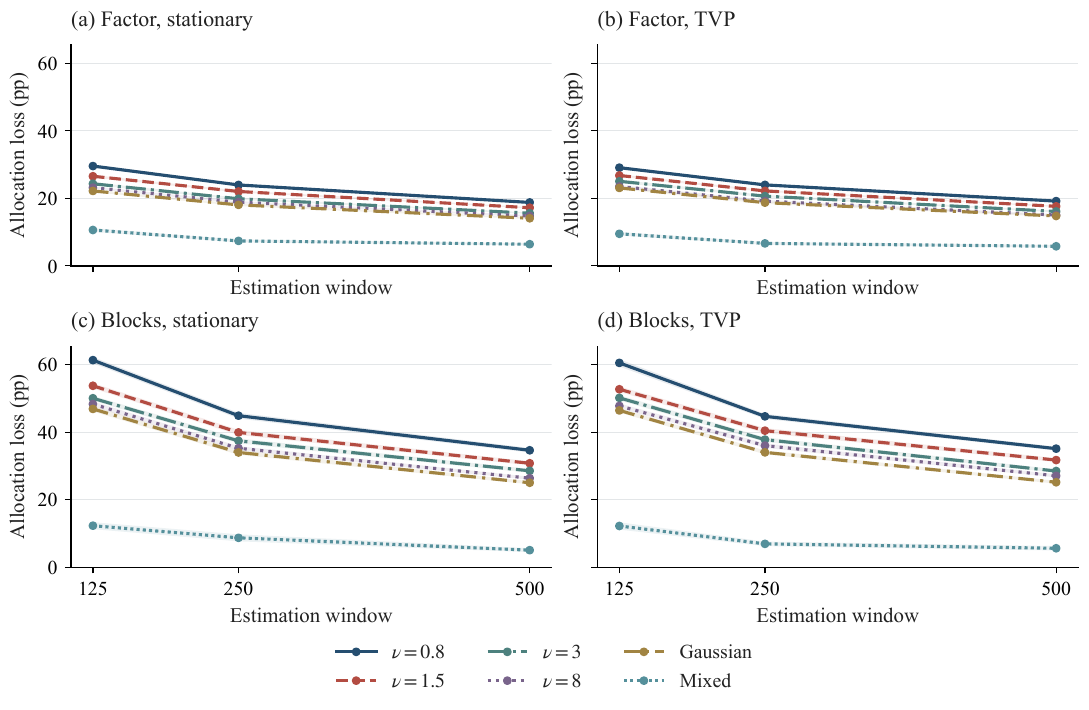}
\caption{\textbf{Allocation recovery across network structure, local history and innovation law.} Curves give mean off-diagonal amplitude allocation error per recipient; shaded bands are $\pm1.96$ Monte Carlo standard errors. Each cell contains 100 paths. The same paths and initial tuning enter all windows.}\label{fig-suite-recovery}
\end{figure}

\FloatBarrier
\subsection{Propagation structure and available history}\label{app-regularized-simulation}
Twenty-six designs extend the factor and block comparisons to null, dense, signed, star and weakly correlated systems. Holding the current-window length fixed, we compare recipient-share error, source classification and the selected Top 10 list's population-strength shortfall at endpoints with increasing available history. Each design uses twelve independent paths, 500 burn-in observations and 3,000 retained observations. Each path uses its own initial 500 observations for calibration and tuning. The current window is $m=250$, with $H=10$, at endpoints 501, 1,001 and 2,999. Available pre-window history expands from 251 to 751 to 1,499 adjacent pairs.

At the final endpoint, mean recipient-share error is 4.44 percentage points, Top 10 population-strength shortfall is 0.75 points, and own/external classification accuracy is 97.7\%. Table~\ref{tbl-causal-point-stages} holds paths fixed while historical information expands. Table~\ref{tbl-causal-point-cases} then separates the propagation designs at the final endpoint.

\paragraph*{Propagation structures.}
The designs span independent, factor, block, dense and star propagation. The factor and block coefficients use Section~\ref{app-heterogeneous-designs}. Independent-block variants set the innovation correlation to identity. Null dynamics are $.5I$. Dense dynamics have diagonal $.2$ and off-diagonal $.6/39$; signed dense dynamics conjugate this matrix by alternating signs.

Star dynamics are $.5I+ae_1(\mathbf1-e_1)'$, with $a$ set so recipient 1 has population external amplitude share $.4,.6$ or $.8$. The channel-reversal design embeds the main text's three-node Gaussian construction in $.2I_{40}$. The four common-correlation designs use $.5I$ and $R=(1-\rho)I+\rho gg'$, with $\rho=.01,.03,.06$ and $g=\mathbf1$, or $\rho=.03$ and alternating signs.

\paragraph*{Innovation laws.}
Gaussian and mixed-margin innovations share their specified copula correlation. The mixed-margin law is nonelliptical: Student indices cycle through $.8,1.5,3,8$, with unit scaled MAD normalization before applying marginal scales. Independent direct-Student factor and block designs reuse Section~\ref{app-suite-pipeline}'s innovation construction with fresh paths. Factor and block scales follow their stated scatter inputs; null and dense scales are $\exp(v_i)$ for $v_i$ equally spaced on $[-.3,.3]$. Star, channel-reversal and common-correlation scales are one. Time-varying cases use $A_t=A_0+.06u_tI$, $u_t=(t-2999)/4000$; direct-Student cases also vary scales as in Section~\ref{app-suite-pipeline}. Each endpoint is evaluated against its contemporaneous population target.

\begin{table}[!htbp]\centering\small\renewcommand{\arraystretch}{1.05}\setlength{\tabcolsep}{4pt}
\caption{RAC accuracy as historical information expands.}\label{tbl-causal-point-stages}
\begin{tabular}{lrrrrr}\toprule
Stage & History pairs & TCI error & Share MAE & Top 10 loss & Class. (\%)\\\midrule
Early & 251 & 5.39 (0.36) & 6.99 (0.34) & 0.93 (0.04) & 93.5\\
Expanding & 751 & 3.64 (0.17) & 4.70 (0.18) & 0.85 (0.04) & 97.0\\
Mature & 1499 & 3.30 (0.10) & 4.44 (0.10) & 0.75 (0.03) & 97.7\\\bottomrule\end{tabular}
\par\smallskip\begin{minipage}{.98\linewidth}\footnotesize\raggedright\textit{Notes.} Each stage contains 26 designs with 12 paths per design, $K=40$, $m=250$ and $H=10$. Endpoints are positions 501, 1,001 and 2,999 (zero based). Means weight designs equally. Parentheses give Monte Carlo standard errors across replication groups. Errors and selection losses are percentage points; classification uses all recipients. The same paths supply all three stages.\end{minipage}
\end{table}

\begin{table}[!htbp]\centering\small\renewcommand{\arraystretch}{1.05}\setlength{\tabcolsep}{4pt}
\caption{RAC recipient accuracy across propagation designs.}\label{tbl-causal-point-cases}
\begin{tabular}{lrrr}\toprule
Design & Share MAE & Top 10 loss & Class. (\%)\\\midrule

\addlinespace[3pt]
\multicolumn{4}{l}{\textit{Baseline and clustered propagation}}\\
Null (G) & 5.11 (1.55) & 0.00 (0.00) & 98.1\\
Null (M) & 2.23 (1.17) & 0.00 (0.00) & 100.0\\
Factor (G) & 2.12 (0.10) & 0.26 (0.05) & 100.0\\
Factor (M) & 2.15 (0.09) & 0.10 (0.02) & 100.0\\
Blocks (G) & 2.73 (0.29) & 0.88 (0.13) & 100.0\\
Blocks (M) & 2.19 (0.11) & 0.81 (0.09) & 100.0\\
Blocks, independent (G) & 6.18 (0.35) & 2.63 (0.60) & 98.8\\
Blocks, independent (M) & 5.49 (0.87) & 1.62 (0.17) & 97.7\\
\addlinespace[3pt]
\multicolumn{4}{l}{\textit{Dense propagation}}\\
Dense (G) & 8.42 (0.76) & 2.83 (0.21) & 91.0\\
Dense (M) & 7.38 (1.21) & 1.66 (0.30) & 100.0\\
Signed dense (G) & 8.93 (0.70) & 3.16 (0.25) & 89.0\\
\addlinespace[3pt]
\multicolumn{4}{l}{\textit{Star and reversal systems}}\\
Star 40 (G) & 4.34 (1.33) & 1.00 (0.52) & 99.2\\
Star 60 (G) & 4.10 (0.78) & 0.00 (0.00) & 97.5\\
Star 80 (G) & 7.88 (1.23) & 0.00 (0.00) & 97.7\\
Star 60 (M) & 2.45 (1.27) & 0.00 (0.00) & 100.0\\
Channel reversal (G) & 4.81 (1.49) & 0.00 (0.00) & 100.0\\
\addlinespace[3pt]
\multicolumn{4}{l}{\textit{Changing dynamics and independent margins}}\\
Varying factor (G) & 2.32 (0.09) & 0.19 (0.05) & 100.0\\
Varying blocks (M) & 2.36 (0.18) & 0.79 (0.13) & 100.0\\
Factor, independent (T) & 5.50 (0.47) & 0.20 (0.01) & 100.0\\
Varying factor, indep. (T) & 4.95 (0.57) & 0.18 (0.01) & 100.0\\
Blocks, independent (T) & 4.67 (0.44) & 1.53 (0.20) & 98.3\\
Varying blocks, indep. (T) & 4.17 (0.26) & 1.74 (0.17) & 98.5\\
\addlinespace[3pt]
\multicolumn{4}{l}{\textit{Common-shock systems}}\\
Common $\rho=.01$ & 6.60 (1.46) & 0.00 (0.00) & 97.3\\
Common $\rho=.03$ & 3.02 (0.60) & 0.00 (0.00) & 84.2\\
Common $\rho=.06$ & 2.33 (0.45) & 0.00 (0.00) & 100.0\\
Signed common $\rho=.03$ & 3.05 (0.56) & 0.00 (0.00) & 91.7\\
\bottomrule\end{tabular}
\par\smallskip\begin{minipage}{.98\linewidth}\footnotesize\raggedright\textit{Notes.} Each row includes all twelve paths at the final endpoint. Parentheses give Monte Carlo standard errors. G denotes Gaussian innovations; M uses a Gaussian copula with Student indices cycling through $.8,1.5,3,8$; T draws independent Student coordinates directly with the marginal scales of Section~\ref{app-suite-pipeline}. Classification includes every recipient. The two $\rho=.03$ designs place all population shares within five points of 50\%.\end{minipage}
\end{table}

\paragraph*{Classification by source category.}
Table~\ref{tbl-classification-resolution} separates externally and own-dominated population recipients. In the Gaussian star with focal share $.6$, the externally dominated recipient is recovered in eight of twelve paths (66.7\%), while 98.3\% of own-dominated recipients are classified correctly. Together, the two source categories give an aggregate accuracy of 97.5\%. At focal share $.8$, all twelve external recipients are recovered. The common-correlation and signed-common designs with $\rho=.03$ place every population share within five percentage points of the boundary; their classification rates are 84.2\% and 91.7\%, respectively. 
\begin{table}[!htbp]\centering\small\renewcommand{\arraystretch}{1.05}\setlength{\tabcolsep}{4pt}
\caption{Classification by population source category and boundary distance.}\label{tbl-classification-resolution}
\begin{tabular}{lrrrrrr}\toprule
& \multicolumn{2}{c}{External dominance} & \multicolumn{2}{c}{Own dominance} & \multicolumn{2}{c}{Near boundary}\\\cmidrule(lr){2-3}\cmidrule(lr){4-5}\cmidrule(lr){6-7} Design & $N_E$ & Correct (\%) & $N_O$ & Correct (\%) & $N_B$ & Correct (\%)\\\midrule

\addlinespace[3pt]
\multicolumn{7}{l}{\textit{Baseline and clustered propagation}}\\
Null (G) & 0 & --- & 480 & 98.1 & 0 & ---\\
Null (M) & 0 & --- & 480 & 100.0 & 0 & ---\\
Factor (G) & 480 & 100.0 & 0 & --- & 0 & ---\\
Factor (M) & 480 & 100.0 & 0 & --- & 0 & ---\\
Blocks (G) & 480 & 100.0 & 0 & --- & 0 & ---\\
Blocks (M) & 480 & 100.0 & 0 & --- & 0 & ---\\
Blocks, independent (G) & 0 & --- & 480 & 98.8 & 0 & ---\\
Blocks, independent (M) & 0 & --- & 480 & 97.7 & 0 & ---\\
\addlinespace[3pt]
\multicolumn{7}{l}{\textit{Dense propagation}}\\
Dense (G) & 480 & 91.0 & 0 & --- & 0 & ---\\
Dense (M) & 480 & 100.0 & 0 & --- & 0 & ---\\
Signed dense (G) & 480 & 89.0 & 0 & --- & 0 & ---\\
\addlinespace[3pt]
\multicolumn{7}{l}{\textit{Star and reversal systems}}\\
Star 40 (G) & 0 & --- & 480 & 99.2 & 0 & ---\\
Star 60 (G) & 12 & 66.7 & 468 & 98.3 & 0 & ---\\
Star 80 (G) & 12 & 100.0 & 468 & 97.6 & 0 & ---\\
Star 60 (M) & 12 & 100.0 & 468 & 100.0 & 0 & ---\\
Channel reversal (G) & 12 & 100.0 & 468 & 100.0 & 0 & ---\\
\addlinespace[3pt]
\multicolumn{7}{l}{\textit{Changing dynamics and independent margins}}\\
Varying factor (G) & 480 & 100.0 & 0 & --- & 0 & ---\\
Varying blocks (M) & 480 & 100.0 & 0 & --- & 0 & ---\\
Factor, independent (T) & 0 & --- & 480 & 100.0 & 0 & ---\\
Varying factor, indep. (T) & 0 & --- & 480 & 100.0 & 0 & ---\\
Blocks, independent (T) & 0 & --- & 480 & 98.3 & 0 & ---\\
Varying blocks, indep. (T) & 0 & --- & 480 & 98.5 & 0 & ---\\
\addlinespace[3pt]
\multicolumn{7}{l}{\textit{Common-shock systems}}\\
Common $\rho=.01$ & 0 & --- & 480 & 97.3 & 0 & ---\\
Common $\rho=.03$ & 480 & 84.2 & 0 & --- & 480 & 84.2\\
Common $\rho=.06$ & 480 & 100.0 & 0 & --- & 0 & ---\\
Signed common $\rho=.03$ & 480 & 91.7 & 0 & --- & 480 & 91.7\\
\bottomrule\end{tabular}
\par\smallskip\begin{minipage}{.98\linewidth}\footnotesize\raggedright\textit{Notes.} Each row pools forty recipients across twelve paths at the final endpoint. $E$ and $O$ denote population external and own dominance; $B$ denotes an absolute distance from 50\% of at most five percentage points. Boundary observations are a subset of the first two groups. A dash indicates an empty group. Design labels match Table~\ref{tbl-causal-point-cases}.\end{minipage}
\end{table}

\FloatBarrier

\section{Financial monitoring evidence}\label{app-monitoring}

\FloatBarrier
\subsection{Data, residuals and tail diagnostics}\label{app-tail-diagnostics}
The VOLARE data were downloaded on October 4, 2026, and all three panels end on September 30, 2026. The stock panel contains forty symbols across 2,953 dates and 118,120 observations, with no missing cells. Membership is fixed at the downloaded panel and does not track historical index membership. 

The currency panel has seven missing cells across six of 4,415 dates; the futures panel has 87 across 67 of 4,401 dates. Each market calendar is the union of its recorded dates. Fits use complete pairs adjacent on that calendar, without imputing missing values or joining observations across a missing date.
Tail diagnostics use endpoint residuals from the fitted coefficient levels. Deviations from each local residual median enter trailing windows on the original calendar, retaining finite-observation counts. 

\paragraph*{Network-level tail regimes.}
For each market and date, $B_t^{(2)}$ counts bias-corrected Hill estimates below two and $I_t^{(2)}$ records whether the count is positive. Network frequencies use dates with all $K$ estimates finite. Table~\ref{tbl-network-tail-grid} crosses windows of 60, 120 and 240 original positions with four tail fractions. Each cell uses its complete-date denominator.
\begin{table}[!htbp]\centering\small\renewcommand{\arraystretch}{1.12}\setlength{\tabcolsep}{5pt}
\caption{Network dates with at least one tail estimate below two.}\label{tbl-network-tail-grid}
\begin{tabular}{llrrrrr}\toprule
Market & $W$ & Complete dates & 5\% tail & 10\% tail & 15\% tail & 25\% tail\\\midrule
Stocks & 60 & 2,393 & 89.2 & 63.9 & 49.1 & 42.4\\
 & 120 & 2,333 & 59.5 & 18.9 & 9.6 & 7.2\\
 & 240 & 2,213 & 18.5 & 1.6 & 0.0 & 0.1\\
Currencies & 60 & 3,742--3,855 & 18.0 & 10.5 & 8.7 & 8.5\\
 & 120 & 3,795 & 5.6 & 0.5 & 2.8 & 0.6\\
 & 240 & 3,675 & 0.0 & 0.0 & 0.0 & 0.0\\
Futures & 60 & 1,628--3,841 & 36.0 & 34.6 & 31.8 & 25.5\\
 & 120 & 3,781 & 27.5 & 19.0 & 16.2 & 13.8\\
 & 240 & 3,661 & 9.4 & 13.2 & 6.3 & 2.6\\\bottomrule\end{tabular}
\par\smallskip\begin{minipage}{.98\linewidth}\footnotesize\raggedright\textit{Notes.} Cells give percentages of complete network dates under bias-corrected Hill. Dates may vary across fractions because a valid estimate requires at least three upper order statistics. At $W=60$, the currency 5\% cell uses 3,742 dates and its other cells 3,855; the corresponding futures counts are 1,628 and 3,841. At longer windows every fraction has the same dates. The full grid, including uncorrected Hill, is in the replication materials.\end{minipage}
\end{table}

\paragraph*{Upper-bound classifications.}
We assess the primary $W=120$, 10\%-tail flags with moving-block resampling \citep{kunsch1989jackknife}. At each point-flagged recipient--date, 4,095 bootstrap samples draw blocks of 5, 10 or 20 original positions. Every sample re-estimates the residual median, upper order statistics and Hill bias correction. Nonfinite observations remain attached to their resampled calendar positions and are excluded when estimating the tail index; invalid bootstrap estimates receive an infinite upper bound. The bounds are conditional on the fitted residual series.

Let $U_{jt}(p)$ denote the empirical $p$-quantile of the bootstrap tail-index estimates, rounded upward to an observed order statistic. The pointwise classification requires $\widehat\nu_{jt}<2$ and $U_{jt}(.95)<2$. The network-adjusted classification replaces $.95$ by $1-.05/K$, allocating the nominal 5\% error budget across nodes within each date. The approximate percentile bounds require local tail regularity and short-range dependence.
\begin{table}[!htbp]\centering\small\renewcommand{\arraystretch}{1.12}\setlength{\tabcolsep}{5pt}
\caption{Block-bootstrap upper-bound classifications of tail signals.}\label{tbl-network-tail-bounds}
\begin{tabular}{lrrrrr}\toprule
 & Point-flagged & \multicolumn{3}{c}{Smallest 95\% upper bound} & Retained\\\cmidrule(lr){3-5}Market & node--dates & Block 5 & Block 10 & Block 20 & node--dates\\\midrule
Stocks & 514 & 3.15 & 2.88 & 2.52 & 0\\
Currencies & 20 & 3.91 & 3.60 & 3.28 & 0\\
Futures & 821 & 2.67 & 2.51 & 2.37 & 0\\\bottomrule\end{tabular}
\par\smallskip\begin{minipage}{.98\linewidth}\footnotesize\raggedright\textit{Notes.} Among primary-specification estimates below two, the middle columns give the smallest one-sided percentile upper bound across recipient--dates, separately by block length. A flag is retained when both the estimate and its upper bound are below two. The retained count is common to all block lengths and to pointwise and within-date Bonferroni-adjusted classifications across the $K$ nodes. Bounds use 4,095 resamples of each 120-position window, retaining calendar gaps.\end{minipage}
\end{table}

\begin{table}[!htbp]\centering\small\renewcommand{\arraystretch}{1.05}\setlength{\tabcolsep}{4pt}
\caption{Innovation-tail diagnostics around financial events.}\label{tbl-primary-tail-events}
\begin{tabular}{llrrr}\toprule
& & \multicolumn{2}{c}{Mean daily node count} & Distinct nodes\\\cmidrule(lr){3-4} Market & Event & $\hat\nu<2$ & $\hat\nu<4$ & Ever $\hat\nu<2$\\\midrule

Stocks & 2018 volatility shock & 0.20 & 28.32 & 3\\
Stocks & COVID & 0.03 & 25.08 & 3\\
Stocks & Ukraine invasion (2022) & 0.19 & 19.21 & 1\\
Stocks & SVB & 0.28 & 28.18 & 2\\
\addlinespace[3pt]
Currencies & China FX & 0.00 & 2.54 & 0\\
Currencies & Brexit & 0.16 & 4.29 & 1\\
Currencies & 2018 volatility shock & 0.00 & 3.52 & 0\\
Currencies & COVID & 0.00 & 3.73 & 0\\
Currencies & Ukraine invasion (2022) & 0.00 & 4.44 & 0\\
Currencies & SVB & 0.00 & 3.71 & 0\\
\addlinespace[3pt]
Futures & China FX & 0.03 & 3.01 & 1\\
Futures & Brexit & 0.02 & 3.16 & 1\\
Futures & 2018 volatility shock & 0.03 & 2.85 & 1\\
Futures & COVID & 0.00 & 2.83 & 0\\
Futures & Ukraine invasion (2022) & 0.32 & 2.62 & 1\\
Futures & SVB & 0.28 & 3.18 & 1\\
\bottomrule\end{tabular}
\par\smallskip\begin{minipage}{.98\linewidth}\footnotesize\raggedright\textit{Notes.} Neighborhoods contain 60 original positions before the anchor and 60 beginning at it. Counts use trailing $W=120$ bias-corrected Hill estimates with a 10\% upper-tail fraction. Mean columns average daily asset counts; the last column counts distinct assets.\end{minipage}
\end{table}

\begin{table}[!htbp]\centering\small\renewcommand{\arraystretch}{1.05}\setlength{\tabcolsep}{4pt}
\caption{Tail-fraction and Hill-correction sensitivity.}\label{tbl-primary-tail-fractions}
\begin{tabular}{lrrrrr}\toprule
Market & Tail fraction & Hill & Corrected & $<2$ (\%) & $<4$ (\%)\\\midrule

Stocks & 5\% & 3.58 & 3.88 & 2.8 & 53.3\\
Stocks & 10\% & 3.15 & 3.59 & 0.6 & 65.1\\
Stocks & 15\% & 2.86 & 3.41 & 0.3 & 76.1\\
Stocks & 25\% & 2.41 & 3.14 & 0.2 & 91.2\\
\addlinespace[3pt]
Currencies & 5\% & 3.79 & 4.10 & 1.3 & 47.5\\
Currencies & 10\% & 3.19 & 3.63 & 0.1 & 64.3\\
Currencies & 15\% & 2.81 & 3.35 & 0.6 & 76.5\\
Currencies & 25\% & 2.34 & 3.04 & 0.1 & 94.5\\
\addlinespace[3pt]
Futures & 5\% & 3.76 & 4.03 & 6.2 & 49.4\\
Futures & 10\% & 3.09 & 3.51 & 4.3 & 63.8\\
Futures & 15\% & 2.79 & 3.32 & 3.3 & 75.7\\
Futures & 25\% & 2.36 & 3.05 & 2.8 & 91.0\\
\bottomrule\end{tabular}
\par\smallskip\begin{minipage}{.98\linewidth}\footnotesize\raggedright\textit{Notes.} Windows contain 120 original positions. Hill columns give medians across valid asset--dates. The final columns give the percentages of valid bias-corrected asset--date estimates below the moment thresholds.\end{minipage}
\end{table}

\begin{figure}[!htbp]\centering
\includegraphics[width=\linewidth]{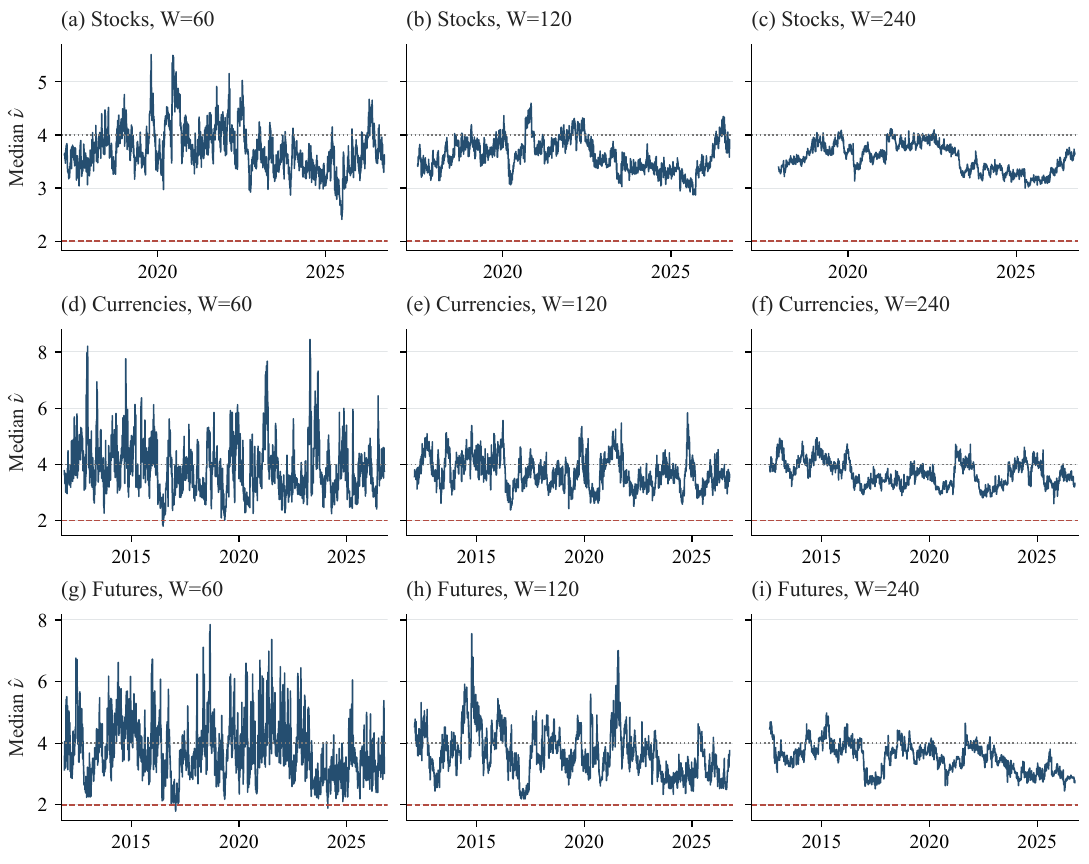}
\caption{\textbf{Trailing-window sensitivity of fitted innovation tails.} Columns use $W=60,120,240$ original positions; rows identify markets. Curves are cross-asset median bias-corrected Hill estimates with a 10\% upper-tail fraction. Reference lines mark two and four.}\label{fig-primary-tail-windows}
\end{figure}

\FloatBarrier
\subsection{Response profiles behind recipient selection}\label{app-stock-monitoring}
\begin{figure}[!htbp]\centering
\includegraphics[width=\linewidth]{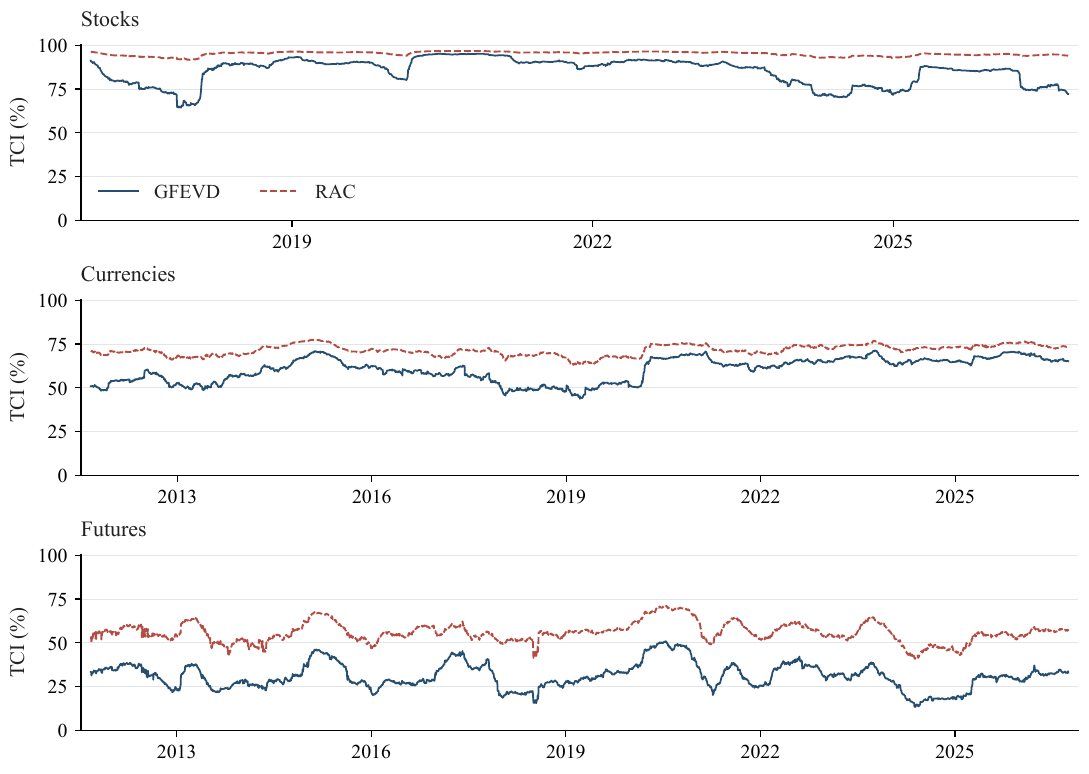}
\caption{\textbf{Daily total connectedness under GFEVD and RAC.} Each panel reports total connectedness in percent for one market over all evaluation dates. Both allocations use the same fitted dynamics.}\label{fig-primary-daily}
\end{figure}

\begin{figure}[!htbp]\centering
\includegraphics[width=\linewidth]{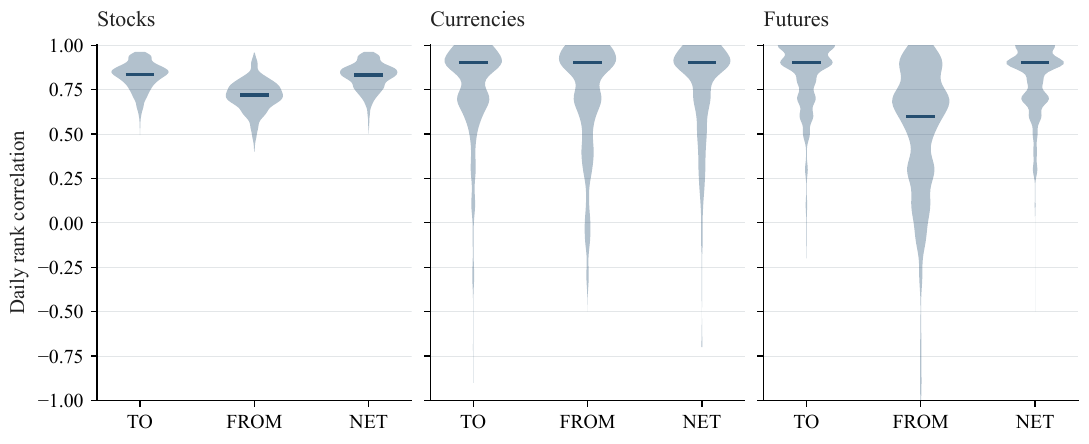}
\caption{\textbf{Daily directional-rank agreement between GFEVD and RAC.} Each Spearman correlation uses every node; distributions cover all evaluation dates and lines mark medians.}\label{fig-primary-ranks}
\end{figure}

Across stocks, currencies and futures, recipients selected only by amplitude have shorter own-response durations and higher external/own duration ratios than recipients selected only by the comparison rule (Table~\ref{tbl-primary-receiver-profiles}). This ordering holds against both GFEVD and the common-rank-scatter squared allocation. Figures~\ref{fig-primary-daily} and \ref{fig-primary-ranks} show daily total connectedness and directional rank agreement.

\begin{table}[!htbp]\centering\small\renewcommand{\arraystretch}{1.05}\setlength{\tabcolsep}{4pt}
\caption{Response duration among selected recipients.}\label{tbl-primary-receiver-profiles}
\begin{tabular}{lllrrr}\toprule
& & & \multicolumn{2}{c}{Effective duration} & External/own\\\cmidrule(lr){4-5} Market & Comparison & Membership & Own & External & ratio\\\midrule

Currencies & GFEVD & Amplitude only & 3.02 & 4.04 & 1.42\\
Currencies & GFEVD & Common & 3.48 & 4.21 & 1.22\\
Currencies & GFEVD & Comparison only & 3.77 & 3.96 & 1.05\\
\addlinespace[3pt]
Currencies & Rank squared & Amplitude only & 2.89 & 3.88 & 1.44\\
Currencies & Rank squared & Common & 3.47 & 4.23 & 1.23\\
Currencies & Rank squared & Comparison only & 3.77 & 4.02 & 1.06\\
\addlinespace[3pt]
Futures & GFEVD & Amplitude only & 3.39 & 6.78 & 2.14\\
Futures & GFEVD & Common & 4.13 & 6.32 & 1.67\\
Futures & GFEVD & Comparison only & 5.82 & 6.30 & 1.09\\
\addlinespace[3pt]
Futures & Rank squared & Amplitude only & 3.23 & 6.62 & 2.20\\
Futures & Rank squared & Common & 4.12 & 6.42 & 1.69\\
Futures & Rank squared & Comparison only & 6.18 & 6.96 & 1.14\\
\addlinespace[3pt]
Stocks & GFEVD & Amplitude only & 3.80 & 6.85 & 1.89\\
Stocks & GFEVD & Common & 4.26 & 6.84 & 1.68\\
Stocks & GFEVD & Comparison only & 4.69 & 6.80 & 1.49\\
\addlinespace[3pt]
Stocks & Rank squared & Amplitude only & 3.69 & 6.86 & 1.95\\
Stocks & Rank squared & Common & 4.25 & 6.84 & 1.68\\
Stocks & Rank squared & Comparison only & 5.02 & 7.11 & 1.47\\
\bottomrule\end{tabular}
\par\smallskip\begin{minipage}{.98\linewidth}\footnotesize\raggedright\textit{Notes.} Daily selection budgets are ten stocks and two assets in each smaller panel. All durations are calculated with the fitted rank scatter. Membership groups are defined relative to RAC; the last column averages recipient--date external/own duration ratios.\end{minipage}
\end{table}

\FloatBarrier
\subsection{External evidence from global banks}\label{app-external-banks}
The public data accompanying \citet{demirer2018estimating} contain daily range-based variance for 96 bank stocks and ten government bonds. We use all 96 bank series over the 2,676 dates from September 12, 2003 to February 7, 2014. The bank panel is complete and strictly positive. Its fixed membership follows the original study's publicly traded subset of the world's largest banks by assets.

We model log bank variance and retain the stock estimator settings: $m=250$, $H=10$, $\gamma=10$, calibration on the first 500 positions and up to 1,500 historical observations preceding the current window.

Evaluation uses each month's last available date after position 500, from August 31, 2005 to February 7, 2014, giving 103 endpoints. Each local fit contains 250 adjacent lag--response pairs. The original study estimates rolling elastic-net VARs; here all four allocations use the common fitted dynamics specified in Section~\ref{app-estimation}.

Table~\ref{tbl-bank-priorities} compares a fixed ten-bank list with a Top 24 list covering one quarter of the panel. At the latter budget, the common-rank-scatter comparison replaces 7.22 names on average. Applying the squared-selected list to the amplitude objective forgoes 0.165 percentage points of mean external share; the reverse comparison forgoes 1.297 points under the squared objective. As in the stock application, transmitter ranks agree more closely than recipient ranks.
\begin{table}[!htbp]\centering\small\renewcommand{\arraystretch}{1.1}\setlength{\tabcolsep}{5pt}
\caption{Global-bank recipient priorities under common fitted dynamics.}\label{tbl-bank-priorities}
\begin{tabular}{lrrrrr}\toprule
Budget & Replaced [95\% CI] & From corr. & To corr. & $L_1$ (pp) & $L_2$ (pp)\\\midrule
\multicolumn{6}{l}{\textit{GFEVD versus RAC}}\\
Top 24 & 8.10 [7.39,8.81] & 0.812 & 0.909 & 0.224 & 1.554\\
Top 10 & 5.40 [4.95,5.85] & 0.812 & 0.909 & 0.292 & 1.874\\
\multicolumn{6}{l}{\textit{Common rank scatter: squared versus amplitude}}\\
Top 24 & 7.22 [6.29,8.16] & 0.864 & 0.957 & 0.165 & 1.297\\
Top 10 & 4.62 [4.05,5.19] & 0.864 & 0.957 & 0.209 & 1.592\\
\multicolumn{6}{l}{\textit{Common covariance: squared versus amplitude}}\\
Top 24 & 8.38 [7.31,9.45] & 0.820 & 0.952 & 0.184 & 1.496\\
Top 10 & 5.32 [4.83,5.81] & 0.820 & 0.952 & 0.254 & 1.617\\
\bottomrule\end{tabular}
\par\smallskip\begin{minipage}{.98\linewidth}\footnotesize\raggedright\textit{Notes.} Averages use all 103 monthly endpoints and all 96 banks. Replaced is the number of names in either list but not the other, counted on one side. Correlations are Spearman coefficients. $L_1$ is the mean external amplitude share forgone by using the squared-selected list; $L_2$ reverses the comparison under the squared objective. Both are percentage-point differences in fitted objective values. The first two panels use the RAC amplitude objective; the third uses covariance-amplitude allocation. Brackets give Bartlett HAC(12) normal-approximation intervals for the time-average number replaced.\end{minipage}
\end{table}

\FloatBarrier
Source-dominance reversals occur on 1,436 of 9,888 bank--month observations (14.5\%) under common rank scatter, with a time-average HAC interval of [2.2\%,26.9\%]. Requiring squared external share below 45\% and amplitude external share above 55\% retains 10.4\% of all observations. The reverse classification does not occur. Table~\ref{tbl-bank-years} locates the reversals in time. Their frequency falls from 32.8\% at impact to 14.5\% over ten horizons. In this panel, accumulation reduces the contemporaneous classification gap. Common-scatter Top 24 lists differ in every annual block, with the final block covering January and February 2014.
\begin{table}[!htbp]\centering\small\renewcommand{\arraystretch}{1.1}\setlength{\tabcolsep}{7pt}
\caption{Global-bank list changes and source dominance over time.}\label{tbl-bank-years}
\begin{tabular}{lrrrrr}\toprule
& & & \multicolumn{3}{c}{Own-to-external reversals (\%)}\\\cmidrule(lr){4-6} Year & Endpoints & Replaced & $H=1$ & $H=10$ & Buffered\\\midrule
2005 & 5 & 7.60 & 95.0 & 79.2 & 66.0\\
2006 & 12 & 6.92 & 83.2 & 45.4 & 32.9\\
2007 & 12 & 7.25 & 54.9 & 18.8 & 11.5\\
2008 & 12 & 10.00 & 12.6 & 1.7 & 0.9\\
2009 & 12 & 8.33 & 9.8 & 0.8 & 0.4\\
2010 & 12 & 7.33 & 16.0 & 1.9 & 1.0\\
2011 & 12 & 6.42 & 16.3 & 6.2 & 3.6\\
2012 & 12 & 5.83 & 20.2 & 4.1 & 2.4\\
2013 & 12 & 6.33 & 24.9 & 10.8 & 7.6\\
2014 & 2 & 2.50 & 24.5 & 12.0 & 6.2\\
\midrule
All & 103 & 7.22 & 32.8 & 14.5 & 10.4\\
\bottomrule\end{tabular}
\par\smallskip\begin{minipage}{.98\linewidth}\footnotesize\raggedright\textit{Notes.} Dynamics and rank scatter are common. Replaced is the mean number of Top 24 recipients selected by one allocation but not the other. The last three columns are percentages of bank--month observations. The $H=1$ and $H=10$ columns report $F_1>1/2>F_2$ at impact and over ten horizons. Buffered requires $F_1>0.55$ and $F_2<0.45$ at $H=10$. The first year begins in August; the final year ends at the February 7 data cutoff.\end{minipage}
\end{table}

\FloatBarrier

\subsection{Source dominance at impact and over time}\label{app-futures-source-case}
At impact, positive-semidefinite scatter implies that absolute aggregation gives an external share at least as large as squared aggregation (Section~\ref{app-impact-inequality}). Table~\ref{tbl-futures-mechanism} separates this contemporaneous comparison from ten-horizon reversals and classifies the latter using source breadth and relative duration. For corn, the reversal frequency increases from 5.1\% at impact to 68.8\% over ten horizons; for E-mini futures it is 74.1\% and 73.5\%, respectively.
\begin{table}[!htbp]\centering\small\renewcommand{\arraystretch}{1.05}\setlength{\tabcolsep}{4pt}
\caption{Futures reversals at impact and across the response path.}\label{tbl-futures-mechanism}
\begin{tabular}{lrrr*{4}{>{\raggedleft\arraybackslash}p{.085\linewidth}}}\toprule
& \multicolumn{2}{c}{Reversals (\%)} & & \multicolumn{4}{c}{Mechanism among $H=10$ reversals (\%)}\\\cmidrule(lr){2-3}\cmidrule(lr){5-8} Contract & $H=1$ & $H=10$ & New (\%) & Breadth & Duration & Either & Joint\\\midrule
All & 41.6 & 70.9 & 51.8 & 59.7 & 0.5 & 22.5 & 17.2\\
C & 5.1 & 68.8 & 97.7 & 26.5 & 2.6 & 49.6 & 21.2\\
CL & 61.9 & 73.4 & 37.1 & 75.6 & 0.0 & 10.5 & 13.9\\
ES & 74.1 & 73.5 & 11.2 & 95.6 & 0.0 & 2.1 & 2.4\\
GC & 62.7 & 84.8 & 33.3 & 57.8 & 0.0 & 9.9 & 32.2\\
NG & 4.3 & 54.1 & 97.6 & 34.7 & 0.0 & 51.9 & 13.4\\
\bottomrule\end{tabular}
\par\smallskip\begin{minipage}{.98\linewidth}\footnotesize\raggedright\textit{Notes.} All entries are percentages. The first two columns use every recipient--date (19,500 in All; 3,900 per contract). New is the percentage of $H=10$ reversals absent at impact. The final four columns partition $H=10$ reversals. Holding measured $r=A_E/A_O$ fixed, breadth is sufficient if $r^2<n_E$, and duration if $r^2<d_E/d_O$: Breadth and Duration mean exactly one is sufficient; Either means each is sufficient separately; Joint means neither alone is sufficient. The comparisons reallocate the same responses with fitted dynamics held fixed.\end{minipage}
\end{table}

\FloatBarrier

\begin{figure}[!htbp]\centering
\includegraphics[width=\linewidth]{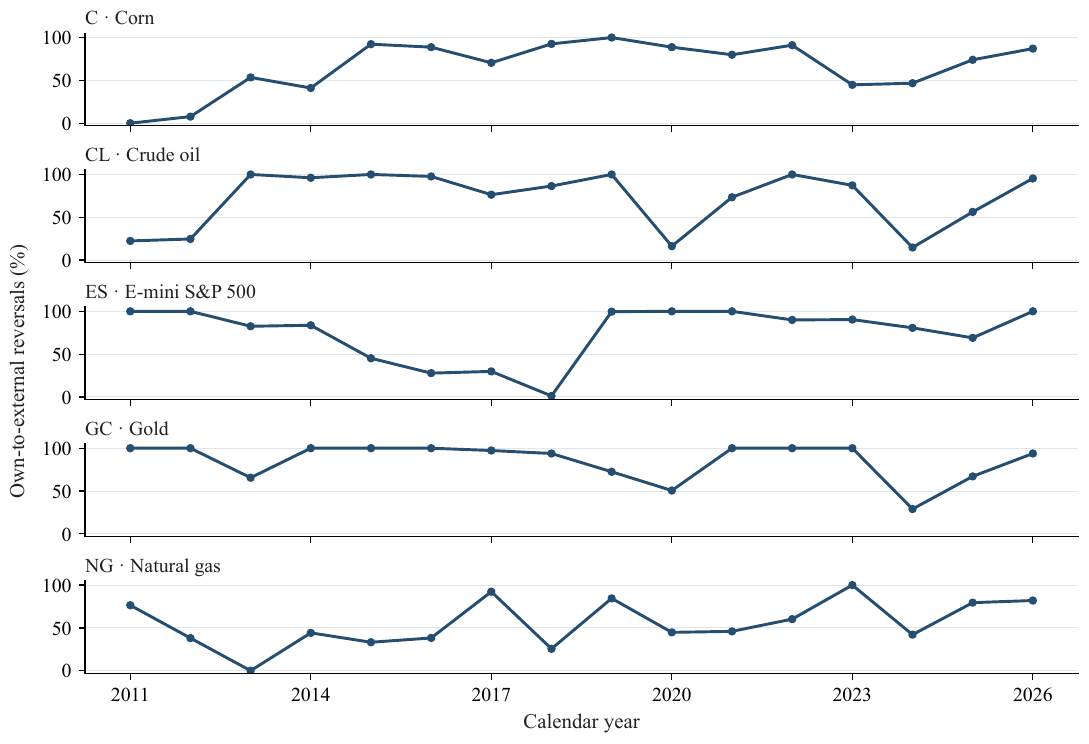}
\caption{\textbf{Annual futures source-dominance reversals.} Each panel gives one contract over every calendar year, with the percentage of dates with own-dominated squared responses and externally dominated amplitude. Dynamics and rank scatter are held fixed.}\label{fig-primary-classification}
\end{figure}

\FloatBarrier
\subsection{Scatter, aggregation and response horizon}\label{app-network-sensitivity}
Changing response aggregation contributes more to recipient-strength and recipient-rank differences than changing scatter, in every market and under both endpoint references (Table~\ref{tbl-primary-decomposition}). For the identity of stock core links, scatter contributes more.

For a GFEVD reference, let $d_A$ measure the distance after changing aggregation alone, $d_S$ the distance after changing scatter alone, and $d_{SA}$ the distance after both changes. The contributions are
\[
\begin{aligned}
\phi_A&=\tfrac12(d_A+d_{SA}-d_S),\\
\phi_S&=\tfrac12(d_S+d_{SA}-d_A).
\end{aligned}
\]
The two contributions sum to the endpoint distance. We sum them over dates before normalizing; the reverse comparison starts from RAC and evaluates the same two changes.

Extending the response horizon lowers incoming-rank agreement in all three panels while the fitted systems remain fixed (Table~\ref{tbl-primary-horizons}).
\par
\begin{table}[!htbp]\centering\small\renewcommand{\arraystretch}{1.05}\setlength{\tabcolsep}{4pt}
\caption{Scatter and aggregation under both endpoint references.}\label{tbl-primary-decomposition}
\begin{tabular}{lllrr}\toprule
Market & Quantity & Reference & Scatter (\%) & Aggregation (\%)\\\midrule

Currencies & Core identity & GFEVD & 48.2 & 51.8\\
Currencies & Core identity & RAC & 45.9 & 54.1\\
Currencies & Incoming rank & GFEVD & 39.0 & 61.0\\
Currencies & Incoming rank & RAC & 34.9 & 65.1\\
Currencies & Incoming strength & GFEVD & 10.3 & 89.7\\
Currencies & Incoming strength & RAC & 5.9 & 94.1\\
\addlinespace[3pt]
Futures & Core identity & GFEVD & 46.6 & 53.4\\
Futures & Core identity & RAC & 45.4 & 54.6\\
Futures & Incoming rank & GFEVD & 45.3 & 54.7\\
Futures & Incoming rank & RAC & 38.3 & 61.7\\
Futures & Incoming strength & GFEVD & 7.6 & 92.4\\
Futures & Incoming strength & RAC & 0.1 & 99.9\\
\addlinespace[3pt]
Stocks & Core identity & GFEVD & 69.4 & 30.6\\
Stocks & Core identity & RAC & 66.7 & 33.3\\
Stocks & Incoming rank & GFEVD & 33.9 & 66.1\\
Stocks & Incoming rank & RAC & 36.8 & 63.2\\
Stocks & Incoming strength & GFEVD & 10.0 & 90.0\\
Stocks & Incoming strength & RAC & 2.7 & 97.3\\
\bottomrule\end{tabular}
\par\smallskip\begin{minipage}{.98\linewidth}\footnotesize\raggedright\textit{Notes.} Shapley contributions use the same fitted dynamics and the four scatter--aggregation allocations. Contributions are summed over dates before normalization.\end{minipage}
\end{table}

\par

\par
\begin{table}[!htbp]\centering\small\renewcommand{\arraystretch}{1.05}\setlength{\tabcolsep}{4pt}
\caption{Response-horizon sensitivity on the same fitted systems.}\label{tbl-primary-horizons}
\begin{tabular}{lrrrr}\toprule
Market & $H$ & TO corr. & FROM corr. & Core replaced (\%)\\\midrule

Currencies & 5 & 0.802 & 0.734 & 34.1\\
Currencies & 20 & 0.726 & 0.703 & 37.2\\
\addlinespace[3pt]
Futures & 5 & 0.838 & 0.639 & 38.2\\
Futures & 20 & 0.833 & 0.492 & 38.8\\
\addlinespace[3pt]
Stocks & 5 & 0.840 & 0.760 & 39.6\\
Stocks & 20 & 0.816 & 0.684 & 42.6\\
\bottomrule\end{tabular}
\par\smallskip\begin{minipage}{.98\linewidth}\footnotesize\raggedright\textit{Notes.} Correlations and core replacement compare GFEVD with RAC over all original dates. Core replacement uses a fixed 10\% external-channel density. The fitted coefficient and scatter estimates are held fixed as the response horizon changes.\end{minipage}
\end{table}

\par

\par
\begin{samepage}
The estimation sensitivity grid varies the current window and the amount of pre-window history. All five specifications use common evaluation dates starting at position 1,001, ensuring at least 500 historical observations before the longest current window. The historical block always precedes the current pairs.
\par
\begin{table}[!htbp]\centering\small\renewcommand{\arraystretch}{1.05}\setlength{\tabcolsep}{4pt}
\caption{Estimation-window and historical-span sensitivity.}\label{tbl-single-primary-sensitivity}
\begin{tabular}{llrrrr}\toprule
Market & Specification & Replaced & Rank move & Reversal (\%) & Margin 5 (\%)\\\midrule

Stocks & Primary & 3.25 & 5.09 & 0.0 & 0.0\\
Stocks & $m=125$ & 2.86 & 4.71 & 0.8 & 0.3\\
Stocks & $m=500$ & 4.24 & 6.13 & 0.0 & 0.0\\
Stocks & History 500 & 3.17 & 4.92 & 0.0 & 0.0\\
Stocks & History 2,000 & 3.51 & 5.52 & 0.0 & 0.0\\
\addlinespace[3pt]
Currencies & Primary & 0.37 & 0.47 & 9.0 & 3.1\\
Currencies & $m=125$ & 0.33 & 0.50 & 12.7 & 5.7\\
Currencies & $m=500$ & 0.28 & 0.40 & 7.7 & 3.1\\
Currencies & History 500 & 0.33 & 0.47 & 9.2 & 2.9\\
Currencies & History 2,000 & 0.36 & 0.47 & 8.9 & 2.9\\
\addlinespace[3pt]
Futures & Primary & 0.57 & 0.76 & 72.8 & 46.3\\
Futures & $m=125$ & 0.61 & 0.83 & 72.9 & 51.0\\
Futures & $m=500$ & 0.55 & 0.69 & 70.6 & 44.2\\
Futures & History 500 & 0.59 & 0.78 & 71.7 & 45.0\\
Futures & History 2,000 & 0.57 & 0.76 & 73.3 & 46.8\\
\bottomrule\end{tabular}
\par\smallskip\begin{minipage}{.98\linewidth}\footnotesize\raggedright\textit{Notes.} All specifications use common dates within each panel: 1,952 stock dates, 3,414 currency dates and 3,400 futures dates, beginning at position 1,001. Primary uses $m=250$ and up to 1,500 historical observations. Alternative windows retain that history cap; alternative histories retain $m=250$. Replaced counts daily Top 10 stock or Top 2 smaller-panel names changed by aggregation alone. Reversal requires $F_1>1/2>F_2$; Margin 5 requires $F_1>0.55$ and $F_2<0.45$.\end{minipage}
\end{table}

\par

\par\end{samepage}

\FloatBarrier
\paragraph*{Lag order.}
We add a second lag while retaining the same estimator, initial-500 calibration and penalty multiplier, $m=250$, the 1,500-position historical cap and $H=10$. All original evaluation dates are retained. A usable response requires both lags on the original calendar. Historical templates retain both own lags and add lag-specific aggregates, rank-factor directions or sparse external lag coordinates. In the historical selection rule, $2K$ replaces $K$ in the factor threshold, and $2(K-1)$ replaces $K-1$ in the sparse sizes and search cost. The bounded coefficient adjustment and residual-scatter construction are unchanged.

Table~\ref{tbl-current-lag-sensitivity} reports the matched-date comparisons. With two lags, attribution replaces 3.01 stock Top 10 names and shifts ranks by 4.83 places on average. Amplitude-only stock recipients have a mean external/own duration ratio of 1.94, compared with 1.44 for energy-only recipients. Futures external share is 27.5\% under energy and 54.7\% under amplitude, with 68.3\% own-to-external reversals.

Of the 10,266 VAR(2) fits, 19 have companion spectral radius at least one; the maximum is 1.025106. The finite-horizon allocations retain these dates.

\begin{table}[!htbp]\centering\small\renewcommand{\arraystretch}{1.05}\setlength{\tabcolsep}{4pt}
\caption{Attribution under alternative lag orders.}\label{tbl-current-lag-sensitivity}
\begin{tabular}{llrrrrrr}\toprule
& & & & \multicolumn{2}{c}{External share (\%)} & & \\
Market & Lags & Replaced & Rank move & Energy & Amplitude & Reversal (\%) & Margin 5 (\%)\\\midrule

Stocks & 1 & 3.14 & 5.05 & 85.4 & 95.1 & 0.0 & 0.0\\
Stocks & 2 & 3.01 & 4.83 & 85.8 & 95.1 & 0.0 & 0.0\\
\addlinespace[3pt]
Currencies & 1 & 0.34 & 0.44 & 60.6 & 71.5 & 10.6 & 4.7\\
Currencies & 2 & 0.33 & 0.43 & 61.1 & 71.4 & 9.2 & 3.1\\
\addlinespace[3pt]
Futures & 1 & 0.54 & 0.71 & 28.0 & 55.9 & 70.9 & 46.4\\
Futures & 2 & 0.53 & 0.65 & 27.5 & 54.7 & 68.3 & 44.9\\
\bottomrule\end{tabular}
\par\smallskip\begin{minipage}{.98\linewidth}\footnotesize\raggedright\textit{Notes.} Both lag orders use all 2,452 stock, 3,914 currency and 3,900 futures evaluation dates, $m=250$ and $H=10$. Within each fit, dynamics and rank scatter are fixed; energy and amplitude differ only in response aggregation. Replaced counts changed Top 10 stock or Top 2 smaller-panel names. Rank move is the mean absolute recipient-rank difference. Reversal requires $F_1>1/2>F_2$; Margin 5 requires $F_1>0.55$ and $F_2<0.45$. Shares and classifications are computed before averaging.\end{minipage}
\end{table}

\FloatBarrier
\subsection{Time-average uncertainty}\label{app-headline-hac}
Let $Z_t$ denote a daily comparison computed from the fitted networks, with stable mean $\mu$. Under a functional central limit theorem and consistent long-run variance estimation, Bartlett HAC estimates the asymptotic variance of $\sqrt N(\bar Z-\mu)$. Table~\ref{tbl-primary-hac} varies the truncation lag for the same sequence $Z_t$ of daily overlaps.
\begin{table}[!htbp]\centering\small\renewcommand{\arraystretch}{1.05}\setlength{\tabcolsep}{4pt}
\caption{Dependence-adjusted time averages of recipient overlap.}\label{tbl-primary-hac}
\begin{tabular}{lrrrr}\toprule
& & & \multicolumn{2}{c}{95\% interval}\\\cmidrule(lr){4-5} Market & HAC lag & Mean (\%) & Lower & Upper\\\midrule

Currencies & 125 & 83.2 & 78.4 & 88.1\\
Currencies & 250 & 83.2 & 77.4 & 89.0\\
Currencies & 500 & 83.2 & 76.8 & 89.6\\
\addlinespace[3pt]
Futures & 125 & 73.2 & 67.9 & 78.5\\
Futures & 250 & 73.2 & 67.6 & 78.8\\
Futures & 500 & 73.2 & 67.3 & 79.1\\
\addlinespace[3pt]
Stocks & 125 & 68.6 & 65.6 & 71.5\\
Stocks & 250 & 68.6 & 65.3 & 71.8\\
Stocks & 500 & 68.6 & 65.6 & 71.5\\
\bottomrule\end{tabular}
\par\smallskip\begin{minipage}{.98\linewidth}\footnotesize\raggedright\textit{Notes.} Nominal 95\% intervals concern time-average overlap between amplitude and squared allocations under common rank scatter. Budgets are ten stocks and two assets in the smaller panels. Daily fitted comparisons use the original calendar sequence.\end{minipage}
\end{table}

\FloatBarrier
\section{External evidence from industrial production}\label{app-external-industry}
The public replication panel of \citet{brunner2023estimating} extends the comparison to 88 industrial sectors over January 1972--December 2019. Its complete, positive, seasonally adjusted production indices include historical backfilling from higher-level series. The authors' transformation, $y_{it}=1200\Delta\log(IP_{it})$, gives 575 monthly growth observations without further imputation.

Estimation uses $m=250$, $\gamma=10$, initial-500 coordinate calibration and up to 1,500 pre-window observations. The 74 monthly endpoints from November 2013 to December 2019 each have 250 adjacent pairs. All fitted systems are stable. Responses use the original three-month horizon, with $H=10$ as a sensitivity check. Our unconditional growth dynamics replace the original study's factor-augmented elastic-net systems; the links allocate growth responses across sources and horizons.

Table~\ref{tbl-industry-priorities} separates scatter from aggregation for Top 22 and Top 10 lists. Top 22 covers one quarter of the panel. At $H=3$, GFEVD and RAC replace 5.27 names on average; holding rank scatter fixed replaces 3.22. The corresponding amplitude and squared objective shortfalls are 0.318 and 0.648 percentage points.

The larger difference concerns the dominant source of production-growth responses. With common rank scatter, 5,965 of 6,512 sector--month observations (91.6\%) are own-dominated under squared responses and externally dominated under amplitude. Requiring external shares below 45\% and above 55\%, respectively, retains 87.0\% of all observations. The reverse classification occurs on 0.0\%. At ten months the own-to-external frequency is 91.0\%. Table~\ref{tbl-industry-years} reports every annual block at both horizons.

At $H=3$, external amplitude spans 67.63 effective source channels on average; own and external responses occupy 1.53 and 1.94 effective horizons. Broad contributions across sectors accompany the change in dominant-source diagnosis. Allowing scatter to change as well, the full GFEVD--RAC comparison reverses source dominance on 43.8\% of observations.
\begin{table}[!htbp]\centering\small\renewcommand{\arraystretch}{1.08}\setlength{\tabcolsep}{5pt}
\caption{Industrial recipient priorities under common fitted dynamics.}\label{tbl-industry-priorities}
\begin{tabular}{lrrrrr}\toprule
$H$ & Budget & Replaced [range] & From corr. & $L_1$ (pp) & $L_2$ (pp)\\\midrule
\multicolumn{6}{l}{\textit{GFEVD versus RAC}}\\
3 & 22 & 5.27 [4,8] & 0.785 & 0.742 & 2.212\\
3 & 10 & 3.77 [2,6] & 0.785 & 0.844 & 2.397\\
10 & 22 & 5.73 [4,7] & 0.749 & 0.823 & 2.333\\
10 & 10 & 3.27 [2,6] & 0.749 & 0.855 & 2.012\\
\multicolumn{6}{l}{\textit{Common rank scatter: squared versus amplitude}}\\
3 & 22 & 3.22 [2,4] & 0.920 & 0.318 & 0.648\\
3 & 10 & 1.72 [0,3] & 0.920 & 0.148 & 0.825\\
10 & 22 & 3.39 [2,5] & 0.877 & 0.360 & 0.857\\
10 & 10 & 1.36 [0,3] & 0.877 & 0.197 & 0.973\\
\multicolumn{6}{l}{\textit{Common covariance: squared versus amplitude}}\\
3 & 22 & 4.65 [3,6] & 0.883 & 0.308 & 1.776\\
3 & 10 & 2.81 [1,4] & 0.883 & 0.186 & 1.151\\
10 & 22 & 5.16 [3,7] & 0.833 & 0.442 & 1.895\\
10 & 10 & 2.84 [2,4] & 0.833 & 0.283 & 1.463\\
\bottomrule\end{tabular}
\par\smallskip\begin{minipage}{.98\linewidth}\footnotesize\raggedright\textit{Notes.} All 88 sectors and 74 monthly endpoints enter. Replaced counts names in either list but not the other, on one side; brackets give the minimum and maximum monthly counts. From corr. is the mean Spearman correlation of recipient ranks. $L_1$ is the mean external amplitude share forgone by using the squared-selected list; $L_2$ reverses the comparison under the squared objective. Both are fitted-objective differences in percentage points. The third panel uses covariance-amplitude allocation. $H$ counts monthly responses, including impact.\end{minipage}
\end{table}

\begin{table}[!htbp]\centering\small\renewcommand{\arraystretch}{1.1}\setlength{\tabcolsep}{7pt}
\caption{Industrial source classification across years and horizons.}\label{tbl-industry-years}
\begin{tabular}{lrrrrrr}\toprule
 & & \multicolumn{3}{c}{$H=3$} & \multicolumn{2}{c}{$H=10$}\\
\cmidrule(lr){3-5}\cmidrule(lr){6-7}
Year & Months & Replaced & Reversal (\%) & Buffered (\%) & Replaced & Reversal (\%)\\\midrule
2013 & 2 & 2.50 & 85.2 & 82.4 & 3.00 & 84.1\\
2014 & 12 & 3.08 & 87.4 & 82.2 & 2.67 & 86.3\\
2015 & 12 & 2.42 & 89.5 & 83.0 & 3.33 & 88.7\\
2016 & 12 & 3.67 & 91.8 & 86.3 & 3.92 & 90.9\\
2017 & 12 & 3.50 & 92.9 & 88.3 & 3.58 & 92.6\\
2018 & 12 & 3.17 & 93.8 & 90.2 & 3.08 & 93.6\\
2019 & 12 & 3.58 & 95.3 & 92.9 & 3.83 & 94.8\\
\midrule
All & 74 & 3.22 & 91.6 & 87.0 & 3.39 & 91.0\\
\bottomrule\end{tabular}
\par\smallskip\begin{minipage}{.98\linewidth}\footnotesize\raggedright\textit{Notes.} Dynamics and rank scatter are common. Replaced is the mean number of Top 22 recipients in one list but not the other. Reversal is the percentage of sector--month observations with $F_1>1/2>F_2$; Buffered requires $F_1>0.55$ and $F_2<0.45$. The 2013 block contains November and December; every later block contains twelve months. All-row percentages pool the 6,512 sector--month observations at each horizon.\end{minipage}
\end{table}

\FloatBarrier

\FloatBarrier

\end{document}